\baselineskip=12pt

\def\adots{\mathinner{\mskip1mu\raise1pt\hbox{.}\mskip2mu\raise4pt\hbox{.}\mskip2mu\raise7pt\vbox{\kern7pt\hbox{.}}\mskip1mu}}

\def\ker{\mathop{\rm ker}\nolimits}
\def\Lie{\mathop{\hbox{\it \$}}\nolimits}
\def\shape{\mathop{\Pi\llap{$\amalg$}}\nolimits}
\def\tr{\mathop{\rm tr}\nolimits}

\mathchardef\bfplus="062B
\mathchardef\bfminus="067B
\font\title=cmbx10 scaled\magstep5
\font\chapter=cmbx10 scaled\magstep4
\font\section=cmbx10 scaled\magstep2

\def\~#1{{\accent"7E #1}}
\def\sqr#1#2{{\vcenter{\hrule height.#2pt \hbox{\vrule width.#2pt height#1pt \kern#1pt \vrule width.#2pt}\hrule height.#2pt}}}
\def\square{\mathchoice\sqr63\sqr63\sqr{4.2}2\sqr{1.5}2}
\def\Square{\mathord{\hskip 1pt\square\hskip 1pt}}
\def\hk#1#2{{\vcenter{\hrule height0.0pt \hbox{\vrule width0.0pt \kern#1pt \vrule width.#2pt height#1pt}\hrule height.#2pt}}}

\def\operp{\mathrel{\bigcirc \hskip -11.5pt \perp}}

\def\1{{\overline 1}}
\def\2{{\overline 2}}

\centerline{\chapter Spin Coefficients for Four-Dimensional}
\centerline{\chapter Neutral Metrics, and Null Geometry}
\vskip 24pt
\noindent {\section Abstract}\hfil\break
Notation for spin coefficients for metrics of neutral signature in four dimensions is introduced. The utility and interpretation of spin coefficients is explored through themes in null geometry familiar from (complex) general relativity. Four-dimensional Walker geometry is exploited to provide examples and the generalization of the real neutral version of Pleba\'nski's second heavenly equation to certain Walker geometries given in [P.R. Law, Y. Matsushita, A Spinor Approach to Walker Geometry, Comm. Math. Phys. 282 (2008) 577--623] is extended further.
\vskip 24pt
\noindent Peter R Law. 4 Mack Place, Monroe, NY 10950, USA. prldb@member.ams.org\hfil\break
\vskip 12pt
\noindent JGP: Real and complex differential geometry; Spinors and twistors
\vskip 12pt
\noindent Key Words and Phrases: neutral geometry, Walker geometry, four dimensions, spinors.
\vfill\eject
\noindent{\section 1. Introduction}
\vskip 12pt
Anticipating the value of spinors in the study of four-dimensional neutral geometry, I presented in [16] the algebraic classification of the Weyl curvature spinors for neutral metrics in four dimensions, complementing the Petrov-like classification of the Weyl curvature tensor already given in [15] (see also [5], \S\S 39 - 40, and [24], \S 18). A spinor approach has proven very natural in describing four-dimensional Walker geometry, see [17]. Much of the two-component spinor formalism expounded in [22] carries over with only obvious modifications to the case of neutral signature. In fact, the spinor formalism appropriate for neutral signature is best regarded as a real version of the spinor formalism for complex space-times and in large part can be formally obtained by complexifying the Lorentzian spinor theory and then regarding all objects as real valued. See [16], however, for a brief exposition of the Clifford algebra background to the two-component spinor formalism for neutral signature and some spinor algebra and [17], Appendix 2, for the spinor form of curvature.

A few elements of the spinor formalism for neutral signature, however, require explicit treatment. One prominent well known example arises from the fact that the Hodge star operator on $\Lambda^2({\bf R}^{2,2})$ is involutary rather than a complex structure as it is on $\Lambda^2({\bf R}^{1,3})$ and thus the spinor representation of duality of two-forms in the neutral case differs appropriately from that in the Lorentz case. Another example arises from the fact that the priming operation of [22], (4.5.17), cannot be applied in the neutral case since the spin spaces are {\bf R}-linear spaces. Thus, a modified priming operation must be introduced in the neutral case with consequent modifications of the notation for the spin coefficients. In this paper, I specify an appropriate priming operation, introduce notation for the spin coefficients, and detail the subsequent equations for future reference. Through familiar themes in (complex) general relativity (GR), I exploit null geometry to explore the utility and interpretation of spin coefficients in four-dimensional neutral geometry. I treat Walker geometry in some detail. In particular, the generalization of the real neutral version of Pleba\'nski's [25] second heavenly equation to certain Walker geometries given in [17] is extended further.

An important notational device in the neutral case is the use of a tilde in place of complex conjugation in the Lorentz case. For each complex conjugate pair of objects in the Lorentzian case, e.g., $\Psi_{ABCD}$ and $\overline\Psi_{A'B'C'D'}$, {\sl denote\/} the neutral analogue of the conjugated object by replacing the bar by a tilde, obtaining for example $\Psi_{ABCD}$ and $\tilde\Psi_{A'B'C'D'}$ respectively. Similarly, for each equation in the Lorentzian case, say $(*)$, one can take the complex conjugate to obtain another equation (which may or may not actually be an independent equation): $(\overline *)$. Instead of taking the complex conjugate of $(*)$ itself, $(\overline *)$ may also be obtained by taking the complex conjugate of each step in the derivation of $(*)$. With the notational convention regarding the use of tilde in the neutral case, it follows that if $(\dag)$ is the neutral analogue of $(*)$ (obtained by modifying appropriately the derivation of $(*)$ for neutral signature), the neutral analogue of $(\overline *)$ will be a tilde version of $(\dag)$ (obtained by replacing complex conjugation in $(\overline *)$ by tildes but also performing the modifications necessary to obtain $(\dag)$ from $(*)$). In this manner, one establishes a formal analogy between the use of the tilde {\sl notation\/} in the neutral case with the {\sl operation\/} of complex conjugation in the Lorentz case. The notational efficiency of this convention is obvious: one can `apply' tilde to equations in the neutral case to obtain valid equations even though tilde is not an operation on spinors per se. Indeed, one should not be tempted to view tilde as an operation with geometric content; tilde is a purely notational convention of utility.

One further advantage of the tilde notation beyond its utility within the neutral signature case is that it makes the notation in the neutral case compatible with that in the Lorentz case under complexification of each. This observation may in fact be taken as a guiding principle for the employment of tilde and also suggests that the priming notation employed in [22] would be better replaced by that employed in the neutral case for consistency under complexification.

Throughout this paper, $(M,g)$ denotes a four-dimensional manifold $M$ equipped with a metric $g$ of neutral signature. I employ the abstract index notation of [22]. The spinor results given subsequently may be taken to apply locally on $(M,g)$ and globally when $(M,g)$ admits a spinor structure, i.e., when $(M,g)$ admits a reduction of the frame bundle to the identity component ${\bf SO^\bfplus(2,2)}$ and the second Stiefel-Whitney class of the tangent bundle vanishes.
\vskip 24pt
\noindent{\section 2. Spin Coefficients}
\vskip 12pt
At least locally, i.e., on suitable open subsets $U$ of $M$, one can construct bundles $S_U$ and $S'_U$ of unprimed and primed spinors respectively. These bundles are both real vector bundles with fibres isomorphic to ${\bf R}^2$ but they are independent of each other. Moreover, $TU \cong S_U \otimes S'_U$. Following the argument of [22] \S 4.4, the Levi-Civita connection $\nabla_a$ of $g$ induces a connection $D$ on $S_U$ (i.e., an operator $D:\Gamma(S_U) \to \Gamma(S_U \otimes T_\bullet U)$ with the usual properties, where $T_\bullet U$ is the cotangent bundle over $U$) which is uniquely characterized by the condition $D\epsilon_{AB} = 0$ (where $D$ is extended in the usual manner to tensor products of $S_U$). Similarly, $\nabla_a$ induces a connection $D'$ on $S'_U$ which is uniquely characterized by the requirement that $D'\epsilon_{A'B'} = 0$. In the Lorentz case, the fact both are induced from $\nabla_a$ and that $S'$ is the complex conjugate of $S$ is manifested by a relationship between the induced connections involving complex conjugation (with $\alpha \in S_U$, $D'(\overline\alpha) = \overline{D\alpha}$); there is no analogue in the neutral case due to the independence of $S$ and $S'$. Nevertheless, as both are induced from $\nabla_a$, both may be denoted by the same symbol $\nabla_{AA'}$, which is also employed interchangeably with $\nabla_a$, and this notation also refers to the induced connections on the tensor products of $S_U$ and $S'_U$.

Let $\{\epsilon_{\bf A}{}^A\} = \{o^A,\iota^A\}$ be a frame field for $S_U$ and $\{\epsilon_{\bf A'}{}^{A'}\} = \{o^{A'},\iota^{A'}\}$ a frame field for $S'_U$ (if one needs to distinguish the latter from the former without indices attached, place tildes over the symbols of the latter) with dual frames $\{\epsilon_A{}^{\bf A}\} = \{-\chi^{-1}\iota_A,\chi^{-1}o_A\}$ and $\{\epsilon_{A'}{}^{\bf A'}\} = \{-\tilde\chi^{-1}\iota_{A'},\tilde\chi^{-1}o_{A'}\}$, respectively, where
$$\epsilon_{AB}\epsilon_{\bf 0}{}^A\epsilon_{\bf 1}{}^B =: \chi \hskip 1.25in \epsilon_{A'B'}\epsilon_{\bf 0'}{}^{A'}\epsilon_{\bf 1'}{}^{B'} =: \tilde\chi.\eqno(2.1)$$
As in [22] \S 4.5, one defines
$$\gamma_{{\bf A}{\bf A}'{\bf C}}{}^{\bf B} := \epsilon_{{\bf A}}{}^A\epsilon_{{\bf A}'}{}^{A'}\epsilon_E{}^{\bf B}\nabla_{AA'}\epsilon_{\bf C}{}^E,\eqno(2.2.{\rm a})$$
but
$$\tilde\gamma_{{\bf A}{\bf A}'{\bf C}'}{}^{{\bf B}'} := \epsilon_{{\bf A}}{}^A\epsilon_{{\bf A}'}{}^{A'}\epsilon_{E'}{}^{{\bf B}'}\nabla_{AA'}\epsilon_{{\bf C}'}{}^{E'}\eqno(2.2.{\rm b})$$
in place of [22] (4.5.3). With $\tilde\gamma_{{\bf A}{\bf A}'{\bf C}'}{}^{{\bf B}'}$ replacing $\overline\gamma_{{\bf A}{\bf A}'{\bf C}'}{}^{{\bf B}'}$, the expression of covariant derivatives of spinor fields in terms of $\gamma_{{\bf A}{\bf A}'{\bf C}}{}^{\bf B}$ and $\tilde\gamma_{{\bf A}{\bf A}'{\bf C}'}{}^{{\bf B}'}$ is as in [22] \S 4.5.

The spin coefficients are just symbols for the various $\gamma_{{\bf A}{\bf A}'{\bf C}}{}^{\bf B}$ and $\tilde\gamma_{{\bf A}{\bf A}'{\bf C}'}{}^{{\bf B}'}$, of which there are 16 each, giving 32 independent real-valued quantities. Following [22], I introduce a {\sl priming operation\/} which makes the notation more efficient.
\vskip 24pt
\noindent {\bf 2.3 The Priming Operation}\hfil\break
The priming operation utilized by [22] (4.5.17) is inapplicable in the neutral case, so define
$$I:\{o^A,\iota^A\} \mapsto \{-\iota^A,o^A\},\eqno(2.3.1)$$
which preserves $\chi$ and is anti-involutary ($I^2 = -1$), i.e., induces an anti-involutary symplectomorphism of $(S,\epsilon_{AB})$. The induced operation on the dual basis is
$$\{-\chi^{-1}\iota_A,\chi^{-1}o_A\} \mapsto \{-\chi^{-1}o_A,-\chi^{-1}\iota_A\}.\eqno(2.3.2)$$
If one defines the analogous operation ($\tilde I$) for frames of $S'$, then $I$ and $\tilde I$ together induce, for the associated null tetrad 
$$\ell^a := o^Ao^{A'} \hskip .5in \tilde m^a := \iota^Ao^{A'} \hskip .5in n^a := \iota^A \iota^{A'} \hskip .5in m^a := o^A\iota^{A'}\eqno(2.3.3)$$
the transformation
$$(\ell)' = n \hskip .5in (n)' = \ell \hskip .5in (m)' = -\tilde m \hskip .5in (\tilde m)' = -m,\eqno(2.3.4)$$
which differs from [22], (4.5.18). Of course, the choice of $I$ (and $\tilde I$) is neither unique nor canonical, merely convenient. One could choose instead $-I$, for example, but as long as one also then chooses $-\tilde I$, (2.3.4) is left unchanged. The choice for $S$ does not of course determine that for $S'$, convenience is the only guide. But in order for the priming operation to `commute' with the tilde notation on the spin spaces themselves, one must choose both $I$ or both $-I$ as the priming operation on the spin spaces. These choices are, also, the only ones yielding transformations of the null tetrad as simple as (2.3.4), i.e., merely interchanging, up to sign, the elements of the null tetrad. Therefore, I take the priming operation in the neutral case to be the pair $I$ and $\tilde I$ so (2.3.1--4) are valid.

Under a (passive) replacement of a basis by its image under $I$ ($\tilde I$), the components of a spinor with respect to the transformed basis are simply related to the components with respect to the initial basis. One can describe the transformation of components by a prime notation as follows. First note that each member of the basis is replaced by the other up to sign. Define
$$[{\bf A}] = \cases{0,\cr 1,\cr}\hskip 1in\hbox{iff}\hskip 1in {\bf A} = \cases{1,\cr 0,\cr}.$$
Since $I$ and $\tilde I$ act in the same manner, it suffices to consider a spinor with unprimed indices to illustrate:
$$\eqalign{\xi^{A\ldots C}{}_{L\ldots N} &= \xi^{\bf A\ldots C}{}_{\bf L\ldots N}\epsilon_{\bf A}{}^A\ldots\epsilon_{\bf C}{}^C\epsilon_L{}^{\bf L}\ldots\epsilon_N{}^{\bf N}\cr
&= (-1)^p\xi^{\bf A\ldots C}{}_{\bf L\ldots N}I\bigl(\epsilon_{[{\bf A}]}{}^A\bigr)\ldots I\bigl(\epsilon_{[{\bf C}]}{}^C\bigr)I\bigl(\epsilon_L{}^{[{\bf L}]}\bigr)\ldots I\bigl(\epsilon_N{}^{[{\bf N}]}\bigr)\cr
&= (-1)^p\xi^{\bf A\ldots C}{}_{\bf L\ldots N}(\epsilon')_{[{\bf A}]}{}^A\ldots(\epsilon')_{[{\bf C}]}{}^C(\epsilon')_L{}^{[{\bf L}]}\ldots(\epsilon')_N{}^{[{\bf N}]}\cr}$$
where $(\epsilon')_{\bf A}{}^A := I\bigl(\epsilon_{\bf A}{}^A\bigr)$ and $p$ is the number of 1's amongst ${\bf A},\ldots,{\bf C},{\bf L},\ldots,{\bf N}$, hence
$$\xi'^{[{\bf A}]\ldots[{\bf C}]}{}_{[{\bf L}]\ldots[{\bf N}]} = (-1)^{p}\xi^{\bf A\ldots C}{}_{\bf L\ldots N}.$$

Since $I^2 = -1$, the square of the priming operation replaces components of odd rank spinors by their negatives and leaves components of even rank spinors unchanged.

In particular, for the spinor connection symbols, for which there are four indices, there is an odd number of 1's (primed and unprimed) iff there is an odd number of 0's (primed and unprimed). Thus, under the priming operation, one includes a negative sign iff there is an odd number of 1's (equivalently, an odd number of 0's) irrespective of whether the index is primed or not. For example,
$$\eqalign{\gamma_{10'1}{}^0 &= \epsilon_1{}^A\epsilon_{0'}{}^{A'}\epsilon_E{}^0\nabla_{AA'}\epsilon_1{}^E\cr
&= -I(\epsilon_0{}^A)I(\epsilon_{1'}{}^{A'})I(\epsilon_E{}^{1'})\nabla_{AA'}(-I(\epsilon_{0'}{}^E))\cr
&= -(\epsilon')_0{}^A(\epsilon')_{1'}{}^{A'}(\epsilon')_E{}^1\nabla_{AA'}(-(\epsilon')_0{}^E)\cr
&= \gamma'_{01'0}{}^1.\cr}$$
The spin coefficients are unchanged under the square of the priming operation.
\vskip 24pt
One may therefore introduce the notation
\vskip 12pt
$$\gamma_{{\bf A}{\bf A}'{\bf B}}{}^{\bf C}\hskip .5in =: \hskip .5in\vcenter{\offinterlineskip
\halign{&\vrule#&\strut\ #\ \cr
\noalign{\medskip}
\noalign{\hrule}
height3pt&\omit&&\omit&&\omit&&\omit&&\omit&\cr
&\hfil${\atop{\bf A}}{\atop{\bf A}'}{\atop{\bf B}}{{\bf C}\atop}$\hfil&&\hfil${\atop 0}{0\atop}$\hfil&&\hfil${\atop 0}{1\atop}$\hfil&&\hfil${\atop 1}{0\atop}$\hfil&&\hfil${\atop 1}{1 \atop}$\hfil&\cr
height3pt&\omit&&\omit&&\omit&&\omit&&\omit&\cr
\noalign{\hrule}
height3pt&\omit&&\omit&&\omit&&\omit&&\omit&\cr
&\hfil$00'$\hfil&&\hfil$\epsilon$\hfil&&\hfil$\kappa$\hfil&&\hfil$-\tau'$\hfil&&\hfil$\gamma'$\hfil&\cr
height3pt&\omit&&\omit&&\omit&&\omit&&\omit&\cr
\noalign{\hrule}
height3pt&\omit&&\omit&&\omit&&\omit&&\omit&\cr
&\hfil$10'$\hfil&&\hfil$\alpha$\hfil&&\hfil$\rho$\hfil&&\hfil$\sigma'$\hfil&&\hfil$-\beta'$\hfil&\cr
height3pt&\omit&&\omit&&\omit&&\omit&&\omit&\cr
\noalign{\hrule}
height3pt&\omit&&\omit&&\omit&&\omit&&\omit&\cr
&\hfil$01'$\hfil&&\hfil$\beta$\hfil&&\hfil$\sigma$\hfil&&\hfil$\rho'$\hfil&&\hfil$-\alpha'$\hfil&\cr
height3pt&\omit&&\omit&&\omit&&\omit&&\omit&\cr
\noalign{\hrule}
height3pt&\omit&&\omit&&\omit&&\omit&&\omit&\cr
&\hfil$11'$\hfil&&\hfil$\gamma$\hfil&&\hfil$\tau$\hfil&&\hfil$-\kappa'$\hfil&&\hfil$\epsilon'$\hfil&\cr
height3pt&\omit&&\omit&&\omit&&\omit&&\omit&\cr
\noalign{\hrule}}}\eqno(2.4)$$
\vskip 12pt
This notation differs from [22], (4.5.16), due to the different priming operation and, because this difference is unavoidable, I have also taken the liberty to dispense with the negative signs in the second column of [22], (4.5.16) (admitted as unfortunate in [22], p. 226).

Since $\nabla_a$ independently induces connections on $S_M$ and $S'_M$, (2.4) itself imposes no constraints on how one denotes the $\tilde\gamma_{{\bf A}{\bf A}'{\bf C}'}{}^{{\bf B}'}$. There is of course the constraint imposed by the decision to employ the priming notation and the particular choice $\tilde I$. Thus, at first thought, the spin coefficients arising from the action of $\nabla_a$ on $S'$ might best be denoted exactly as in (2.4), only with tildes over them. This choice, however, does not actually reflect the tilde convention introduced in \S 1. To follow that protocol suggests introducing notation which parallels the Lorentzian case. I emphasize that this choice is purely one of utility but has the advantages described in \S 1. The $\tilde\gamma_{{\bf A}{\bf A}'{\bf B}'}{}^{{\bf C}'}$ are of course entirely independent of the $\gamma_{{\bf A}{\bf A}'{\bf C}}{}^{\bf B}$.

Hence, define
\vskip 12pt
$$\tilde\gamma_{{\bf A}{\bf A}'{\bf B}'}{}^{{\bf C}'}\hskip .5in =: \hskip .5in\vcenter{\offinterlineskip
\halign{&\vrule#&\strut\ #\ \cr
\noalign{\medskip}
\noalign{\hrule}
height3pt&\omit&&\omit&&\omit&&\omit&&\omit&\cr
&\hfil${\atop{\bf A}}{\atop{\bf A}'}{\atop{\bf B}'}{{\bf C}'\atop}$\hfil&&\hfil${\atop 0'}{0'\atop}$\hfil&&\hfil${\atop 0'}{1'\atop}$\hfil&&\hfil${\atop 1'}{0'\atop}$\hfil&&\hfil${\atop 1'}{1' \atop}$\hfil&\cr
height3pt&\omit&&\omit&&\omit&&\omit&&\omit&\cr
\noalign{\hrule}
height3pt&\omit&&\omit&&\omit&&\omit&&\omit&\cr
&\hfil$00'$\hfil&&\hfil$\tilde\epsilon$\hfil&&\hfil$\tilde\kappa$\hfil&&\hfil$-\tilde\tau'$\hfil&&\hfil$\tilde\gamma'$\hfil&\cr
height3pt&\omit&&\omit&&\omit&&\omit&&\omit&\cr
\noalign{\hrule}
height3pt&\omit&&\omit&&\omit&&\omit&&\omit&\cr
&\hfil$10'$\hfil&&\hfil$\tilde\beta$\hfil&&\hfil$\tilde\sigma$\hfil&&\hfil$\tilde\rho'$\hfil&&\hfil$-\tilde\alpha'$\hfil&\cr
height3pt&\omit&&\omit&&\omit&&\omit&&\omit&\cr
\noalign{\hrule}
height3pt&\omit&&\omit&&\omit&&\omit&&\omit&\cr
&\hfil$01'$\hfil&&\hfil$\tilde\alpha$\hfil&&\hfil$\tilde\rho$\hfil&&\hfil$\tilde\sigma'$\hfil&&\hfil$-\tilde\beta'$\hfil&\cr
height3pt&\omit&&\omit&&\omit&&\omit&&\omit&\cr
\noalign{\hrule}
height3pt&\omit&&\omit&&\omit&&\omit&&\omit&\cr
&\hfil$11'$\hfil&&\hfil$\tilde\gamma$\hfil&&\hfil$\tilde\tau$\hfil&&\hfil$-\tilde\kappa'$\hfil&&\hfil$\tilde\epsilon'$\hfil&\cr
height3pt&\omit&&\omit&&\omit&&\omit&&\omit&\cr
\noalign{\hrule}}}\eqno(2.5)$$
\vskip 12pt
Now define
$$\displaylines{D := \ell^a\nabla_a = o^Ao^{A'}\nabla_{AA'} \hskip 1.25in D' := n^a\nabla_a = \iota^A\iota^{A'}\nabla_{AA'}\cr
\noalign{\vskip -12pt}
\hfill\llap(2.6)\cr
\noalign{\vskip -12pt}
\delta := m^a\nabla_a = o^A\iota^{A'}\nabla_{AA'} \hskip 1.25in \triangle := \tilde m^a\nabla_a = \iota^Ao^{A'}\nabla_{AA'}.\cr}$$
Note that while the priming operation (2.3) does interchange $D$ and $D'$, it switches $\delta$ and $-\triangle$ (i.e., $\triangle = -\delta'$). For this reason, in the neutral case it makes sense to have an independent notation for $\tilde m^a\nabla_a$, reflecting the fact that $\tilde m^a$ is of course not determined by $m^a$.

With the notation for spin coefficients now set, it is a simple matter to transcribe relevant equations from [22], Ch. 4, to their neutral analogues. In brief, in any equation involving spin coefficients in [22], one need only replace $\kappa$, $\rho$, $\sigma$, $\tau$, $\sigma'$, $\rho'$, $\alpha'$ and $\beta'$ by their negatives (and likewise for the tilde versions of these quantities), replace each occurrence of $\delta'$ in [22] by $\triangle$, and take any other independent steps which might be necessary to pass to the neutral analogue. Nevertheless, it will be useful to record here once and for all the important equations from [22], \S\S 4.5, 4.11, \& 4.13 in the form they take in the neutral context.  

From the definitions (2.2), and then expressing the right-hand-side in terms of the null tetrad rather than the spin bases,
\vskip 12pt
$$\displaylines{\vcenter{\offinterlineskip
\halign{&\vrule#&\strut\ #\ \cr
\noalign{\hrule}
height3pt&\omit&&\omit&&\omit&&\omit&\cr
&\hfil$\epsilon$\hfil&&\hfil$\kappa$\hfil&&\hfil$-\tau'$\hfil&&\hfil$\gamma'$\hfil&\cr
height3pt&\omit&&\omit&&\omit&&\omit&\cr
\noalign{\hrule}
height3pt&\omit&&\omit&&\omit&&\omit&\cr
&\hfil$\alpha$\hfil&&\hfil$\rho$\hfil&&\hfil$\sigma'$\hfil&&\hfil$-\beta'$\hfil&\cr
height3pt&\omit&&\omit&&\omit&&\omit&\cr
\noalign{\hrule}
height3pt&\omit&&\omit&&\omit&&\omit&\cr
&\hfil$\beta$\hfil&&\hfil$\sigma$\hfil&&\hfil$\rho'$\hfil&&\hfil$-\alpha'$\hfil&\cr
height3pt&\omit&&\omit&&\omit&&\omit&\cr
\noalign{\hrule}
height3pt&\omit&&\omit&&\omit&&\omit&\cr
&\hfil$\gamma$\hfil&&\hfil$\tau$\hfil&&\hfil$-\kappa'$\hfil&&\hfil$\epsilon'$\hfil&\cr
height3pt&\omit&&\omit&&\omit&&\omit&\cr
\noalign{\hrule}}}
\hskip .5in=\hskip .5in\chi^{-1}\times
\vcenter{\offinterlineskip
\halign{&\vrule#&\strut\ #\ \cr
\noalign{\hrule}
height3pt&\omit&&\omit&&\omit&&\omit&\cr
&\hfil$-\iota_BDo^B$\hfil&&\hfil$o_BDo^B$\hfil&&\hfil$-\iota_BD\iota^B$\hfil&&\hfil$o_BD\iota^B$\hfil&\cr
height3pt&\omit&&\omit&&\omit&&\omit&\cr
\noalign{\hrule}
height3pt&\omit&&\omit&&\omit&&\omit&\cr
&\hfil$-\iota_B\triangle o^B$\hfil&&\hfil$o_B\triangle o^B$\hfil&&\hfil$-\iota_B\triangle\iota^B$\hfil&&\hfil$o_B\triangle \iota^B$\hfil&\cr
height3pt&\omit&&\omit&&\omit&&\omit&\cr
\noalign{\hrule}
height3pt&\omit&&\omit&&\omit&&\omit&\cr
&\hfil$-\iota_B\delta o^B$\hfil&&\hfil$o_B\delta o^B$\hfil&&\hfil$-\iota_B\delta \iota^B$\hfil&&\hfil$o_B\delta\iota^B$\hfil&\cr
height3pt&\omit&&\omit&&\omit&&\omit&\cr
\noalign{\hrule}
height3pt&\omit&&\omit&&\omit&&\omit&\cr
&\hfil$-\iota_BD'o^B$\hfil&&\hfil$o_BD'o^B$\hfil&&\hfil$-\iota_BD'\iota^B$\hfil&&\hfil$o_BD'\iota^B$\hfil&\cr
height3pt&\omit&&\omit&&\omit&&\omit&\cr
\noalign{\hrule}}}\cr
\noalign{\vskip 6pt}
\hfill\llap(2.7)\cr
\noalign{\vskip 6pt}
= \chi^{-1}\tilde\chi^{-1}\times
\vcenter{\offinterlineskip
\halign{&\vrule#&\strut\ #\ \cr
\noalign{\hrule}
height3pt&\omit&&\omit&&\omit&&\omit&\cr
&\hfil${1 \over 2}(n_aD\ell^a + m_aD\tilde m^a + \tilde\chi D\chi)$\hfil&&\hfil$-m_aD\ell^a$\hfil&&\hfil$-\tilde m_aDn^a$\hfil&&\hfil${1 \over 2}(\ell_aDn^a + \tilde m_aDm^a + \tilde\chi D\chi)$\hfil&\cr
height3pt&\omit&&\omit&&\omit&&\omit&\cr
\noalign{\hrule}
height3pt&\omit&&\omit&&\omit&&\omit&\cr
&\hfil${1 \over 2}(n_a\triangle\ell^a + m_a\triangle\tilde m^a + \tilde\chi\triangle\chi)$\hfil&&\hfil$-m_a\triangle\ell^a$\hfil&&\hfil$-\tilde m_a\triangle n^a$\hfil&&\hfil${1 \over 2}(\ell_a\triangle n^a + \tilde m_a\triangle m^a + \tilde\chi \triangle\chi)$\hfil&\cr
height3pt&\omit&&\omit&&\omit&&\omit&\cr
\noalign{\hrule}
height3pt&\omit&&\omit&&\omit&&\omit&\cr
&\hfil${1 \over 2}(n_a\delta\ell^a + m_a\delta\tilde m^a + \tilde\chi\delta\chi)$\hfil&&\hfil$-m_a\delta\ell^a$\hfil&&\hfil$-\tilde m_a\delta n^a$\hfil&&\hfil${1 \over 2}(\ell_a\delta n^a + \tilde m_a\delta m^a + \tilde\chi \delta\chi)$\hfil&\cr
height3pt&\omit&&\omit&&\omit&&\omit&\cr
\noalign{\hrule}
height3pt&\omit&&\omit&&\omit&&\omit&\cr
&\hfil${1 \over 2}(n_aD'\ell^a + m_aD'\tilde m^a + \tilde\chi D'\chi)$\hfil&&\hfil$-m_aD'\ell^a$\hfil&&\hfil$-\tilde m_aD'n^a$\hfil&&\hfil${1 \over 2}(\ell_aD'n^a + \tilde m_aD'm^a + \tilde\chi D'\chi)$\hfil&\cr
height3pt&\omit&&\omit&&\omit&&\omit&\cr
\noalign{\hrule}}}\cr}$$
\vskip 12pt
and
\vskip 12pt
$$\displaylines{\vcenter{\offinterlineskip
\halign{&\vrule#&\strut\ #\ \cr
\noalign{\hrule}
height3pt&\omit&&\omit&&\omit&&\omit&\cr
&\hfil$\tilde\epsilon$\hfil&&\hfil$\tilde\kappa$\hfil&&\hfil$-\tilde\tau'$\hfil&&\hfil$\tilde\gamma'$\hfil&\cr
height3pt&\omit&&\omit&&\omit&&\omit&\cr
\noalign{\hrule}
height3pt&\omit&&\omit&&\omit&&\omit&\cr
&\hfil$\tilde\beta$\hfil&&\hfil$\tilde\sigma$\hfil&&\hfil$\tilde\beta'$\hfil&&\hfil$-\tilde\alpha'$\hfil&\cr
height3pt&\omit&&\omit&&\omit&&\omit&\cr
\noalign{\hrule}
height3pt&\omit&&\omit&&\omit&&\omit&\cr
&\hfil$\tilde\alpha$\hfil&&\hfil$\tilde\rho$\hfil&&\hfil$\tilde\sigma'$\hfil&&\hfil$-\tilde\beta'$\hfil&\cr
height3pt&\omit&&\omit&&\omit&&\omit&\cr
\noalign{\hrule}
height3pt&\omit&&\omit&&\omit&&\omit&\cr
&\hfil$\tilde\gamma$\hfil&&\hfil$\tilde\tau$\hfil&&\hfil$-\tilde\kappa'$\hfil&&\hfil$\tilde\epsilon'$\hfil&\cr
height3pt&\omit&&\omit&&\omit&&\omit&\cr
\noalign{\hrule}}}
\hskip .5in=\hskip .5in\tilde\chi^{-1}\times
\vcenter{\offinterlineskip
\halign{&\vrule#&\strut\ #\ \cr
\noalign{\hrule}
height3pt&\omit&&\omit&&\omit&&\omit&\cr
&\hfil$-\iota_{B'}Do^{B'}$\hfil&&\hfil$o_{B'}Do^{B'}$\hfil&&\hfil$-\iota_{B'}D\iota^{B'}$\hfil&&\hfil$o_{B'}D\iota^{B'}$\hfil&\cr
height3pt&\omit&&\omit&&\omit&&\omit&\cr
\noalign{\hrule}
height3pt&\omit&&\omit&&\omit&&\omit&\cr
&\hfil$-\iota_{B'}\triangle o^{B'}$\hfil&&\hfil$o_{B'}\triangle o^{B'}$\hfil&&\hfil$-\iota_{B'}\triangle\iota^{B'}$\hfil&&\hfil$o_{B'}\triangle \iota^{B'}$\hfil&\cr
height3pt&\omit&&\omit&&\omit&&\omit&\cr
\noalign{\hrule}
height3pt&\omit&&\omit&&\omit&&\omit&\cr
&\hfil$-\iota_{B'}\delta o^{B'}$\hfil&&\hfil$o_{B'}\delta o^{B'}$\hfil&&\hfil$-\iota_{B'}\delta \iota^{B'}$\hfil&&\hfil$o_{B'}\delta\iota^{B'}$\hfil&\cr
height3pt&\omit&&\omit&&\omit&&\omit&\cr
\noalign{\hrule}
height3pt&\omit&&\omit&&\omit&&\omit&\cr
&\hfil$-\iota_{B'}D'o^{B'}$\hfil&&\hfil$o_{B'}D'o^{B'}$\hfil&&\hfil$-\iota_{B'}D'\iota^{B'}$\hfil&&\hfil$o_{B'}D'\iota^{B'}$\hfil&\cr
height3pt&\omit&&\omit&&\omit&&\omit&\cr
\noalign{\hrule}}}\cr
\noalign{\vskip 6pt}
\hfill\llap(2.8)\cr
\noalign{\vskip 6pt}
= \chi^{-1}\tilde\chi^{-1}\times
\vcenter{\offinterlineskip
\halign{&\vrule#&\strut\ #\ \cr
\noalign{\hrule}
height3pt&\omit&&\omit&&\omit&&\omit&\cr
&\hfil${1 \over 2}(n_aD\ell^a + \tilde m_aD m^a + \chi D\tilde\chi)$\hfil&&\hfil$-\tilde m_aD\ell^a$\hfil&&\hfil$-m_aDn^a$\hfil&&\hfil${1 \over 2}(\ell_aDn^a + m_aD\tilde m^a + \chi D\tilde\chi)$\hfil&\cr
height3pt&\omit&&\omit&&\omit&&\omit&\cr
\noalign{\hrule}
height3pt&\omit&&\omit&&\omit&&\omit&\cr
&\hfil${1 \over 2}(n_a\triangle\ell^a + \tilde m_a\triangle m^a + \chi\triangle\tilde\chi)$\hfil&&\hfil$-\tilde m_a\triangle\ell^a$\hfil&&\hfil$-m_a\triangle n^a$\hfil&&\hfil${1 \over 2}(\ell_a\triangle n^a + m_a\triangle \tilde m^a + \chi \triangle\tilde\chi)$\hfil&\cr
height3pt&\omit&&\omit&&\omit&&\omit&\cr
\noalign{\hrule}
height3pt&\omit&&\omit&&\omit&&\omit&\cr
&\hfil${1 \over 2}(n_a\delta\ell^a + \tilde m_a\delta m^a + \chi\delta\tilde\chi)$\hfil&&\hfil$-\tilde m_a\delta\ell^a$\hfil&&\hfil$-m_a\delta n^a$\hfil&&\hfil${1 \over 2}(\ell_a\delta n^a + m_a\delta \tilde m^a + \chi \delta\tilde\chi)$\hfil&\cr
height3pt&\omit&&\omit&&\omit&&\omit&\cr
\noalign{\hrule}
height3pt&\omit&&\omit&&\omit&&\omit&\cr
&\hfil${1 \over 2}(n_aD'\ell^a + \tilde m_aD' m^a + \chi D'\tilde\chi)$\hfil&&\hfil$-\tilde m_aD'\ell^a$\hfil&&\hfil$-m_aD'n^a$\hfil&&\hfil${1 \over 2}(\ell_aD'n^a + m_aD'\tilde m^a + \chi D'\tilde\chi)$\hfil&\cr
height3pt&\omit&&\omit&&\omit&&\omit&\cr
\noalign{\hrule}}}\cr}$$
\vskip 12pt
Equivalently,
$$\vcenter{ \halign{$\hfil#$&&${}#\hfil$&\qquad$\hfil#$\cr
D o^A &= \epsilon o^A + \kappa \iota^A & Do^{A'} &= \tilde\epsilon o^{A'} + \tilde\kappa\iota^{A'}\cr
\delta o^A &= \beta o^A + \sigma\iota^A & \triangle o^{A'} &= \tilde\beta o^{A'} + \tilde\sigma\iota^{A'}\cr
\triangle o^A &= \alpha o^A + \rho\iota^A & \delta o^{A'} &= \tilde\alpha o^{A'} + \tilde\rho\iota^{A'}\cr
D' o^A &= \gamma o^A + \tau\iota^A & D' o^{A'} &= \tilde\gamma o^{A'} + \tilde\tau\iota^{A'}\cr
D\iota^A &= -\tau' o^A + \gamma'\iota^A & D\iota^{A'} &= -\tilde\tau' o^{A'} + \tilde\gamma'\iota^{A'}\cr
\delta \iota^A &= \rho' o^A - \alpha'\iota^A & \triangle \iota^{A'} &= \tilde\rho' o^{A'} - \tilde\alpha' \iota^{A'}\cr
\triangle \iota^A &= \sigma' o^{A'} - \beta'\iota^A & \delta \iota^{A'} &= \tilde\sigma' o^{A'} - \tilde\beta'\iota^{A'}\cr
D'\iota^A &= -\kappa' o^A + \epsilon'\iota^A & D'\iota^{A'} &= -\tilde\kappa' o^{A'} + \tilde\epsilon'\iota^{A'},\cr}}\eqno(2.9)$$
and,
$$\vcenter{\openup1\jot \halign{$\hfil#$&&${}#\hfil$&\qquad$\hfil#$\cr
D\ell^a &= (\epsilon + \tilde\epsilon)\ell^a + \tilde\kappa m^a + \kappa\tilde m^a & Dn^a &= (\gamma' + \tilde\gamma')n^a - \tau'm^a - \tilde\tau'\tilde m^a\cr
\triangle \ell^a &= (\alpha + \tilde\beta)\ell^a + \tilde\sigma m^a + \rho\tilde m^a & \triangle n^a &= -(\tilde\alpha' + \beta')n^a + \sigma' m^a + \tilde\rho'\tilde m^a\cr
\delta\ell^a &= (\tilde\alpha + \beta)\ell^a + \tilde\rho m^a + \sigma \tilde m^a & \delta n^a &= -(\alpha' + \tilde\beta')n^a + \rho' m^a + \tilde\sigma'\tilde m^a\cr
D'\ell^a &= (\gamma + \tilde\gamma)\ell^a + \tilde\tau m^a + \tau\tilde m^a & D'n^a &= (\epsilon' + \tilde\epsilon')n^a - \kappa' m^a - \tilde\kappa'\tilde m^a\cr
Dm^a &= (\epsilon + \tilde\gamma')m^a - \tilde\tau'\ell^a + \kappa n^a & D\tilde m^a &= (\gamma' + \tilde\epsilon)\tilde m^a - \tau'\ell^a + \tilde\kappa n^a\cr
\triangle m^a &= (\alpha - \tilde\alpha')m^a + \tilde\rho'\ell^a + \rho n^a & \triangle \tilde m^a &= (\tilde\beta - \beta')\tilde m^a + \sigma'\ell^a + \tilde\sigma n^a\cr
\delta m^a &= (\beta - \tilde\beta')m^a + \tilde\sigma'\ell^a + \sigma n^a & \delta \tilde m^a &= (\tilde\alpha - \alpha')\tilde m^a + \rho'\ell^a + \tilde\rho n^a\cr
D' m^a &= (\gamma + \tilde\epsilon')m^a - \tilde\kappa'\ell^a + \tau n^a & D'\tilde m^a &= (\epsilon' + \tilde\gamma)\tilde m^a - \kappa'\ell^a + \tilde\tau n^a.\cr}}\eqno(2.10)$$
\vskip 6pt
Note that each set of equations (2.9) and (2.10) is invariant under each of the priming operation and the application of tilde.

When the frame $\{o^A,\iota^A\}$ is a spin frame ($\chi = 1$), then, with $\nabla_{\bf AA'} := \epsilon_{\bf A}{}^A\epsilon_{\bf A'}{}^A\nabla_{AA'}$,
$$0 = \nabla_{\bf AA'}\epsilon_{\bf BC} = \nabla_{\bf AA'}\left(\epsilon_{\bf B}{}^A\,\epsilon_{A{\bf C}}\right) = \epsilon_{A{\bf C}}\nabla_{\bf AA'}\epsilon_{\bf B}{}^A - \epsilon_{{\bf B}A}\nabla_{\bf AA'}\epsilon^A{}_{\bf C} = \gamma_{\bf AA'BC} - \gamma_{\bf AA'CB},$$
i.e., $\epsilon_{\bf DC}\gamma_{\bf AA'B}{}^{\bf D} = \epsilon_{\bf EB}\gamma_{\bf AA'C}{}^{\bf E}$, whence $\epsilon_{01}\gamma_{{\bf AA'}0}{}^0 = \epsilon_{10}\gamma_{{\bf AA'}1}{}^1$, i.e., from (2.4)
$$\epsilon = -\gamma' \hskip .75in \alpha = \beta' \hskip .75in \beta = \alpha' \hskip .75in \gamma = -\epsilon'.\eqno(2.11.1))$$
Similarly, if $\{o^{A'},\iota^{A'}\}$ is a spin frame ($\tilde\chi = 1$), then
$$\tilde\epsilon = -\tilde\gamma' \hskip .75in \tilde\alpha = \tilde\beta' \hskip .75in \tilde\beta = \tilde\alpha' \hskip .75in \tilde\gamma = -\tilde\epsilon'.\eqno(2.11.2)$$

The neutral analogues of [22], (4.5.30--43), are obtained by simply replacing complex conjugation signs by tilde. [22], (4.5.43), and its tilde version enable one to compute the spin coefficients from the Infeld-van der Waerden symbols.
\vskip 24pt
\noindent {\section 3. Curvature Spinors and Spin Coefficients}
\vskip 12pt
The curvature spinors are defined exactly as in [22] \S 4.6, allowing for the differences arising from neutral signature, e.g., the Hodge star operator on two-forms is an involution rather than almost complex structure. As usual complex conjugation signs are replaced by tilde. The conventions I follow for the choices involved in defining the Riemann tensor and derived quantities are detailed in [17], Appendix One. My conventions coincide with those of [22] except that my Ricci tensor is the negative of theirs. This choice requires a simple modification of [22] (4.6.21--23), namely multiplying one side of each equation by -1, to preserve the definitions of the curvature spinors themselves. This one difference of sign choice therefore only impacts equations which relate $\Phi_{ABC'D'}$ to $R_{ab}$ and/or $\Lambda$ to $S$ (the Ricci scalar curvature) and not any equation involving curvature spinors derived directly from expressions involving the Riemann tensor, e.g., [22], (4.11.12). For convenience, the Appendix records curvature formulae used often in this paper; see also [17], Appendix Two.

[22], (4.11.1--10), carry over unchanged except for replacing complex conjugation signs by tilde and omitting the minus signs on the right-hand sides of (4.11.10) (due to my convention for the Ricci tensor).

In place of [22], (4.11.11), one obtains six {\sl independent} equations (though the third is the prime of the second, the fourth the tilde of the second, and the fifth the prime of the fourth and the tilde of the third):
$$\eqalign{[D',D] &= (\gamma + \tilde\gamma)D - (\gamma' + \tilde\gamma')D' + (\tau + \tilde\tau')\triangle + (\tau' + \tilde\tau)\delta\cr
[\delta,D] &= (\beta + \tilde\alpha + \tilde\tau')D - \kappa D' + \sigma\triangle + (\tilde\rho - \epsilon - \tilde\gamma')\delta\cr
[D',\triangle] &= (\beta' + \tilde\alpha' + \tilde\tau)D' - \kappa'D - \sigma'\delta - (\tilde\rho' - \epsilon' - \tilde\gamma)\triangle\cr
[D,\triangle] &= \tilde\kappa D' - (\tau' + \tilde\beta + \alpha)D - \tilde\sigma\delta - (\rho - \tilde\epsilon - \gamma')\triangle\cr
[\delta,D'] &= \tilde\kappa'D - (\tau + \tilde\beta' + \alpha')D' + \tilde\sigma'\triangle + (\rho' - \tilde\epsilon' - \gamma)\delta\cr
[\triangle,\delta] &= (\tilde\rho' - \rho')D + (\rho - \tilde\rho)D' + (\alpha' - \tilde\alpha)\triangle + (\alpha - \tilde\alpha')\delta.\cr}\eqno(3.1)$$
This set of equations is invariant under the priming operation and application of tilde.

The equations of [22], (4.11.12), are grouped into pairs related to each other via the priming operation. One must therefore take account of the effect of (2.3) on the dyad components of the curvature spinors; in particular (using the standard notation of [22], (4.11.6--10)):
$$\displaylines{\Psi_0' = \Psi_4 \qquad \Psi_1' = -\Psi_3 \qquad \Psi_2' = \Psi_2 \qquad \Psi_3' = -\Psi_1 \qquad \Psi_4' = \Psi_0\cr
\noalign{\vskip -12pt}
\hfill\llap(3.2)\cr
\noalign{\vskip -12pt}
\Phi_{00}'= \Phi_{22} \qquad \Phi_{01}' = -\Phi_{21} \qquad \Phi_{02}' = \Phi_{20}\cr
\Phi_{10}' = -\Phi_{12} \qquad \Phi_{11}' = \Phi_{11} \qquad \Phi_{12}' = -\Phi_{10}\cr
\Phi_{20}' = \Phi_{02} \qquad \Phi_{21}' = -\Phi_{01} \qquad \Phi_{22}' = \Phi_{00}.\cr}$$
Note that all quantities involved are unchanged by applying the priming operation twice. Moreover,
$$\Pi := \chi\tilde\chi\Lambda,\eqno(3.3)$$
is invariant under the priming operation since $\chi$ and $\tilde\chi$ are. Therefore, in the neutral case, [22], (4.11.5), gives rise to the following analogue of [22], (4.11.12):
$$\eqalignno{&&(3.4)\cr
\triangle\kappa - D\rho &= \rho^2 + \sigma\tilde\sigma - \tilde\kappa\tau + \kappa(\tau' + 2\alpha + \tilde\beta + \beta') - \rho(\epsilon + \tilde\epsilon) + \Phi_{00}&({\rm a})\cr
-\delta\kappa' - D'\rho' &= (\rho')^2 + \sigma'\tilde\sigma' - \tilde\kappa'\tau' + \kappa'(\tau + 2\alpha' + \tilde\beta' + \beta) - \rho'(\epsilon' + \tilde\epsilon') + \Phi_{22}&({\rm a}')\cr
\delta\kappa - D\sigma &= \sigma(\rho + \tilde\rho - \tilde\gamma' + \gamma' - 2\epsilon) - \kappa(\tau - \tilde\tau' - \tilde\alpha - \alpha' - 2\beta) + \Psi_0 &({\rm b})\cr
-\triangle\kappa' - D'\sigma' &= \sigma'(\rho' + \tilde\rho' - \tilde\gamma + \gamma - 2\epsilon') - \kappa'(\tau' - \tilde\tau - \tilde\alpha' - \alpha - 2\beta') + \Psi_4 &({\rm b}')\cr
D'\kappa - D\tau &= \rho(\tau + \tilde\tau') + \sigma(\tilde\tau + \tau') - \tau(\tilde\gamma' + \epsilon) + \kappa(\tilde\gamma + 2\gamma - \epsilon') + \Psi_1 + \Phi_{01}&({\rm c})\cr
D\kappa' - D'\tau' &= \rho'(\tau' + \tilde\tau) + \sigma'(\tilde\tau' + \tau) - \tau'(\tilde\gamma + \epsilon') + \kappa'(\tilde\gamma' + 2\gamma' - \epsilon) - \Psi_3 - \Phi_{21} &({\rm c}')\cr
\triangle\sigma - \delta\rho &= \tau(\rho - \tilde\rho) + \kappa(\tilde\rho' - \rho') - \rho(\tilde\alpha + \beta) + \sigma(2\alpha - \tilde\alpha' + \beta') - \Psi_1 + \Phi_{01} &({\rm d})\cr
- \delta\sigma' + \triangle\rho' &= \tau'(\rho' - \tilde\rho') + \kappa'(\tilde\rho - \rho) - \rho'(\tilde\alpha' + \beta') + \sigma'(2\alpha' - \tilde\alpha + \beta) + \Psi_3 - \Phi_{21} &({\rm d}')\cr
D'\sigma - \delta\tau &= -\rho'\sigma - \tilde\sigma'\rho + \tau^2 - \kappa\tilde\kappa' - \tau(\beta - \tilde\beta') + \sigma(2\gamma - \epsilon' + \tilde\epsilon') + \Phi_{02}&({\rm e})\cr
D\sigma' + \triangle\tau' &= -\rho\sigma' - \tilde\sigma\rho' + (\tau')^2 - \kappa'\tilde\kappa - \tau'(\beta' - \tilde\beta) + \sigma'(2\gamma' - \epsilon + \tilde\epsilon) + \Phi_{20} &({\rm e}')\cr
\triangle\tau - D'\rho &= \rho\tilde\rho' + \sigma\sigma' - \tau\tilde\tau + \kappa\kappa' - \rho(\gamma + \tilde\gamma) + \tau(\alpha - \tilde\alpha') - \Psi_2 - 2\Pi&({\rm f})\cr
-\delta\tau' - D\rho' &= \rho'\tilde\rho + \sigma'\sigma - \tau'\tilde\tau' + \kappa'\kappa - \rho'(\gamma' + \tilde\gamma') + \tau'(\alpha' - \tilde\alpha) - \Psi_2 - 2\Pi&({\rm f}')\cr
D'\beta - \delta\gamma &= \tau\rho' + \kappa'\sigma - \tilde\kappa'\epsilon - \alpha\tilde\sigma' + \beta(\tilde\epsilon' - \rho' + \gamma) + \gamma(\tilde\beta' + \alpha' + \tau) - \Phi_{12} &({\rm g})\cr
D\beta' + \triangle\gamma' &= \tau'\rho + \kappa\sigma' - \tilde\kappa\epsilon' - \alpha'\tilde\sigma + \beta'(\tilde\epsilon - \rho + \gamma') + \gamma'(\tilde\beta + \alpha + \tau') + \Phi_{10}&({\rm g}')\cr
\triangle\epsilon - D\alpha &= -\tau'\rho - \kappa\sigma' - \tilde\kappa\gamma + \beta\tilde\sigma - \alpha(\tilde\epsilon - \rho + \gamma') + \epsilon(\tilde\beta + \alpha + \tau') - \Phi_{10}&({\rm h})\cr
-\delta\epsilon' - D'\alpha' &= -\tau\rho' - \kappa'\sigma - \tilde\kappa'\gamma' + \beta'\tilde\sigma' - \alpha'(\tilde\epsilon' - \rho' + \gamma) + \epsilon'(\tilde\beta' + \alpha' + \tau) + \Phi_{12}&({\rm h}')\cr
D\beta - \delta\epsilon &= \kappa(\rho' + \gamma) + \sigma(\tau' - \alpha) + \beta(\tilde\gamma' - \tilde\rho) - \epsilon(\tilde\tau' + \tilde\alpha) + \Psi_1&({\rm i})\cr
D'\beta' + \triangle\epsilon' &= \kappa'(\rho + \gamma') + \sigma'(\tau - \alpha') + \beta'(\tilde\gamma - \tilde\rho') - \epsilon'(\tilde\tau + \tilde\alpha') - \Psi_3&({\rm i}')\cr
\triangle\gamma - D'\alpha &= \kappa'(\epsilon - \rho) + \sigma'(\beta - \tau) + \alpha(\tilde\rho' - \tilde\gamma) - \gamma(\tilde\tau + \tilde\alpha') - (\gamma\beta' + \alpha\epsilon') + \Psi_3&({\rm j})\cr
-\delta\gamma' - D\alpha' &= \kappa(\epsilon' - \rho') + \sigma(\beta' - \tau') + \alpha'(\tilde\rho - \tilde\gamma') - \gamma'(\tilde\tau' + \tilde\alpha) - (\gamma'\beta + \alpha'\epsilon) - \Psi_1 &({\rm j}')\cr
D\gamma - D'\epsilon &= \tau\tau' - \kappa\kappa' - \beta(\tau' + \tilde\tau) - \alpha(\tilde\tau' + \tau) - \epsilon(\gamma + \tilde\gamma) + \gamma(\gamma' + \tilde\gamma') + \Psi_2 + \Phi_{11} - \Pi &({\rm k})\cr
D'\gamma' - D\epsilon' &= \tau\tau' - \kappa\kappa' - \beta'(\tau + \tilde\tau') - \alpha'(\tilde\tau + \tau') - \epsilon'(\gamma' + \tilde\gamma') + \gamma'(\gamma + \tilde\gamma) + \Psi_2 + \Phi_{11} - \Pi &({\rm k}')\cr
\triangle\beta - \delta\alpha &= \rho\rho' - \sigma\sigma' - \alpha\tilde\alpha - \beta\tilde\alpha' + \alpha(\beta + \alpha') + \gamma(\rho - \tilde\rho) + \epsilon(\tilde\rho' - \rho') + \Psi_2 - \Phi_{11} - \Pi &({\rm l})\cr
-\delta\beta' + \triangle\alpha' &= \rho\rho' - \sigma\sigma' - \alpha'\tilde\alpha' - \beta'\tilde\alpha + \alpha'(\beta' + \alpha) + \gamma'(\rho' - \tilde\rho') + \epsilon'(\tilde\rho - \rho) + \Psi_2 - \Phi_{11} - \Pi &({\rm l}')\cr}$$
Additionally, there are the tilde versions of all these equations, with $\tilde\Psi_{\bf i}$ replacing $\Psi_{\bf i}$ and $\tilde\Phi_{A'B'CD} = \Phi_{CDA'B'}$ replacing $\Phi_{ABC'D'}$, so that
$$\displaylines{\Phi_{00} \qquad \Phi_{11} \qquad \Phi_{22} \qquad\hbox {are unchanged but}\cr
\noalign{\vskip -12pt}
\hfill\llap(3.5)\cr
\noalign{\vskip -12pt}
\Phi_{01}\ \leftrightarrow\ \Phi_{10} \qquad \Phi_{02}\ \leftrightarrow\ \Phi_{20} \qquad \Phi_{12}\ \leftrightarrow\ \Phi_{21},\cr}$$
when applying tilde.

If the dyads are in fact spin frames, by virtue of (2.11.1--2),
$$({\rm h})\ \Leftrightarrow ({\rm g}') \qquad ({\rm h}')\ \Leftrightarrow ({\rm g}) \qquad ({\rm j})\ \Leftrightarrow ({\rm i}') \qquad ({\rm j}')\ \Leftrightarrow ({\rm i}) \qquad ({\rm k})\ \Leftrightarrow ({\rm k}') \qquad ({\rm l})\ \Leftrightarrow ({\rm l}'),\eqno(3.6)$$
which also motivates the grouping of terms  in (j) and $({\rm j}')$ as for spin frames the final bracketed pair vanishes.
\vskip 24pt
\noindent {\section 4. Computation of Spin Coefficients}
\vskip 12pt
Computation of the spin coefficients depends on the information one has to hand, for example:\hfil\break
1) if one has an explicit pair of spinor bases and a tensor basis, one can compute the Infeld-van der Waerden symbols and use the analogue of [22], (4.5.38); this method tends to be tedious;\hfil\break
2) if one knows the Christoffel symbols and the null tetrad in some coordinate system, one can employ (2.7) and (2.8) or (2.10);\hfil\break
3) if one knows a null tetrad, one can compute the spin coefficients with respect to the corresponding spin frames, whence (2.11.1--2) hold, without explicit knowledge of the connection, by using the analogues of [22], (4.13.44). First note that these equations are independent of which of the two customary definitions of the wedge product are employed. (As indicated in [17], Appendix One, I employ that definition of wedge product between $p$ and $q$ forms which is $(p+q)!/p!q!$ times the square bracket operation of [22], \S 3.3.) The neutral analogue of [22], (4.13.44), may be obtained as in [22], \S 4.13 or directly from (2.10) as follows. With $\phi_a$ denoting any member of the null tetrad of one-forms, express $d\phi = 2\nabla_{[a}\phi_{b]}$ formally as a linear combination of the basis of two-forms induced by the null tetrad, and then compute the coefficients of the linear combination using (2.10). The results are as follows, where the second equality in each set of equations is obtained by applying (2.11.1--2) and yields the direct analogue of [22], (4.13.44).
$$\eqalignno{&&(4.1)\cr
d\ell &= (\tilde\tau + \tilde\beta + \alpha)\ell\wedge m + (\tau + \tilde\alpha + \beta)\ell\wedge\tilde m - (\epsilon + \tilde\epsilon)\ell\wedge n + (\tilde\rho - \rho)m\wedge \tilde m - \tilde\kappa m\wedge n - \kappa\tilde m\wedge n\cr
&= (\tilde\tau + \tilde\alpha' + \beta')\ell\wedge m + (\tau + \tilde\beta' + \alpha')\ell\wedge\tilde m + (\gamma' + \tilde\gamma')\ell\wedge n + (\tilde\rho - \rho)m\wedge \tilde m - \tilde\kappa m\wedge n - \kappa\tilde m\wedge n\cr
dm &= (\gamma + \tilde\epsilon'+ \tilde\rho')\ell\wedge m + \tilde\sigma'\ell\wedge\tilde m + (\tau + \tilde\tau')\ell\wedge n + (\beta - \tilde\beta')m\wedge \tilde m - (\rho + \epsilon + \tilde\gamma') m\wedge n - \sigma m\wedge n\cr
&= (\tilde\rho' -\epsilon' - \tilde\gamma)\ell\wedge m + \tilde\sigma'\ell\wedge\tilde m + (\tau + \tilde\tau')\ell\wedge n + (\alpha' - \tilde\alpha)m\wedge \tilde m + (\gamma' + \tilde\epsilon - \rho) m\wedge n - \sigma \tilde m\wedge n\cr
d\tilde m &= \sigma'\ell\wedge m + (\tilde\gamma + \epsilon'+ \rho')\ell\wedge\tilde m + (\tilde\tau + \tau')\ell\wedge n + (\beta' - \tilde\beta)m\wedge \tilde m - \tilde\sigma m \wedge n -(\tilde\rho + \tilde\epsilon + \gamma')\tilde m\wedge n\cr
&= \sigma'\ell\wedge m - (\tilde\epsilon' + \gamma - \rho')\ell\wedge\tilde m + (\tilde\tau + \tau')\ell\wedge n + (\alpha - \tilde\alpha')m\wedge \tilde m - \tilde\sigma m \wedge n + (\epsilon + \tilde\gamma' - \tilde\rho)\tilde m\wedge n\cr
dn &= -\kappa'\ell\wedge m -\tilde\kappa'\ell\wedge\tilde m + (\epsilon' + \tilde\epsilon')\ell\wedge n + (\rho' - \tilde\rho')m\wedge \tilde m + (\tilde\alpha' + \beta' + \tau') m\wedge n + (\alpha' + \tilde\beta' + \tilde\tau')\tilde m\wedge n\cr
&= -\kappa'\ell\wedge m -\tilde\kappa'\ell\wedge\tilde m - (\gamma + \tilde\gamma)\ell\wedge n + (\rho' - \tilde\rho')m\wedge \tilde m + (\tilde\beta + \alpha + \tau') m\wedge n + (\beta + \tilde\alpha + \tilde\tau')\tilde m\wedge n\cr}$$
\vskip 24pt
\noindent {\section 5. Walker Geometry}
\vskip 12pt
Let $(M,g,[\pi^{A'}])$ be a Walker geometry, i.e., a four-dimensional neutral manifold $(M,g)$ for which the projective spinor field $[\pi^{A'}]$ defines a parallel field of $\alpha$-planes, see [17] for all details relating to Walker geometry. This parallel distribution is automatically integrable (the integral surfaces are called $\alpha$-surfaces). The parallel condition imposed on the distribution is just
$$\pi_{A'}\nabla_b\pi^{A'} = 0,\eqno(5.1)$$
for any local scaled representative (LSR) $\pi^{A'}$ of $[\pi^{A'}]$. Walker [32] showed that one can always introduce local coordinates $(u,v,x,y)$ with respect to which the metric $g$ takes a simple canonical form:
$$(g_{\bf ab}) = \pmatrix{{\bf 0}_2&I_2\cr I_2&W\cr},\qquad\hbox{where}\qquad W = \pmatrix{a&c\cr c&b\cr}\qquad \hbox{is arbitrary}.\eqno(5.2)$$
In [17], it was shown that, after possibly invoking the interchange $u \leftrightarrow v$ and $x \leftrightarrow y$, there is a unique (up to sign) LSR $\pi^{A'}$ such that [17] (2.8) holds:
$$\partial_u = \alpha^A\pi^{A'} \qquad \partial_v = \beta^A\pi^{A'} \qquad dx = \alpha_A\pi_{A'} \qquad dy = \beta_A\pi_{A'}$$
where $\{\alpha^A,\beta^A\}$ is an unprimed spin frame. Walker coordinates $(u,v,x,y)$ for which these conditions hold will be said to be {\sl oriented}. For oriented Walker coordinates, the following {\sl Walker null tetrad\/} was introduced in [17] (2.11):
$$\vcenter{\openup1\jot \halign{$\hfil#$&&${}#\hfil$&\qquad$\hfil#$\cr
\ell^a &:= \partial_u = \alpha^A\pi^{A'} & \tilde m^a &:= \partial_v = \beta^A\pi^{A'}\cr
n^a &:= -{a \over 2}\partial_u - {c \over 2}\partial_v + \partial_x = \beta^A\xi^{A'} & m^a &:= {c \over 2}\partial_u + {b \over 2}\partial_v - \partial_y = \alpha^A\xi^{A'},\cr}}\eqno(5.3)$$
where the associated spin frames $\{\alpha^A,\beta^A\}$ and $\{\pi^{A'},\xi^{A'}\}$ are unique up to an overall change of sign. 

The Infeld-van-der-Waerden symbols connecting the coordinate basis and the spin frames are:
$$\vcenter{\openup3\jot \halign{$\hfil#$&&${}#\hfil$&\qquad$\hfil#$\cr
\sigma_1{}^{\bf AA'} &= \pmatrix{1&0\cr 0&0\cr} & \sigma_2{}^{\bf AA'} &= \pmatrix{0&0\cr 1&0\cr}\cr
\sigma_3{}^{\bf AA'} &= \pmatrix{a/2&0\cr c/2&1\cr} & \sigma_4{}^{\bf AA'} &= \pmatrix{c/2&-1\cr b/2&0\cr}\cr
\sigma_{\bf AA'}{}^1 &= \pmatrix{1&c/2\cr 0&-a/2\cr} & \sigma_{\bf AA'}{}^2 &= \pmatrix{0&b/2\cr 1&-c/2\cr}\cr
\sigma_{\bf AA'}{}^3 &= \pmatrix{0&0\cr 0&1\cr} & \sigma_{\bf AA'}{}^4 &= \pmatrix{0&-1\cr 0&0\cr}.\cr}}\eqno(5.4)$$
From the Infeld-van-der-Waerden symbols (5.4) one could compute the spin coefficients using [22], (4.5.38). Since the Christoffel symbols for Walker coordinates are known, [17], A1.1, one can more rapidly compute the covariant derivatives of the null tetrad:
\vskip 6pt
$$\vcenter{\openup3\jot \halign{$\hfil#$&&${}#\hfil$&\qquad$\hfil#$\cr
D\ell^a &= D\tilde m^a = Dn^a = D m^a = 0 & \triangle\ell^a &= \triangle\tilde m^a = \triangle n^a = \triangle m^a = 0;\cr
D'\ell^a &= {1 \over 2}(a_1\ell^a + c_1\tilde m^a) & D'\tilde m^a &= {1 \over 2}(a_2\ell^a + c_2\tilde m^a);\cr
\delta\ell^a &= -{1 \over 2}(c_1\ell^a + b_1\tilde m^a) & \delta\tilde m^a &= -{1 \over 2}(c_2\ell^a + b_2\tilde m^a);\cr}}$$
$$\eqalignno{D'n^a &= {1 \over 4}(2c_3 - 2a_4 + ba_2 + ca_1 - ac_1 - cc_2)\tilde m^a - {a_1 \over 2}n^a + {a_2 \over 2}m^a;&(5.5)\cr
D'm^a &= {1 \over 4}(2c_3 - 2a_4 + ba_2 + ca_1 - ac_1 - cc_2)\ell^a + {c_1 \over 2}n^a - {c_2 \over 2}m^a;\cr
\delta n^a &= {1 \over 4}(2c_4 - 2b_3 - cc_1 - bc_2 + ab_1 + cb_2)\tilde m^a + {c_1 \over 2}n^a - {c_2 \over 2}m^a;\cr
\delta m^a &= {1 \over 4}(2c_4 - 2b_3 - cc_1 - bc_2 + ab_1 + cb_2)\ell^a - {b_1 \over 2}n^a + {b_2 \over 2}m^a.\cr}$$ 
and then compute the spin coefficients using (2.10). Or, using the covariant version of (5.2), one can  compute the spin coefficients from (4.1). One obtains:
\vskip 6pt
$$\vcenter{\openup2\jot \halign{$\hfil#$&&${}#\hfil$&\qquad$\hfil#$\cr
\epsilon &= 0 & \epsilon' &= \displaystyle{c_2 - a_1 \over 4} & \tilde\epsilon &= 0 & \tilde\epsilon' &= \displaystyle-\left({a_1+c_2 \over 4}\right)\cr
\alpha &= 0 & \alpha' &= \displaystyle{b_2 - c_1 \over 4} & \tilde\alpha &= \displaystyle-\left({b_2 + c_1 \over 4}\right) & \tilde\alpha' &= 0\cr
\beta &= \displaystyle{b_2 - c_1 \over 4} & \beta' &= 0 & \tilde\beta &= 0 & \tilde\beta' &= \displaystyle-\left({b_2 + c_1 \over 4}\right)\cr
\gamma &= \displaystyle{a_1 - c_2 \over 4} & \gamma' &= 0 & \tilde\gamma &= \displaystyle{a_1 + c_2 \over 4} & \tilde\gamma' &= 0\cr
\kappa &= 0 & \kappa' &= \displaystyle-{a_2 \over 2} & \tilde\kappa &= 0 & \tilde\kappa' &= \displaystyle{2a_4 - 2c_3 -ba_2 - ca_1 + ac_1 +cc_2 \over 4}\cr
\rho &= 0 & \rho' &= \displaystyle-{c_2 \over 2} & \tilde\rho &= 0 & \tilde\rho' &= 0\cr
\sigma &= \displaystyle-{b_1 \over 2} & \sigma' &= 0 & \tilde\sigma &= 0 & \tilde\sigma' &= \displaystyle{2c_4 - 2b_3 - cc_1 - bc_2 + ab_1 + cb_2 \over 4}\cr
\tau &= \displaystyle{c_1 \over 2} & \tau' &= 0 & \tilde\tau &= 0 & \tilde\tau' &= 0\cr}}\eqno(5.6)$$
Observe that
$$\alpha' + \tilde\alpha + \tau = \beta + \tilde\beta' + \tau = 0 \hskip 1in \gamma - \tilde\gamma - \rho' = \epsilon' - \tilde\epsilon' + \rho' = 0.\eqno(5.7)$$
From (2.9), the spin frames associated to the Walker null tetrad satisfy:
$$D\alpha^A = D\beta^A = 0 = D\pi^{A'} = D\xi^{A'} \hskip 1in \triangle\alpha^A = \triangle\beta^A = 0 =\triangle\pi^{A'} = \triangle\xi^{A'}$$
$$\vcenter{\openup2\jot \halign{$\hfil#$&&${}#\hfil$&\qquad$\hfil#$\cr
D'\alpha^A &= \displaystyle\left({a_1 - c_2 \over 4}\right)\alpha^A + {c_1 \over 2}\beta^A & \delta\alpha^A &= \displaystyle\left({b_2 - c_1 \over 4}\right)\alpha^A - {b_1 \over 2}\beta^A\cr
D'\beta^A &= \displaystyle{a_2 \over 2}\alpha^A + \left({c_2 - a_1 \over 4}\right)\beta^A & \delta\beta^A &= \displaystyle-{c_2 \over 2}\alpha^A - \left({b_2 - c_1 \over 4}\right)\beta^A\cr
D'\pi^{A'} &= \displaystyle{a_1 + c_2 \over 4}\pi^{A'} & \delta\pi^{A'} &= \displaystyle-\left({b_2 + c_1 \over 4}\right)\pi^{A'}\cr}}\eqno(5.8)$$
$$\eqalign{D'\xi^{A'} &= {2c_3 - 2a_4 + ba_2 + ca_1 - ac_1 - cc_2 \over 4}\pi^{A'} - \left({a_1 + c_2 \over 4}\right)\xi^{A'}\cr
\delta\xi^{A'} &= {2c_4 - 2b_3 - cc_1 - bc_2 + ab_1 + cb_2 \over 4}\pi^{A'} + \left({b_2 + c_1 \over 4}\right)\xi^{A'}\cr}$$
In particular: in accordance with [17], (2.3),
$$\nabla_b\pi^{A'} = P_b\pi^{A'}\qquad\hbox{where}\qquad P_b := {a_1 + c_2 \over 4}\ell_b + {b_2 + c_1 \over 4}\tilde m_b;\eqno(5.9)$$
the condition for $\pi^{A'}$ to be parallel agrees with [17], (3.2); the condition for $\xi^{A'}$ to be parallel agrees with [17], (4.26--28); and if $\pi^{A'}$ is parallel and $\xi_{A'}\nabla_b\xi^{A'} = 0$ ($\tilde\kappa' = \tilde\sigma' = 0$; in particular, for a double Walker geometry, see [17] \S 4) then $\xi^{A'}$ is in fact parallel too, i.e, the primed spin frame is parallel.

The commutators (3.1) are
$$\vcenter{\openup2\jot \halign{$\hfil#$&&${}#\hfil$&\qquad$\hfil#$\cr
[D',D] &= \displaystyle{a_1 \over 2}D + {c_1 \over 2}\triangle & [\delta,D] &= \displaystyle-{c_1 \over 2}D - {b_1 \over 2}\triangle & [\delta,D'] &= \tilde\kappa' D + \tilde\sigma'\triangle\cr
[D',\triangle] &= \displaystyle{a_2 \over 2}D + {c_2 \over 2}\triangle & [D,\triangle] &= 0 & [\triangle,\delta] &= \displaystyle{c_2 \over 2}D + {b_2 \over 2}\triangle.\cr}}\eqno(5.10)$$
The equation $[D,\triangle] = 0$ of course expresses the integrability of the Walker distribution. In fact, since the Walker spin frames are parallel with respect to $D$ and $\triangle$, $D$ and $\triangle$ commute on arbitrary spinors. In the case of a double Walker geometry ([17], \S 4), in which $\xi^{A'}$ is an LSR for the projective spinor field defining the complementary parallel distribution, one sees that integrability of that complementary distribution, $[\delta,D'] = 0$, is given, locally, by $\tilde\kappa' = \tilde\sigma' = 0$, which conditions, as already noted, are precisely the conditions for $\xi_{A'}\nabla_b\xi^{A'} = 0$, i.e., the local condition for the complementary distribution to be parallel.

Substituting the expressions (5.6) for the spin coefficients into (3.4) confirms the curvature calculations of [17]. Setting the vanishing spin coefficients of (5.6) to zero in (3.4) yields the formulae:
$$\Psi_0 = -D\sigma\qquad ({\rm b}) \hskip 1.25in \Psi_4 = -\triangle\kappa'\qquad ({\rm b}')$$
$$\eqalign{\Psi_1 &= D\beta\qquad ({\rm i})/({\rm j}')\cr
&= -\left({D\tau + \triangle\sigma \over 2}\right)\qquad ({\rm c})\& ({\rm d})\cr} \hskip 1in \eqalign{\Psi_3 &= \triangle\gamma\qquad ({\rm j})/({\rm i}')\cr
&= {\triangle\rho' - D\kappa' \over 2}\qquad ({\rm c}')\&({\rm d}')\cr}$$
$$\Psi_2 = {D\gamma + \triangle(\beta - \tau) \over 3} = {D(\gamma + \rho') + \triangle\beta \over 3}\qquad ({\rm f})/({\rm f}')\&({\rm k})/({\rm k}')\&({\rm l})/({\rm l}').$$
$$\tilde\Psi_0 = 0\qquad \widetilde{({\rm b})} \hskip 1.25in \tilde\Psi_1 = 0\qquad \widetilde{({\rm i})}/\widetilde{({\rm j}')}\ \hbox{or}\ \widetilde{({\rm c})}\&\widetilde{({\rm d})}$$
$$\tilde\Psi_2 = {D\tilde\gamma - \triangle\tilde\alpha \over 3} = {S \over 12}\qquad \widetilde{({\rm f})}/\widetilde{({\rm f}')}\&\widetilde{({\rm k})}/\widetilde{({\rm k}')}\&\widetilde{({\rm l})}/\widetilde{({\rm l}')}$$
$$\eqalignno{\tilde\Psi_3 &= \delta\tilde\gamma - D'\tilde\alpha\qquad \widetilde{({\rm j})}/\widetilde{({\rm i}')}&(5.11)\cr
&= -\left({D\tilde\kappa' + \triangle\tilde\sigma' \over 2}\right)\qquad \widetilde{({\rm c}')}\&\widetilde{({\rm d}')}\cr}$$
$$\tilde\Psi_4 = 2(\tilde\sigma'\tilde\epsilon' - \tilde\kappa'\tilde\beta') - \delta\tilde\kappa' - D'\tilde\sigma'\qquad \widetilde{({\rm b}')}.$$
$$S = 4\bigl(D\gamma + \triangle(\beta+2\tau)\bigr) = 4\bigl(D(\gamma - 2\rho') + \triangle\beta\bigr)\qquad ({\rm f})/({\rm f}')\&({\rm k})/({\rm k}')\&({\rm l})/({\rm l}').$$
$$\Phi_{00} = 0\qquad ({\rm a})/\widetilde{({\rm a})} \hskip .75in \Phi_{10} = 0\qquad ({\rm g}')/({\rm h})/\widetilde{({\rm c})}\&\widetilde{({\rm d})} \hskip .75in \Phi_{20} = 0\qquad ({\rm e}')/\widetilde{({\rm e})}$$
$$\eqalign{\Phi_{22} &= -\triangle\tilde\kappa'\qquad \widetilde{({\rm a}')}\cr
&= 2(\rho'\epsilon' - \kappa'\alpha') -\delta\kappa' - D'\rho' \qquad ({\rm a}')\cr} \hskip .4in \eqalign{\Phi_{02} &= D\tilde\sigma'\qquad \widetilde{({\rm e}')}\cr
&= D'\sigma - \delta\tau + 2(\tau\beta - \sigma\gamma)\qquad ({\rm e})\cr}$$
$$\eqalign{\Phi_{01} &= D\tilde\beta'\qquad \widetilde{({\rm g}')}/\widetilde{({\rm h})}\cr &={\triangle\sigma - D\tau\over 2}\qquad ({\rm c}) \& ({\rm d})\cr} \hskip 1in \eqalign{\Phi_{21} &= \triangle\tilde\gamma\qquad \widetilde{({\rm g})}/\widetilde{({\rm h}')}\cr &= -\left({D\kappa' + \triangle\rho' \over 2}\right)\qquad ({\rm c}')\&({\rm d}')\cr}$$
$$\eqalign{\Phi_{11} &= {D\gamma - \triangle\beta \over 2}\qquad ({\rm k})/({\rm k}') \& ({\rm l})({\rm l}')\cr &= {D\tilde\gamma + \triangle\tilde\alpha \over 2}\qquad \widetilde{({\rm k})}/\widetilde{({\rm k}')}\&\widetilde{({\rm l})}/\widetilde{({\rm l}')}\cr}\hskip 1in \eqalign{ \Phi_{12} &= \tau\rho'+\kappa'\sigma-D'\beta+\delta\gamma\qquad ({\rm g})/({\rm h}')\cr &= {\triangle\tilde\sigma' - D\tilde\kappa' \over 2}\qquad \widetilde{({\rm c}')}\&\widetilde{({\rm d}')}\cr}$$
\vskip 24pt
\noindent {\section 6. Null Geometry}
\vskip 12pt
The material in this section is heavily influenced by [23], \S\S 7.1--3. In order to explore the utility and interpretation of the spin coefficients, consider the following notions of {\sl null distribution\/} on $(M,g)$:\hfil\break
\noindent type I, a distribution of null one-dimensional subspaces;\hfil\break
\noindent type II, a distribution of totally null two-dimensional subspaces;\hfil\break
\noindent type III, a three-dimensional distribution $\cal H$ for which the distribution ${\cal H}^\perp$ is a distribution of type I.

Note that type II are self orthogonal. To motivate the definition of null distributions of type III, let $\cal H$ be a three-dimensional distribution on which $g$ is degenerate. At any point $p \in M$, ${\cal D} := {\cal H} \cap {\cal H}^\perp$ is totally null and thus at most two dimensional by virtue of the signature of $g$; if two-dimensional, $\cal D$ must be self orthogonal but is also orthogonal to any vector in $\cal H$ linearly independent of $\cal D$ itself. Thus, ${\cal D}$ is one dimensional, and ${\cal D} = {\cal H}^\perp$ (by dimensional considerations).

Given any null vector $\ell$ in ${\bf R}^{2,2}$, the construction of a Witt decomposition (see [29], [4], \S 1.3) provides another null vector $n$ such that $s(\ell,n) = 1$ and a decomposition ${\bf R}^{2,2} = \langle \ell \rangle_{\bf R} \oplus \langle n \rangle_{\bf R} \oplus Z$, where $Z = \left(\langle \ell \rangle_{\bf R} \oplus \langle n \rangle_{\bf R}\right)^\perp \cong {\bf R}^{1,1}$. Hence, one can choose linearly independent null vectors $m$, $\tilde m \in Z$ satisfying $s(m,\tilde m) = -1$, i.e., one can construct a null tetrad $\{\ell,n,m,\tilde m\}$ from $\ell$ so that $\langle \ell \rangle_{\bf R}^\perp = \langle \ell,m,\tilde m \rangle_{\bf R} = \langle \ell \rangle_{\bf R} \operp Z$.
\vskip 12pt
\noindent {\bf 6.1. Null Distributions of Type I}
\vskip 12pt
To begin, let $\cal D$ be a null distribution of type I. By one-dimensionality, $\cal D$ is automatically integrable and one can choose Frobenius coordinates $(u,x,y,z)$ with integral manifolds specified by constant values of $x$, $y$, $z$, with $u$ parameterizing the integral manifolds, and $\partial_u$ a local section of the distribution. I shall refer to the integral manifolds of $\cal D$ as {\sl integral curves\/}, though of course they are integral curves in the technical sense only with respect to a specific local section of $\cal D$.

Let $(v,r,s,t)$ be another Frobenius chart for $\cal D$ with overlapping domain with that of $(u,x,y,z)$. On their common domain the Jacobian for $(v,r,s,t) = \Phi(u,x,y,z)$ takes the form
$$J = \pmatrix{{\partial v \over \partial u}&{\partial v \over \partial x}&{\partial v \over \partial y}&{\partial v \over \partial z}\cr 0&&&\cr 0&&*&\cr 0&&&\cr}\eqno(6.1.1)$$
with $\partial v/\partial u \not= 0$ (for invertibility of $J$). In particular,
$$\partial_u = {\partial v \over \partial u}\,\partial_v;\qquad\hbox{equivalently } \partial_v = {\partial u \over \partial v}\,\partial_u,\eqno(6.1.2)$$
i.e., the local section of $\cal D$ defined by $(v,r,s,t)$ is, of course, a smooth non-vanishing multiple of that provided by $(u,x,y,z)$. One may regard reparametrization of the integral curves, i.e., replacing $u$ by $f(u)$, for some smooth invertible $f$, as the special case of (6.1.1--2) with $v = f(u)$ and $(r,s,t)=(x,y,z)$.

Given a local section $\ell^a$ of $\cal D$, one can write, at least locally, $\ell^a = o^Ao^{A'}$, for some spinors $o^A \in S_U$, $o^{A'} \in S'_U$, on some suitable open subset $U$ of $M$. Locally, one can then construct smooth spin frames $\{o^A,\iota^A\}$ and $\{o^{A'},\iota^{A'}\}$ and from these a local null tetrad. In relation to $\cal D$, the freedom in this construction is
$$o^A \mapsto \lambda o^A \hskip .5in o^{A'} \mapsto \tilde\lambda o^{A'} \hskip .5in \iota^A \mapsto \lambda^{-1}\iota^A + \mu o^A \hskip .5in \iota^{A'} \mapsto \tilde\lambda^{-1}\iota^{A'} + \tilde\mu o^{A'},\eqno(6.1.3)$$
where $\lambda$, $\tilde\lambda \in {\bf R}^*$ and $\mu$, $\tilde\mu \in {\bf R}$, i.e., the freedom of
$$(D,\tilde D) := \left(\pmatrix{\lambda&\mu\cr 0&\lambda^{-1}\cr},\pmatrix{\tilde\lambda&\tilde\mu\cr 0&\tilde\lambda^{-1}\cr}\right) \in {\bf SL(2;R)} \times {\bf SL(2;R)}.\eqno(6.1.4)$$
Under (6.1.3--4), the associated null tetrad transforms as:
$$\vcenter{\openup1\jot \halign{$\hfil#$&&${}#\hfil$&\qquad$\hfil#$\cr
\ell^a &\mapsto \lambda\tilde\lambda\ell^a & n^a &\mapsto (\lambda\tilde\lambda)^{-1}n^a + \lambda^{-1}\tilde\mu\tilde m^a + \mu\tilde\lambda^{-1} m^a + \mu\tilde\mu\ell^a\cr
m^a &\mapsto \lambda\tilde\lambda^{-1} m^a + \lambda\tilde\mu \ell^a & \tilde m^a &\mapsto \lambda^{-1}\tilde\lambda\tilde m^a + \mu\tilde\lambda \ell^a.\cr}}\eqno(6.1.5)$$
The spin coefficients are not invariant under (6.1.3--4). The simplest transformations are for those spin coefficients of the form $o_A \diamondsuit o^A$ and $o_{A'}\diamondsuit o^{A'}$, where $\diamondsuit = D$, $D'$, $\delta$ or $\triangle$. One computes
$$\vcenter{\openup1\jot \halign{$\hfil#$&&${}#\hfil$&\qquad$\hfil#$\cr
\kappa &\mapsto \lambda^3\tilde\lambda\kappa & \tilde\kappa &\mapsto \tilde\lambda^3\lambda\tilde\kappa\cr
\rho &\mapsto \lambda\tilde\lambda\rho + \lambda^2\tilde\lambda\mu\kappa & \tilde\rho &\mapsto \tilde\lambda\lambda\tilde\rho + \tilde\lambda^2\lambda\tilde\mu\tilde\kappa\cr
\sigma &\mapsto \lambda^3\tilde\lambda^{-1}\sigma + \lambda^3\tilde\mu\kappa & \tilde\sigma &\mapsto \tilde\lambda^3\lambda^{-1}\tilde\sigma + \tilde\lambda^3\mu\tilde\kappa\cr
\tau &\mapsto \lambda\tilde\lambda^{-1}\tau + \lambda\tilde\mu\rho + \lambda^2\tilde\lambda^{-1}\mu\sigma + \lambda^2\mu\tilde\mu\kappa & \tilde\tau &\mapsto \tilde\lambda\lambda^{-1}\tilde\tau + \tilde\lambda\mu\tilde\rho + \tilde\lambda^2\lambda^{-1}\tilde\mu\tilde\sigma + \tilde\lambda^2\tilde\mu\mu\tilde\kappa;\cr}}\eqno(6.1.6)$$
which can be written more compactly in the form
$$\pmatrix{\kappa\cr \sigma\cr} \mapsto \lambda^3\,({^\tau\!{\tilde D}})\pmatrix{\kappa\cr \sigma\cr} \hskip 1in \pmatrix{\rho\cr \tau\cr} \mapsto \lambda\,({^\tau\!{\tilde D}})\pmatrix{\rho\cr \tau\cr} + \lambda^2\mu\,({^\tau\!{\tilde D}})\pmatrix{\kappa\cr \sigma,\cr}$$
and their tilde versions, where $^\tau\!\tilde D$ denotes the transpose of $\tilde D$.
\vskip 24pt
\noindent {\bf 6.1.7 Definition}\hfil\break
If $N$ is a submanifold of $M$ for which the induced metric $h:= g\vert_N$ is degenerate then the Koszul construction of a Levi-Civita connection from $h$ fails (but cf. [15]). Moreover, $g$ fails to define a preferred complementary bundle to $TN$ within $TM$ (as $TN^\perp \cap TN \not= \langle {\bf 0} \rangle_{\bf R}$), so there is no prefered projection of $\nabla_X Y$ from $TM$ to $TN$, for (local) sections $X$, $Y$ of $TN$. One can still define, however, the shape operator (i.e., second fundamental form) $\shape$ as follows: for (local) smooth sections $X$ and $Y$ of $TN$,
$$\shape(X,Y) := [\nabla_X Y] \in TM/TN.$$
One easily checks that $\shape$ is bilinear over $TN$ as a ${\cal C}^\infty(N)$-module and symmetric. 

When $\shape$ is zero, then $\nabla_XY$ lies within $TN$ and thus defines a linear connection on $TN$. For any curve within $N$, parallel propagation along that curve is unambiguous; in particular, the notion of the geodesic through $p \in N$ in direction $X \in T_pN$ is also unambiguous and lies within $N$ (at least in a neighbourhood of $p$). In this sense, $N$ is totally geodesic. Hence, one can, whether $h$ is nondegenerate or degenerate, define $N$ to be {\sl auto-parallel\/} if $\shape$ vanishes, which ([31]) is equivalent to: parallel propagation of $T_pN$ along curves in $N$ remains in $TN$, whence the definition is consistent with [13], \S VII.8.

A distribution $\cal D$ will be called {\sl auto-parallel\/} if invariant under parallel transport along any curve tangent to $\cal D$. By the result in [31] mentioned in the previous paragraph, a distribution $\cal D$ is auto-parallel iff $\nabla_YX \in {\cal D}$ for all local sections $X$, $Y$ of $\cal D$. It follows immediately that when $\cal D$ is auto-parallel, $[X,Y] = \nabla_XY - \nabla_YX \in {\cal D}$, whence $\cal D$ is integrable. The integral manifolds of $\cal D$ are auto-parallel as defined in the previous paragraph. One can also define the shape operator $\shape$ of $\cal D$ by
$$\shape:{\cal D} \times {\cal D} \to {TM \over {\cal D}} \hskip 1in (X,Y) \mapsto \left[\nabla_XY\right]$$
whence $\cal D$ is auto-parallel iff $\shape$ vanishes.

A null distribution $\cal D$ of type I is auto-parallel iff, for any local section $\ell^a$,
$$\ell^b\nabla_b\ell^a \propto \ell^a.\eqno(6.1.8)$$
Some authors, e.g., [23], regard a curve whose tangent vector satisfies (6.1.8) as a geodesic, while others, e.g., [21], refer to such curves as {\sl pre-geodesic}, restricting the geodesic condition to those parametrized curves for which the proportionality constant in (6.1.8) is zero. I shall employ the more restricted usage of {\sl geodesic}, though for null geodesics the less restrictive notion is perhaps more appropriate since it is invariant under conformal rescalings of the metric.

With $\ell^a = o^Ao^{A'}$, (6.1.8) is equivalent to
$$o_ADo^A = 0\ \hbox{ and }\ o_{A'}Do^{A'} = 0,\ \hbox{i.e., to}\ \ \kappa = \tilde\kappa = 0.\eqno(6.1.9)$$
(6.1.6) confirms that the vanishing of each of $\kappa$ and $\tilde\kappa$ is indeed invariant under (6.1.3--4) and that then $\rho$, $\tilde\rho$, $\sigma$, and $\tilde\sigma$ merely scale under (6.1.3--4). It follows from (2.9), that (6.1.9), whence (6.1.8), is equivalent to each of:
$$\displaylines{Do^A = \epsilon o^A\ \hbox{ and }\ Do^{A'} = \tilde\epsilon o^{A'};\cr
\hfill o_Ao^{B'}\nabla_{BB'}o^A = \rho o_B\ \hbox{ and }\ o_{A'}o^B\nabla_{BB'}o^{A'} = \tilde\rho o_{B'};\hfill\llap(6.1.10)\cr
o_Ao^B\nabla_{BB'}o^A = \sigma o_{B'}\ \hbox{ and }\ o_{A'}o^{B'}\nabla_{BB'}o^{A'} = \tilde\sigma o_B.\cr}$$
Thus, in the auto-parallel case, (6.1.10) provides definitions of the spin coefficients $\epsilon$, $\tilde\epsilon$, $\rho$, $\tilde\rho$, $\sigma$, and $\tilde\sigma$ purely in terms of $o^A$ and $o^{A'}$.

From (2.10), in general
$$D\ell^a = (\epsilon + \tilde\epsilon)\ell^a + \tilde\kappa m^a + \kappa\tilde m^a.\eqno(6.1.11)$$
For any pre-geodesic curve, with parameter $u$ say, as is well known it is always possible to reparametrize the curve (say $v = g(u)$) so that the geodesic condition $\ell^b\nabla_b\ell^a = 0$ holds for $\ell^a = \partial_v$. Indeed, for Frobenius coordinates $(u,x,y,z)$ and $\ell^a = \partial_u =: \partial_1$, $\ell^b\nabla_b\ell^a = \Gamma^i_{11}\partial_i$, so the auto-parallel condition is $\Gamma^i_{11} = 0$, $i=2$, 3, 4, i.e., $\ell^b\nabla_b\ell^a = \Gamma^1_{11}\ell^a$. Reparametrization to achieve affine parametrization is given by
$$v = \int^u\exp\left(\int^t\Gamma^1_{11}ds\right)dt$$
along integral curves. Using the hypersurface $u=0$ as a common initial condition, as $\Gamma^1_{11}$ is smooth on the domain of the Frobenius coordinates, one can construct $v$ to be a smooth function on $U$ which restricts to each integral curve as an affine parameter and, in particular, as a strictly increasing function of $u$ (but as $\Gamma^1_{11}$ is a function of $(u,x,y,z)$, then so is $v$; i.e., the level surfaces of $v$ do not coincide with those of $u$, other than $v=0$ and $u=0$ by construction). As
$${\partial (v,x,y,z) \over \partial(u,x,y,z)} = \pmatrix{{\partial v \over \partial u}&{\partial v \over \partial x}&{\partial v \over \partial y}&{\partial v \over \partial z}\cr 0&1&0&0\cr 0&0&1&0\cr 0&0&0&1\cr},$$
is nonsingular ($v=g(u)$ is bijective), $(v,x,y,z)$ are legitimate coordinates on $U$; indeed, Frobenius coordinates for $\cal D$ on $U$.

Hence, without loss of generality, if $\cal D$ is auto-parallel, one can choose, locally, Frobenius coordinates $(v,x,y,z)$ with $v$ an affine parameter for each integral manifold in the domain of the Frobenius chart.

Note that in general, for auto-parallel $\cal D$ ($\kappa = \tilde\kappa = 0$) with $\ell^a = o^Ao^{A'} = \partial_u$, $Do^A = 0$ ($\kappa = \epsilon = 0$) and $Do^{A'} = 0$ ($\tilde\kappa = \tilde\epsilon = 0$) are independent conditions. But with affine parametrization $\ell^a = o^Ao^{A'} = \partial_v$ ($\kappa = \tilde\kappa = 0 = \epsilon + \tilde\epsilon$), each of $Do^A = 0$ and $Do^{A'} = 0$ implies the other, which is also obvious from
$$0 = D\ell^a = (Do^A)o^{A'} + o^A(Do^{A'}).$$
Moreover, if one can scale $o^A$ by $\lambda$ so that $Do^A = 0$, and one scales $o^{A'}$ by $\lambda^{-1}$ in order to leave $\ell^a$ unchanged and equal to $\partial_v$, then the previous equation entails that $Do^{A'} = 0$ too.

Explicitly, when $\kappa = \tilde\kappa = 0$ in (6.1.11), affine parametrization of integral curves forces $\epsilon + \tilde\epsilon = 0$. Under (6.1.3--4), $\epsilon \mapsto \lambda\tilde\lambda\epsilon + \tilde\lambda D\lambda - \lambda^2\tilde\lambda\mu\kappa$. So, with $\kappa = 0$, choosing $\lambda = \exp(-\int\epsilon)$ (where the integral is along integral curves of $\ell^a$) yields a spin frame for which $\epsilon = 0$ and $Do^A = 0$. One can similarly rescale $o^{A'}$ provided $\tilde\kappa = 0$ and thereby, under the freedom of (6.1.3--4), achieve spin frames so that both $\epsilon$ and $\tilde\epsilon$ vanish. Since under these rescalings each of $o^A$ and $o^{A'}$ are $D$-parallel, $\ell^a$ is $D$-parallel; in fact, since $\epsilon + \tilde\epsilon = 0$, $\lambda = \exp\int\epsilon = \exp\int(-\tilde\epsilon) = \tilde\lambda^{-1}$, whence $\lambda\tilde\lambda = 1$, and $o^Ao^{A'} = \ell^a = \partial_v$ for the new spin frames and the given affine parameter $v$. 

In summary, for auto-parallel $\cal D$, one can reparametrize the integral curves to be geodesics (with affine parameter $v$ also a Frobenius coordinate) and then, under the freedom of (6.1.3--4), choose $o^A$ and $o^{A'}$ parallel along the integral curves so that $o^Ao^{A'} = \ell^a  = \partial_v$. Given $v$, these choices restrict the freedom in (6.1.3--4) to $\tilde\lambda = \lambda^{-1}$ and both constant along integral curves (the freedom in $v$ is that of affine reparametrization, which further allows $\tilde\lambda = a\lambda$, $a$ constant).

One can exploit the remaining freedom in (6.1.3--4) to impose further specializations on the description of an auto-parallel null $\cal D$ of type I. Since $\epsilon = \tilde\epsilon = 0$, then, by (2.9 \& 11), $D\iota^A = 0\ \Leftrightarrow\ \tau' = 0$ and $D\iota^{A'} = 0\ \Leftrightarrow\ \tilde\tau' = 0$. Restrict given $\iota^A$ and $\iota^{A'}$ to the hypersurface $v = 0$, and then parallely propagate them off the hypersurface along the integral curves of $\cal D$ so they are $D$-constant (equivalently; for given $\iota^A$, $D\iota^A = -\tau'o^A$; one is free to replace $\iota^A$ by $\iota^A + \mu o^A$, whence $D(\iota^A + \mu o^A) = (-\tau' + D\mu)o^A$; now choose $\mu$ to solve $D\mu = \tau'$ along integral curves with $\mu = 0$ on the hypersurface $v=0$). Of course, the null tetrad is now parallel along the integral curves and $\mu$ and $\tilde\mu$ in (6.1.3--4) are restricted to be constant along integral curves.

It will prove convenient to refer to the choice of Frobenius coordinates $(v,x,y,z)$ with $v$ an affine parameter for the integral curves of $\cal D$ and with spin frames parallel along integral curves and satisfying $\ell^a := o^Ao^{A'} = \partial_v$ as the {\sl affine-parallel\/} (A-P) condition. Note that this condition is just a convenient description of any auto-parallel type I null distribution, not a restriction on such.
\vskip 24pt
\noindent {\bf 6.1.12 Definition}\hfil\break
A distribution $\cal D$ is called {\sl parallel\/} if invariant under parallel translation ([31]). From [31], $\cal D$ is parallel iff $\nabla_YX \in {\cal D}$ for all local sections $X$ of $\cal D$ and local vector fields $Y$. A null distribution of type I is parallel iff $\nabla_Y\ell \propto \ell$, for any local section $\ell^a$ of $\cal D$. If one has completed $\ell^a$ to a null tetrad, by (2.10) the condition for $\cal D$ to be parallel is
$$\kappa = \tilde\kappa = \sigma = \tilde\sigma = \rho = \tilde\rho = \tau = \tilde\tau = 0.\eqno(6.1.12{\rm a})$$
By (6.1.6), the vanishing of all these quantities is indeed invariant under (6.1.3--4). Any specializations imposed in the auto-parallel case can also be imposed here. Note that Walker [32] provides a convenient form for the metric when $(M,g)$ admits a parallel null distribution of type I: it is possible to choose Frobenius coordinates $(v,x,y,z)$ for which ${\cal D} = \langle \partial_v \rangle_{\bf R}$ and with respect to which the metric has components:
$$(g_{\rm ab}) = \pmatrix{0&0&0&1\cr 0&a&b&d\cr 0&b&c&e\cr 1&d&e&f\cr},\eqno(6.1.12{\rm b})$$
with $a,\ldots,f$ smooth functions (noting that $\nabla_1\partial_1 = \Gamma^k_{11}\partial_k$, then $v$ is affine iff $\Gamma^k_{11} = 0$, each $k$; writing out the Lagrangian for this metric form, one readily checks that no geodesic equation has a $(\dot v)^2$ term, i.e., $\Gamma^k_{11}$ is indeed zero, whence this metric form entails $v$ is indeed an affine parameter) and $a$, $b$, $c$, $d$, and $e$ independent of $v$.
\vskip 24pt
Return to the case of an arbitrary null distribution $\cal D$ of type I. ${\cal D}^\perp$ is  a three-dimensional distribution. Equivalently, ${\cal D}^\perp$ is a rank three subbundle of $TM$ which contains the rank one subbundle $\cal D$. Denote the quotient bundle ${\cal D}^\perp/{\cal D}$ by $\cal Q$ and let $\pi:{\cal D}^\perp \to {\cal Q}$ be the projection. At any point $p \in M$, ${\cal D}^\perp_p = \langle \ell,m,\tilde m \rangle_{\bf R}$, ${\cal Q}_p$ can be represented by ${\cal M}_p := \langle m^a,\tilde m^a \rangle_{\bf R} \leq T_pM$ with $\pi_p(a\ell^a + bm^a + c\tilde m^a)$ represented by $bm^a + c\tilde m^a \in {\cal M}_p$.

From (6.1.11), $D\ell^a \in {\cal D}^\perp$ and $\tilde\kappa m^a + \kappa\tilde m^a$ is a measure of the extent to which $\cal D$ fails to be auto-parallel. Since $\cal D$ only defines $\ell^a$ up to scale, consider the projection of $D\ell^a$ in the quotient $\cal Q$, which may be taken as a `curvature' of a given integral curve. Representing the projection by
$$k^a := \tilde\kappa m^a + \kappa\tilde m^a,\hbox{ note } k^ak_a = -2\kappa\tilde\kappa = (D\ell^a)(D\ell_a).\eqno(6.1.13)$$
By (6.1.5--6), the projection is not invariant under (6.1.3--4) unless one imposes $(\lambda\tilde\lambda)^2 = 1$, i.e., $\tilde\lambda = \pm\lambda^{-1}$.

If $\ell^a$ is a local section of $\cal D$ and $V^a$ a nonvanishing vector field (on the same open domain) nowhere proportional to $\ell^a$, then $\Lie_\ell V^a = 0$ is equivalent to the commutativity of the induced local one-parameter groups of diffeomorphisms $\phi^t$ and $\psi^s$ constituting the flows of $\ell^a$ and $V^a$, respectively. Thus, for any integral curve $C$ of $\ell^a$, the parameter change along $C$ between two integral curves of $V^a$ passing through $C$ is the same as between those two integral curves of $V^a$ along any other integral curve of $\ell^a$ also passing through those two integral curves of $V^a$. Any such $V^a$ is called a {\sl connecting vector field\/} for $\cal D$. In terms of Frobenius coordinates for $\cal D$, a connecting vector field $V^a$ has constant components along a specific integral curve and thus `connects' that integral curve to another in these coordinates.

It suffices for subsequent purposes to suppose $V^a$ is defined along $C_0$ satisfying $\Lie_{\ell}V^a = 0$. Now, $0 = [\ell,V] = \nabla_{\ell}V - \nabla_V\ell$, i.e.,
$$DV^a = V^b\nabla_b\ell^a.\eqno(6.1.14)$$
As $0 = \ell_a\ell^a$, then $0 = \ell_a\nabla_b\ell^a$ and consequently
$$\ell_aDV^a = 0,\qquad\hbox{i.e.,}\qquad DV^a \perp \ell^a.\eqno(6.1.15)$$
Writing
$$V^a = \eta\ell^a + \zeta\tilde m^a + \tilde\zeta m^a + \nu n^a,\eqno(6.1.16)$$
under (6.1.3--4), one computes
$$\vcenter{\openup1\jot \halign{$\hfil#$&&${}#\hfil$&\qquad$\hfil#$\cr
\eta &\mapsto \eta(\lambda\tilde\lambda)^{-1} - \zeta\mu\tilde\lambda^{-1} - \tilde\zeta\tilde\mu\lambda^{-1} + \nu\mu\tilde\mu & \nu &\mapsto \nu\lambda\tilde\lambda\cr
\zeta &\mapsto \lambda\tilde\lambda^{-1}\zeta - \nu\lambda\tilde\mu & \tilde\zeta &\mapsto \tilde\lambda\lambda^{-1}\tilde\zeta - \nu\tilde\lambda\mu.\cr}}\eqno(6.1.17)$$
Taking components of (6.1.14) with respect to the null tetrad yields:
$$\vcenter{\openup1\jot \halign{$\hfil#$&&${}#\hfil$&\qquad$\hfil#$\cr
D\nu &= -\zeta\tilde\kappa - \tilde\zeta\kappa - \nu(\gamma' + \tilde\gamma') & D\eta &= \zeta(\alpha + \tilde\beta + \tau') + \tilde\zeta(\tilde\alpha + \beta + \tilde\tau') + \nu(\gamma + \tilde\gamma)\cr
D\zeta &= \zeta(\rho + \epsilon - \tilde\epsilon) + \tilde\zeta\sigma + \nu(\tau + \tilde\tau') & D\tilde\zeta &= \zeta\tilde\sigma  + \tilde\zeta(\tilde\rho + \tilde\epsilon - \epsilon) + \nu(\tilde\tau + \tau').\cr}}\eqno(6.1.18)$$
These equations can be simplified by using the freedom in (6.1.3) to choose $\iota^A$ and $\iota^{A'}$ so as to satisfy $D\iota^A = 0$ and $D \iota^{A'} = 0$ (by generalizing the argument which ensures this condition in the establishment of the A-P condition), i.e., by (2.9 \& 11), $\epsilon = \tau' = \tilde\epsilon = \tilde\tau' = 0$. With these choices, (6.1.18) becomes
$$D\pmatrix{\eta\cr \zeta\cr \noalign{\smallskip}\tilde\zeta\cr \nu\cr} = \pmatrix{0&\alpha+\tilde\beta&\tilde\alpha+\beta&\gamma+\tilde\gamma\cr 0&\rho&\sigma&\tau\cr 0&\tilde\sigma&\tilde\rho&\tilde\tau\cr \noalign{\smallskip}0&-\tilde\kappa&-\kappa&0\cr}\pmatrix{\eta\cr \zeta\cr \noalign{\smallskip}\tilde\zeta\cr \nu\cr} =: M_{\cal D}\pmatrix{\eta\cr \zeta\cr \noalign{\smallskip}\tilde\zeta\cr \nu\cr}.\eqno(6.1.19)$$
Observe that if $\nu$ is zero, a pair of coupled equations for $\zeta$ and $\tilde\zeta$ results involving only $\sigma$, $\rho$, $\tilde\sigma$ and $\tilde\rho$. But, for nonzero $\zeta$ and $\tilde\zeta$, $\nu$ cannot remain zero unless both $\kappa$ and $\tilde\kappa$ vanish, which is precisely the condition for $\cal D$ to be auto-parallel. The condition $\nu \equiv 0$ is equivalent to $V^a\ell_a = 0$; (6.1.19) indicates that, in general, $V^a$ will not remain orthogonal to $\ell^a$ along $C_0$ even if it is so at some point.

Suppose then that $\cal D$ is indeed auto-parallel. Choose $\iota^A$ and $\iota^{A'}$ $D$ parallel so that (6.1.19) is valid. Then $\nu \equiv $ constant solves (6.1.19); in particular, one can choose connecting vector fields $V^a$ satisfying $V^a \in {\cal D}^\perp$, which I shall call {\sl orthogonal} ({\sl abreast\/} in the Lorentzian context). This condition is obviously a geometric condition on $V^a$ imposed by $\cal D$. At any point $p$ on $C_0$,
$$q^a := \zeta\tilde m^a + \tilde\zeta m^a \in {\cal M}_p, \eqno(6.1.20)$$
represents $\pi_p(V^a)$ and this representation is invariant under (6.1.3--4), see (6.1.17) and (6.1.5). The scalar product $g$ on $T_pM$ induces a scalar product $h_p$ on ${\cal Q}_p$, whence $\pi_p$ preserves scalar products, and a scalar product $s_p$ on ${\cal M}_p$, with
$$({\cal Q}_p,h_p) \cong ({\cal M}_p,s_p) \cong {\bf R}^{1,1};\eqno(6.1.21)$$
in particular, $q^aq_a = V^aV_a$. 

To further describe orthogonal connecting vector fields, impose the A-P condition on $\cal D$. Then, as $\lambda = \tilde\lambda^{-1}$ now in (6.1.3), $\ell^a$ and $D$ are fixed, while $\zeta \mapsto \lambda^2\zeta$. Since $\lambda$ is also $D$ constant, $D\zeta$ transforms just as $\zeta$ does under the remaining freedom allowed under (6.1.3--4). Similarly $D\tilde\zeta$ transforms as $\tilde\zeta$ does. (Or, as $\lambda = \tilde\lambda^{-1}$, $\zeta \mapsto \lambda^2\zeta$, $\tilde\zeta \mapsto \tilde\lambda^2\tilde\zeta$, from (6.1.19) \& (6.1.6) $D\zeta \mapsto \lambda^2D\zeta$ and $D\tilde\zeta \mapsto \tilde\lambda^2D\tilde\zeta$.) Thus, $(D\zeta)\tilde m^a + (D\tilde\zeta)m^a$ is invariant under the remaining freedom in (6.1.3--4) allowed under the A-P condition, modulo $\ell^a$, just as $\zeta\tilde m^a + \tilde\zeta m^a$ is.
Recalling from (6.1.15) that $DV^a \in {\cal D}^\perp$, the projection of $DV^a$ into $\cal Q$ may be represented by $(D\zeta)m^a + (D\tilde\zeta)\tilde m^a = Dq^a$, with
$\zeta$ and $\tilde\zeta$ satisfying
$$D\pmatrix{\zeta\cr \noalign{\smallskip}\tilde\zeta\cr} = \pmatrix{\rho&\sigma\cr \noalign{\smallskip}\tilde\sigma&\tilde\rho\cr}\pmatrix{\zeta\cr \noalign{\smallskip}\tilde\zeta\cr} =: P_{\cal D}(v)\pmatrix{\zeta\cr \noalign{\smallskip}\tilde\zeta\cr}.\eqno(6.1.22)$$
The tangent planes ${\cal M}_v := \langle \tilde m^a,m^a \rangle_{\bf R}$, parametrized along the integral curve $C_0$ by $v$, can be identified with each other under parallel transport along $C_0$ since $\tilde m^a$ and $m^a$ are $D$ parallel under the A-P condition. Denoting this identification by $\cal M$, one can then view (6.1.22) as a nonautonomous linear system of first order ODE's in ${\cal M} \cong {\bf R}^{1,1}$.

The matrix $P_{\cal D}(v)$ defining the linear system (6.1.22) is generic. It's eigenvalues are
$${\rho+\tilde\rho \pm\sqrt{(\rho-\tilde\rho)^2+4\sigma\tilde\sigma} \over 2},\eqno(6.1.23)$$
which will be a complex conjugate pair when $\sigma\tilde\sigma < 0$ and $\vert 4\sigma\tilde\sigma \vert > (\rho-\tilde\rho)^2$, but otherwise both are real. In general, one would not expect (asymptotic) stability ([1], \S 23) about the zero solution. (6.1.23) suggests that $\rho$, $\tilde\rho$ and $\sigma\tilde\sigma$ contain geometric information about auto-parallel $\cal D$. Indeed, by linear algebra
$$\nabla_a\ell^b = n_aD\ell^b + \ell_aD'\ell^b - \tilde m_a\delta\ell^b - m_a\triangle\ell^b.\eqno(6.1.24)$$
Under the present assumptions the null tetrad is parallel along integral curves; in particular, $D\ell^b = 0$. Moreover, since $\ell^b$ is null, $0 = D'(\ell_b\ell^b) = \ell_bD'\ell^b$. It follows that, noting (2.10),
$$\nabla_a\ell^a = -(\tilde m_a\delta\ell^a + m_a\triangle\ell^a) = \rho+\tilde\rho.\eqno(6.1.25)$$
Moreover, $\nabla_{[a}\ell_{b]} = \ell_{[a}D'\ell_{b]} - \tilde m_{[a}\delta\ell_{b]} - m_{[a}\triangle\ell_{b]}$. Substituting from (2.10) in the right-hand side or simply using (4.1), one has in the present circumstances
$$\nabla_{[a}\ell_{b]} = (\tilde\tau + \alpha + \tilde\beta)\ell_{[a}m_{b]} + (\tau + \tilde\alpha + \beta)\ell_{[a}\tilde m_{b]} - (\rho-\tilde\rho)m_{[a}\tilde m_{b]},$$
whence
$$\nabla_{[a}\ell_{b]}\nabla^{[a}\ell^{b]} = -{(\rho-\tilde\rho)^2 \over 2} \leq 0.\eqno(6.1.26)$$
From (6.1.25--26), one obtains geometric formulae for $\rho$ and $\tilde\rho$ in terms of the local section $\ell^a$ of $\cal D$ (where $\ell^a = o^Ao^{A'} = \partial_v$ and $Do^A = Do^{A'} = 0$). On the other hand, using (2.10),
$$\eqalign{\nabla_{(a}\ell_{b)} &= \ell_{(a}D'\ell_{b)} - \tilde m_{(a}\delta\ell_{b)} - m_{(a}\triangle\ell_{b)}\cr
&= (\gamma+\tilde\gamma)\ell_a\ell_b + (\tilde\tau - \alpha - \tilde\beta)\ell_{(a}m_{b)} + (\tau - \tilde\alpha - \beta)\ell_{(a}\tilde m_{b)} - \tilde\sigma m_am_b - \sigma\tilde m_a\tilde m_b - (\rho+\tilde\rho)m_{(a}\tilde m_{b)},\cr}$$
whence
$$\nabla_{(a}\ell_{b)}\nabla^{(a}\ell^{b)} = 2\sigma\tilde\sigma + {(\rho+\tilde\rho)^2 \over 2}.\eqno(6.1.27)$$
From (6.1.25 \& 27), one therefore obtains a geometric formula for $\sigma\tilde\sigma$. (Note that when $\cal D$ is not auto-parallel, by (6.1.9) at least one of $o_ADo^A = 0$ and $o_{A'}Do^{A'} = 0$ fails. If $o_ADo^A \not= 0$, then $Do^A$ and $o_Bo_{A'}\nabla^ao^B$ are not proportional to $o^A$ and $o_Ao^B\nabla_bo^A$ is not proportional to $o_{B'}$. If one then chooses $\iota^A$ proportional to one of the first two, or $\iota^{B'}$ to the third, then $\epsilon = 0$, $\rho = 0$, and $\sigma = 0$ respectively, i.e., the vanishing of these quantities has no geometrical significance in the non-auto-parallel case. Similar remarks apply to the tilde quantities)

Note from (6.1.6) that, with $\kappa = \tilde\kappa = 0$, and $\lambda = \tilde\lambda^{-1}$ (as it does under the A-P condition), while $\rho$ and $\tilde\rho$ are invariant, only the product $\sigma\tilde\sigma$ is invariant rather than the individual factors. Thus, while the equations in (6.1.22) are invariant under the remaining freedom allowed in (6.1.3--4) under the A-P condition, the individual terms are not.

If $v \mapsto B(v)$ is a continuous mapping from ${\bf R}$ to ${\bf R}(2)$ such that $B(w)B(v) = B(v)B(w)$, for all $v$, $w \in {\bf R}$, then the nonautonomous linear system 
$$D\pmatrix{x\cr y\cr} = B(v)\pmatrix{x\cr y\cr},$$
has solution with initial condition $(x_0,y_0)$ (see [30], pp. 234--235):
$$\pmatrix{x(v)\cr y(v)\cr} = \exp\left(\int^v_0 B(w)dw\right)\pmatrix{x_0\cr y_0\cr}.$$
In particular, this result is valid whenever $B(v) = f(v)B$, for some continuous function $f$ and constant matrix $B$. Observe that, with
$$\displaylines{D(v) := {\rho+\tilde\rho \over 2}\pmatrix{1&0\cr 0&1\cr} =: {\rho+\tilde\rho \over 2}E_1 \hskip 1in S(v) := {\rho - \tilde\rho \over 2}\pmatrix{1&0\cr 0&-1\cr} =: {\rho - \tilde\rho \over 2}E_3\cr
\noalign{\vskip -6pt}
\hfill\llap(6.1.28)\cr
\noalign{\vskip -6pt}
I(v) := {\tilde\sigma - \sigma \over 2}\pmatrix{0&-1\cr1&0\cr} =: {\tilde\sigma - \sigma \over 2}E_2 \hskip 1in H(v) := {\tilde\sigma+\sigma \over 2}\pmatrix{0&1\cr 1&0\cr} = {\tilde\sigma+\sigma \over 2}E_4\cr}$$
\vskip 6pt
$$P_{\cal D}(v) := \pmatrix{\rho&\sigma\cr \tilde\sigma&\tilde\rho\cr} = D(v) + S(v) + I(v) + H(v).\eqno(6.1.29)$$
When $P_{\cal D}(v)$ reduces to $D(v)$, $I(v)$, $H(v)$, or $S(v)$, (6.1.22) has solution, 
$$\displaylines{\pmatrix{\zeta(v)\cr \noalign{\smallskip}\tilde\zeta(v)\cr} = \exp\left(\int^v_0{\rho+\tilde\rho \over 2}\right)\pmatrix{\zeta_0\cr \noalign{\smallskip}\tilde\zeta_0\cr}\cr
\noalign{\vskip 9pt}
\pmatrix{\zeta(v)\cr \noalign{\smallskip}\tilde\zeta(v)\cr} = R\left(\int^v_0{\tilde\sigma - \sigma \over 2}\right)\pmatrix{\zeta_0\cr \noalign{\smallskip}\tilde\zeta_0\cr}\cr
\hfill\llap(6.1.30)\cr
\pmatrix{\zeta(v)\cr \noalign{\smallskip}\tilde\zeta(v)\cr} = L\left(\int^v_0{\tilde\sigma+\sigma \over 2}\right)\pmatrix{\zeta_0\cr \noalign{\smallskip}\tilde\zeta_0\cr}\cr
\noalign{\vskip 9pt}
\pmatrix{\zeta(v)\cr \noalign{\smallskip}\tilde\zeta(v)\cr} = \pmatrix{\exp\left(\int^v_0{\rho-\tilde\rho \over 2}\right)\zeta_0\cr \noalign{\smallskip}\exp\left(\int^v_0{\tilde\rho-\rho \over 2}\right)\tilde\zeta_0\cr}\cr}$$
respectively, where
$$R(t) := \pmatrix{\cos t&-\sin t\cr \sin t&\cos t\cr} \hskip 1.25in L(t) := \pmatrix{\cosh t&\sinh t\cr \sinh t&\cosh t\cr}.$$
These four transformations have simple interpretations. The first effects a dilation of the $(\zeta,\tilde\zeta)$-plane; an expansion if $\int^v_0(\rho + \tilde\rho) > 0$ and a contraction if $\int^v_0(\rho + \tilde\rho) < 0$. The second transformation is a rotation of the $(\zeta,\tilde\zeta)$-plane, the third a hyperbolic rotation of the $(\zeta,\tilde\zeta)$-plane, and the fourth scales the $(\zeta,\tilde\zeta)$ coordinates inversely, which is equivalent, via a simple coordinate transformation, to a hyperbolic rotation.

Since $E_1$ commutes with each of $E_2$, $E_3$, and $E_4$, one also obtains the solutions:
$$\pmatrix{\zeta(v)\cr \noalign{\smallskip}\tilde\zeta(v)\cr} = \exp\left(\int^v_0\,\rho\right)R\left(\int^v_0\,\tilde\sigma\right)\pmatrix{\zeta_0\cr \noalign{\smallskip}\tilde\zeta_0\cr},$$
when $\rho=\tilde\rho$ and $\sigma = -\tilde\sigma$;
$$\pmatrix{\zeta(v)\cr \noalign{\smallskip}\tilde\zeta(v)\cr} = \pmatrix{\exp\left(\int^v_0\,\rho\right)\zeta_0\cr \noalign{\smallskip}\exp\left(\int^v_0\,\tilde\rho\right)\tilde\zeta_0\cr},\eqno(6.1.31)$$
when $\sigma = \tilde\sigma = 0$; and
$$\pmatrix{\zeta(v)\cr \noalign{\smallskip}\tilde\zeta(v)\cr} = \exp\left(\int^v_0\,\rho\right)L\left(\int^v_0\,\sigma\right)\pmatrix{\zeta_0\cr \noalign{\smallskip}\tilde\zeta_0\cr},$$
when $\rho=\tilde\rho$ and $\sigma=\tilde\sigma$.
These solutions provide some insight into how $\rho$, $\tilde\rho$, $\sigma$, and $\tilde\sigma$ effect the behaviour of orthogonal connecting vector fields along an integral curve of auto-parallel $\cal D$.

Noting that the induced scalar product of $q^a$ with itself is $-2\zeta\tilde\zeta$, a simple computation yields
$$D(\zeta\tilde\zeta) = (\rho+\tilde\rho)\zeta\tilde\zeta + \sigma\tilde\zeta^2 + \tilde\sigma\zeta^2,\eqno(6.1.32)$$
indicating, as one would expect from the preceding considerations, that the flow does not preserve the induced scalar-product structure of $\cal M$, which in the $(\zeta,\tilde\zeta)$-coordinates is that of the hyperbolic plane ${\bf R}^2_{\rm hb}$ (see [29] for notation) unless $\sigma = \tilde\sigma = 0$ and $\rho = -\tilde\rho$, i.e., $P_{\cal D}(v) = S(v) = \rho E_2$, in which case the flow amounts to a hyperbolic rotation with respect to the coordinates $t := \zeta + \tilde\zeta$ and $z := \zeta - \tilde\zeta$, with respect to which ${\cal M} \cong {\bf R}^{1,1}$, i.e., the flow is indeed induced by an isometry of ${\cal M} \cong {\bf R}^{1,1}$ in this case.

By Liouville's Theorem for linear systems of ODEs ([1], \S 27.6), the flow of the system (6.1.22) will preserve two-forms in $\cal M$, whence the induced notion of volume, precisely when $\tr\bigl(P_{\cal D}(v)\bigr) = 0$, i.e., when $\rho+\tilde\rho = 0$, which is readily seen to be consistent with the special cases of (6.1.30), in which the solutions themselves are linear transformations of the initial condition, whence preservation of two-forms requires the linear transformation to be an element of {\bf SL(2;R)}.

If $\cal D$ is parallel, then by (6.1.12a), for any connecting vector field $V^a$, $\zeta$, $\tilde\zeta$ and $\nu$ in (6.1.16) are all constant along $C_0$; in particular, for orthogonal $V^a$, $\pi(V^a)$ is constant along $C_0$. Within the domain of the Frobenius coordinates, with the A-P condition imposed, the collection of tangent planes $M_v$ along $C_0$ forms a trivial bundle over $C_0$ and the orthogonal connecting vector fields are represented by the constant sections of this bundle.
\vskip 24pt
\noindent {\bf 6.1.33 Example}\hfil\break
Most of the Walker metrics exploited in the recent literature to illustrate and study various geometric themes (see citations in [17]) consist of a single set of Walker coordinates $(u,v,x,y)$. For these spaces, a natural null distribution of type I is provided by ${\cal D} := \langle \partial_u \rangle_{\bf R}$. From [17] (A1.8), $\nabla_{\partial_u}\partial_u = 0$, i.e., $\cal D$ is auto-parallel but not, in general, parallel (unless $c_1 = b_1 = 0$, i.e., $b$ and $c$ are independent of $u$). In particular, $u$ is an affine parameter for $\cal D$. Choosing the spin frames $\{o^A,\iota^A\}$ and $\{o^{A'},\iota^{A'}\}$ to coincide with the Walker spin frames $\{\alpha^A,\beta^A\}$ and $\{\pi^{A'},\xi^{A'}\}$, respectively, is an allowed choice for $\cal D$ and they satisfy the A-P-condition; in particular, these spin frames are $D$-parallel, see \S 5.

The spin coefficients (5.6) for this choice of spin frames confirms that $\cal D$ is auto-parallel (6.1.9), that the integral curves are affinely parametrized (6.1.11), and that of course the spin frames are $D$-parallel. The components of connecting vector fields $V^a$ satisfy (6.1.19)
$$D\pmatrix{\eta\cr \zeta\cr \noalign{\smallskip}\tilde\zeta\cr \nu\cr} = \pmatrix{0&0&-{c_u \over 2}&{a_u \over 2}\cr 0&0&-{b_u \over 2}&{c_u \over 2}\cr 0&0&0&0\cr \noalign{\smallskip}0&0&0&0\cr}\pmatrix{\eta\cr \zeta\cr \noalign{\smallskip}\tilde\zeta\cr \nu\cr},\eqno(6.1.33{\rm a})$$
which can be directly integrated along an integral curve $C_0$ with respect to the Walker coordinate $u$ to yield
$$\nu = \nu_0 \hskip .5in \tilde\zeta = \tilde\zeta_0 \hskip .5in \zeta(u) = {b_0 - b(u) \over 2}\tilde\zeta_0 + {c(u) - c_0 \over 2}\nu_0 \hskip .5in \eta(u) = {c_0 - c(u) \over 2}\tilde\zeta_0 + {a(u) - a_0 \over 2}\nu_0,\eqno(6.1.33{\rm b})$$
where the subscript 0 denotes values at $u=0$.

For orthogonal connecting vector fields, set $\nu_0 = 0$ in (6.1.33b). The projection into ${\cal Q}_u$ is represented in $M_u$ by
$$\pmatrix{\zeta(u)\cr \noalign{\smallskip}\tilde\zeta(u)\cr} = \pmatrix{1&{b_0 - b(u) \over 2}\cr 0&1\cr}\pmatrix{\zeta_0 \cr \noalign{\smallskip}\tilde\zeta_0\cr},\eqno(6.1.33{\rm c})$$
which amounts, in $\cal M$, to a linear shear of $(\zeta_0,\tilde\zeta_0)$ depending on $u$.

The distribution ${\cal D}' := \langle \partial_v \rangle_{\bf R}$ is also null of type I, also auto-parallel, and $v$ is an affine parameter. Spin frames which satisfy the A-P condition for ${\cal D}'$ are
$$o^A := \beta^A \hskip .75in \iota^A := -\alpha^A \hskip .75in o^{A'} := \pi^{A'} \hskip .75in \iota^{A'} := \xi^{A'},\eqno(6.1.33{\rm d})$$
where $\alpha^A$, $\beta^A$, $\pi^{A'}$, and $\xi^{A'}$ are as above. Since these spin frames for ${\cal D}'$ are obtained by rearranging the elements of the spin frames for $\cal D$, it is a simple matter to derive the spin coefficients for these spin frames. If quantities with respect to these spin frames are denoted by hatted quantities, one computes that $\hat D = \triangle$, $\hat D' = -\delta$, $\hat\delta = D'$, and $\hat\triangle = -D$. From (2.7--8), one readily computes the spin coefficients with respect to the spin frames (6.1.33c), in terms of the given Walker coordinates $(u,v,x,y)$, to be:
\vskip 6pt
$$\vcenter{\openup1\jot \halign{$\hfil#$&&${}#\hfil$&\qquad$\hfil#$\cr
\epsilon &= 0 & \epsilon' &= \displaystyle{c_1 - b_2 \over 4} & \tilde\epsilon &= 0 &\tilde\epsilon' &= \displaystyle-\left({c_1+b_2 \over 4}\right)\cr
\alpha &= 0 & \alpha' &= \displaystyle{c_2 - a_1 \over 4} & \tilde\alpha &= \displaystyle{c_2 + a_1 \over 4} & \tilde\alpha' &= 0\cr
\beta &= \displaystyle{c_2 - a_1 \over 4} & \beta' &= 0 & \tilde\beta &= 0 & \tilde\beta' &= \displaystyle{c_2 + a_1 \over 4}\cr
\gamma &= \displaystyle{b_2 - c_1 \over 4} & \gamma' &= 0 & \tilde\gamma &= \displaystyle{c_1 + b_2 \over 4} & \tilde\gamma' &= 0\cr
\kappa &= 0 & \kappa' &= \displaystyle{b_1 \over 2} & \tilde\kappa &= 0 & \tilde\kappa' &= \displaystyle{2c_4 - 2b_3 - cc_1 - bc_2 + ab_1 + cb_2 \over 4}\cr
\rho &= 0 & \rho' &= \displaystyle-{c_1 \over 2} & \tilde\rho &= 0 & \tilde\rho' &= 0\cr
\sigma &= \displaystyle-{a_2 \over 2} & \sigma' &= 0 & \tilde\sigma &= 0 & \tilde\sigma' &= \displaystyle{2c_3 - 2a_4 + ba_2 + ca_1 - ac_1 -cc_2 \over 4}\cr
\tau &= \displaystyle-{c_2 \over 2} & \tau' &= 0 & \tilde\tau &= 0 & \tilde\tau' &= 0\cr}}$$
One notes that precisely the same spin coefficients vanish and one obtains in place of (6.1.33a)
$$D\pmatrix{\eta\cr \zeta\cr \noalign{\smallskip}\tilde\zeta\cr \nu\cr} = \pmatrix{0&0&{c_v \over 2}&{b_v \over 2}\cr 0&0&-{a_v \over 2}&-{c_v \over 2}\cr 0&0&0&0\cr \noalign{\smallskip}0&0&0&0\cr}\pmatrix{\eta\cr \zeta\cr \noalign{\smallskip}\tilde\zeta\cr \nu\cr}$$
Thus, the connecting vector fields for ${\cal D}'$ are little different to those for $\cal D$ and the orthogonal connecting vector fields are completely analogous, but with $a$ as a function of $v$ taking the role of $b$ as a function of $u$. This result is, of course, consistent with the fundamental symmetry of Walker coordinates, [17], (A1.7)
\vskip 24pt
Returning to an arbitrary null distribution $\cal D$ of type I, the influence of curvature can be made explicit by considering the Jacobi equation on $C_0$. For any connecting vector field $V^a$ along $C_0$, by (6.1.14)
$$\eqalign{D^2V^a &= D(V^b\nabla_b\ell^a)\cr
&= \ell^c\nabla_cV^b\nabla_b\ell^a\cr
&= V^b\nabla_b\ell^c\nabla_c\ell^a + [\ell,V]^b\nabla_b\ell^a + R_{bcd}{}^a\ell^bV^c\ell^d\cr
&= V^b\nabla_bD\ell^a + R_{bcd}{}^a\ell^bV^c\ell^d.\cr}$$
Assuming $\cal D$ is auto-parallel and $(v,x,y,z)$ are Frobenius coordinates with $v$ an affine parameter for the integral curves of $\cal D$, then the first summand on the right in the last equality vanishes and one obtains the Jacobi equation:
$$D^2V^a = R_{bcd}{}^a\ell^bV^c\ell^d\eqno(6.1.34)$$
Assuming now that the A-P condition has been imposed, substituting (6.1.16) into (6.1.34) yields
$$(D^2\eta)\ell^a + (D^2\zeta)\tilde m^a + (D^2\tilde\zeta)m^a + (D^2\nu)n^a = R_{bcd}{}^a\ell^b(\eta\ell^c + \zeta\tilde m^c + \tilde\zeta m^c + \nu n^c)\ell^d.$$
Taking components with respect to the null tetrad yields:
$$\displaylines{D^2\nu = 0 \hskip .5in D^2\eta = \zeta(\tilde\Psi_1 + \Phi_{10}) + \tilde\zeta(\Psi_1 + \Phi_{01}) + \nu(\Psi_2 + \tilde\Psi_2 + 2\Phi_{11} - 2\Lambda)\cr
\noalign{\vskip -6pt}
\hfill\llap(6.1.35)\cr
\noalign{\vskip -6pt}
D^2\zeta = -\zeta\Phi_{00} - \tilde\zeta\Psi_0 - \nu(\Psi_1 + \Phi_{01}) \hskip .75in
D^2\tilde\zeta = -\zeta\tilde\Psi_0 - \tilde\zeta\Phi_{00} - \nu(\tilde\Psi_1 + \Phi_{10}).\cr}$$

Putting 
$$N_{\cal D} := \pmatrix{0&-(\tilde\Psi_1+\Phi_{10})&-(\Psi_1 + \Phi_{01})&2\Lambda - 2\Phi_{11} - \Psi_2 - \tilde\Psi_2\cr 0&\Phi_{00}&\Psi_0&\Psi_1 + \Phi_{01}\cr 0&\tilde\Psi_0&\Phi_{00}&\tilde\Psi_1+\Phi_{10}\cr 0&0&0&0\cr},\eqno(6.1.36)$$
and denoting by $Z$ the column vector in (6.1.19), (6.1.35) and (6.1.19) in the auto-parallel case together yield $-N_{\cal D}Z = D^2Z = D(M_{\cal D}Z) = (DM_{\cal D})Z + M_{\cal D}DZ = (DM_{\cal D})Z + M_{\cal D}^2Z$, whence
$$DM_{\cal D} = -(M_{\cal D}^2+N_{\cal D}).\eqno(6.1.37)$$
(6.1.37) consists of four pairs of equations, with each equation in a pair the tilde of the other, plus one equation which is invariant under tilde:
$$\vcenter{\openup1\jot \halign{$\hfil#$&&${}#\hfil$&\qquad$\hfil#$\cr
D(\alpha+\tilde\beta) &= -\rho(\alpha+\tilde\beta) - \tilde\sigma(\tilde\alpha+\beta) + \tilde\Psi_1 + \Phi_{10} & D(\tilde\alpha+\beta) &= -\tilde\rho(\tilde\alpha+\beta) - \sigma(\alpha+\tilde\beta) + \Psi_1 + \Phi_{01}\cr
D\rho &= -(\rho^2 + \sigma\tilde\sigma) - \Phi_{00} & D\tilde\rho &= -(\tilde\rho^2 + \sigma\tilde\sigma) - \Phi_{00}\hfill\cr
D\sigma &= -(\rho+\tilde\rho)\sigma - \Psi_0 & D\tilde\sigma &= -(\rho+\tilde\rho)\tilde\sigma - \tilde\Psi_0\cr
D\tau &= -(\tau\rho+\tilde\tau\sigma) - \Psi_1 - \Phi_{01} & D\tilde\tau &= -(\tilde\tau\tilde\rho+\tau\tilde\sigma) - \tilde\Psi_1 - \Phi_{10}\cr}}\eqno(6.1.38)$$
$$D(\gamma + \tilde\gamma) = \tau(\alpha+\tilde\beta) + \tilde\tau(\tilde\alpha+\beta) + \Psi_2 + \tilde\Psi_2 + 2\Phi_{11} - 2\Lambda.$$
Recalling that in the current circumstances, $\kappa = \tilde\kappa = \epsilon = \tilde\epsilon = \tau' = \tilde\tau' = 0$, and (2.11) holds, then one readily checks that the equations of (6.1.38) are instances of (3.4); namely: $({\rm h}) + \widetilde{({\rm i})}$ and $\widetilde{({\rm h})} + ({\rm i})$; (a) and $\widetilde{({\rm a})}$; (b) and $\widetilde{({\rm b})}$; (c) and $\widetilde{({\rm c})}$; $({\rm k}) + \widetilde{({\rm k})}$; respectively.

Observe from the third pair that
$$\sigma \equiv 0\ \Rightarrow\ \Psi_0 = 0 \hskip 1.25in \tilde\sigma \equiv 0\ \Rightarrow\ \tilde\Psi_0 = 0,\eqno(6.1.39)$$
i.e., the vanishing of $\sigma$($\tilde\sigma$) implies $o^A$($o^{A'}$) is a Weyl principal spinor (WPS) of $\Psi_{ABCD}$($\Psi_{A'B'C'D'}$). In (6.1.33), $\tilde\sigma = 0$, $\tilde\Psi_0 = 0$ and $\pi^{A'}$ is indeed a WPS (in fact, of multiplicity at least two, by [17] 2.5) while $\sigma = -b_1/2$ and $\Psi_0 = b_{11}/2$, i.e., in this case $\Psi_0 = -\sigma_1$, and $\alpha^A$ is not, in general, a WPS. Unlike the situation for Lorentzian signature, the vanishing of $\Phi_{00}$, $\rho$ and $\tilde\rho$ along $C_0$ only entails the vanishing of $\sigma\tilde\sigma$ and not of $\sigma$ and $\tilde\sigma$ separately, as exemplified by Walker geometry.

When $\cal D$ is additionally parallel, then (6.1.38) reduce to
$$\vcenter{\openup1\jot \halign{$\hfil#$&&${}#\hfil$&\qquad$\hfil#$\cr
\Psi_0 &= \tilde\Psi_0 = \Phi_{00} = 0 & \Psi_1 + \Phi_{01} &= 0 = \tilde\Psi_0 + \Phi_{10}\cr
D(\alpha+\tilde\beta) &= 0 = D(\tilde\alpha+\beta) & D(\gamma+\tilde\gamma) &= \Psi_2 + \tilde\Psi_2 - 2\Lambda.\cr}}$$
As already noted, these equations are instances of (3.4). Moreover, when $\cal D$ is parallel, (3.4)(f)\&($\tilde{\rm f}$) give, respectively,
$$\Psi_2 + 2 \Lambda = 0 = \tilde\Psi_2 + 2\Lambda.$$

Observe from the second pair of (6.1.38) that
$$D(\rho-\tilde\rho) = \tilde\rho^2 - \rho^2 = (\tilde\rho-\rho)(\tilde\rho+\rho).\eqno(6.1.40)$$
Hence, the propagation of $\rho-\tilde\rho$ is unaffected by curvature. Moreover, $\rho-\tilde\rho \equiv 0$ is a solution along $C_0$, whence the unique solution which vanishes at any point, i.e., if $\rho=\tilde\rho$ at any point on $C_0$, this equality is valid everywhere along $C_0$.

Defining, for two Jacobi fields $V^a$ and $W^a$ along $C_0$,
$$\Sigma(V,W) := {V_aDW^a - W_aDV^a  \over 2}\eqno(6.1.41)$$
one readily computes
$$D\bigl(\Sigma(V,W)\bigr) = {1 \over 2}(V_aD^2W^a - W_aD^2V^a) = {1 \over 2}(R_{abcd}\ell^aW^b\ell^cV^d - R_{abcd}\ell^aV^b\ell^cW^d) = 0,$$
whence $\Sigma$ is a well defined symplectic form on the space of Jacobi fields along $C_0$. Putting $W^a =: \omega\ell^a + \tilde\xi m^a + \xi\tilde m^A + \chi n^a$, then from (6.1.16 \& 19), one computes
$$\Sigma(V,W) = {(\rho-\tilde\rho)(\zeta\tilde\xi - \tilde\zeta\xi) + (\tilde\tau + \alpha + \tilde\beta)(\nu\xi - \zeta\chi) + (\tau+\tilde\alpha+\beta)(\nu\tilde\xi - \tilde\zeta\chi) \over 2}.\eqno(6.1.42)$$

Suppose $N_{\cal D}$ vanishes along $C_0$ so that $DZ = M_{\cal D}Z$ and $DM_{\cal D} = -M_{\cal D}^2$. Since the null tetrad is parallel along $C_0$, each tangent space $T_vM$ along $C_0$ can be identified with ${\bf R}^4$ by identifying the null tetrad in $T_vM$ with the standard basis of ${\bf R}^4$. Then $DZ = M_{\cal D}Z$ may be viewed as a nonautonomous linear system of first-order ODE's in ${\bf R}^4$. Let $\{Z_1,Z_2,Z_3,Z_4\}$ be the fundamental system of solutions to this system satisfying $Z_i(0) = e_i$, $i=1,\ldots,4$, the elements of the standard basis. Let $\Xi$ denote the matrix whose $j$'th column consists of the components of $Z_j$ with respect to the standard basis, whence $\det(\Xi)$ is the Wronskian of the system of ODE's. So, $D\Xi = M_{\cal D}\Xi$ but $D^2\Xi = {\bf 0}$, whence the entries of $\Xi$ are linear in $v$. Write $\Xi = \Xi_1v + \Xi_2$, for some elements $\Xi_1$, $\Xi_2$ of ${\bf R}(4)$. Now $\Xi_2 = \Xi(0) = 1_4$, by construction. Also, $\Xi_1 = D\Xi = M_{\cal D}\Xi$ entails $\Xi_1 = M_{\cal D}(0)\Xi(0) = M_{\cal D}(0)$ and
$$M_{\cal D} = \Xi_1\Xi^{-1} = M_{\cal D}(0)\bigl(M_{\cal D}(0)v + 1_4\bigr)^{-1},\eqno(6.1.43)$$
which determines the entries of $M_{\cal D}$ along $C_0$ in terms of their initial values.
\vskip 24pt
\noindent {\bf 6.1.44 Example}\hfil\break
Consider the distribution ${\cal D} = \langle \partial_u \rangle_{\bf R}$  of (6.1.33). From (5.6 \& 11), (6.1.19 \& 36), one finds
$$M_{\cal D} = {1 \over 2}\pmatrix{0&0&-c_1&a_1\cr 0&0&-b_1&c_1\cr0&0&0&0\cr 0&0&0&0\cr} \hskip 1.25in N_{\cal D} = {1 \over 2}\pmatrix{0&0&c_{11}&-a_{11}\cr 0&0&b_{11}&-c_{11}\cr 0&0&0&0\cr 0&0&0&0\cr}.\eqno(6.1.44{\rm a})$$
As $M_{\cal D}^2 = {\bf 0}$, (6.1.37) reduces to $DM_{\cal D} = -N_{\cal D}$, which is obvious from (6.1.44a). Clearly, $N_{\cal D}$ vanishes along $C_0$ iff  $a$, $b$, and $c$ are linear functions of $u$ along $C_0$. If indeed
$${1 \over 2}\pmatrix{-c&a\cr -b&c\cr} = Au+B,\qquad A,\ B \in {\bf R}(2),$$
then $M_{\cal D} = \left({{\bf 0}_2 \atop {\bf 0}_2}{A \atop {\bf 0}_2}\right)$, which is constant along $C_0$. Hence, in (6.1.43), $M_{\cal D} = M_{\cal D}(0)$ and $M_{\cal D}v+1_4 = \left({1_2 \atop {\bf 0}_2}{Av \atop 1_2}\right)$ which has inverse $\left({1_2 \atop {\bf 0}_2}{-Av \atop 1_2}\right)$, whence $M_{\cal D}(M_{\cal D}v+1_4)^{-1} = M_{\cal D}$.

One also computes that $\Sigma$ vanishes identically for this example.
\vskip 24pt
Now consider orthogonal connecting vector fields $V^a$, so $\nu \equiv 0$. As before, the projection $\pi(V^a)$ can be represented by the $q^a$ of (6.1.20). The behaviour of $q^a$ is given by (6.1.22). As before, one identifies each $M_v$ along $C_0$ with ${\bf R}^2$ by identifying $\{\tilde m^a,m^a\}$ with the standard basis of ${\bf R}^2$. If one puts
$$X := \pmatrix{\zeta\cr \noalign{\smallskip}\tilde\zeta\cr} \hskip 1.25in Q_{\cal D} := \pmatrix{\Phi_{00}&\Psi_0\cr \tilde\Psi_0&\Phi_{00}\cr},\eqno(6.1.45)$$
then from (6.1.22 \& 35) $DX = P_{\cal D}X$ and $D^2X = -Q_{\cal D}X$ and one derives
$$DP_{\cal D} = -(P_{\cal D}^2+Q_{\cal D}),\eqno(6.1.46)$$ 
exactly as for (6.1.37). These equations (6.1.37) are just the middle pair of the four pairs in (6.1.38). To the observations (6.1.39--40) one can add the following observation. If $A(v)$ is any quantity along $C_0$ satisfying $DA = (\tr P_{\cal D})A$, then
$$D\bigl((\rho-\tilde\rho)A\bigr) = (\tilde\rho^2-\rho^2)A + (\rho^2 - \tilde\rho^2)A = 0,\eqno(6.1.47)$$
i.e., $(\rho-\tilde\rho)A$ is constant along $C_0$.

For orthogonal connecting vector fields $V^a$ and $W^a$, (6.1.42) reduces to
$$\Sigma(V,W) = {1 \over 2}(\rho-\tilde\rho)(\zeta\tilde\xi - \tilde\zeta\xi).\eqno(6.1.48)$$
In particular, $\Sigma$ is well defined on sections of the bundle $\cal Q$ over $C_0$. Since $A := (\zeta\tilde\xi - \tilde\zeta\xi)/2$ is a constant multiple of the Wronskian of (6.1.22) it satisfies $\dot A = (\tr P_{\cal D})A$, whence the constancy of $\Sigma(V,W)$ along $C_0$ is, in this case, an instance of (6.1.47).

If $Q_{\cal D} \equiv 0$ on $C_0$, one can repeat the argument leading to (6.1.43) with $\Xi$ now the matrix of fundamental solutions of $DX = P_{\cal D}X$ satisfying $\Xi(0) = 1_2$, obtaining
$$P_{\cal D}(v) = P_{\cal D}(0)\bigl(P_{\cal D}(0)v + 1_2\bigr)^{-1}.\eqno(6.1.49)$$
When $P_{\cal D}$ is nonsingular ($\rho\tilde\rho \not= \sigma\tilde\sigma$), then $DP_{\cal D} = -P_{\cal D}^2$ is equivalent to $P_{\cal D}^{-1}(DP_{\cal D})P_{\cal D}^{-1} = -1_2$; from ${\bf 0}_2 = D(P_{\cal D}^{-1}P_{\cal D}) = D(P_{\cal D}^{-1})P_{\cal D} + P_{\cal D}^{-1}DP_{\cal D}$, one therefore deduces $D(P_{\cal D}^{-1}) = -P_{\cal D}^{-1}(DP_{\cal D})P_{\cal D}^{-1} = 1_2$, whence $P_{\cal D}^{-1} = 1_2v+C$, for some $C \in {\bf R}(2)$, viz., $C = [P_{\cal D}(0)]^{-1}$. This result, $P = (1_2v+C)^{-1}$, is equivalent to (6.1.49) when $P_{\cal D}(0)$ is invertible.
As
$$C = \bigl[P(0)\bigr]^{-1} = {1 \over d}\pmatrix{\tilde\rho_0&-\sigma_0\cr -\tilde\sigma_0&\rho_0\cr},\qquad d := \rho_0\tilde\rho_0 - \sigma_0\tilde\sigma_0,$$
then $P(v) = (1_2v + C)^{-1} = P_{\cal D}(0)\bigl(P_{\cal D}(0)v+1_2\bigr)$ amounts to
$$\displaylines{\rho = {\rho_0 + vd \over 1 + v(\rho_0+\tilde\rho_0) + v^2d} \hskip 1.25in \tilde\rho = {\tilde\rho_0 + vd \over 1 + v(\rho_0+\tilde\rho_0) + v^2d}\cr
\noalign{\vskip -6pt}
\hfill\llap(6.1.50)\cr
\noalign{\vskip -6pt}
\sigma = {\sigma_0 \over 1 + v(\rho_0+\tilde\rho_0) + v^2d} \hskip 1.25in \tilde\sigma = {\tilde\sigma_0 \over 1 + v(\rho_0+\tilde\rho_0) + v^2d}\cr}$$
For the example of (6.1.33 \& 44), $P_{\cal D}$ is in fact singular. Now $Q_{\cal D} \equiv 0$ along $C_0$ iff $b$ is a linear function of $u$, and then $DP_{\cal D} = 0$. Substituting $P_{\cal D}(0) = \bigl({0 \atop 0}{\sigma_0 \atop 0}\bigr)$ into (6.1.49) yields
$$P_{\cal D}(u) = \pmatrix{0&\sigma_0\cr 0&0\cr}\pmatrix{1&\sigma_0u\cr 0&1\cr}^{-1} = \pmatrix{0&\sigma_0\cr 0&0\cr}\pmatrix{1&-\sigma_0u\cr 0&1\cr} = \pmatrix{0&\sigma_0\cr 0&0\cr},$$
which is (6.1.50) upon substituting $\rho_0 = \tilde\rho_0 = \tilde\sigma_0 = d = 0$.

Another way to express the geometric content of (6.1.25--27) is as follows. Recalling (6.1.21), along $C_0$ one has, in general, the scalar product spaces $({\cal Q}_v,h_v)$. Let $[X]$ and $[Y]$ be elements of ${\cal Q}_v = {\cal D}^\perp_v/{\cal D}_v$, so $[X]$ has representatives of the form $X^a = \zeta\tilde m^a + \tilde\zeta m^a + \eta\ell^a$, with $\eta$ arbitrary, say, and $[Y]$ of the form $Y^a = \xi\tilde m^a + \tilde\xi m^a + \omega\ell^a$, with $\omega$ arbitrary, say. If $e_{abcd}$ denotes the volume form of $(M,g)$, consider $e_{abcd}X^aY^bn^c\ell^d$. Under (6.1.5), $e_{abcd}X^aY^bn^c\ell^d \mapsto e_{abcd}X^aY^b(\lambda\tilde\lambda)^{-1}n^c(\lambda\tilde\lambda)\ell^d = e_{abcd}X^aY^bn^c\ell^d$, the other terms in the transformation of $n^c$ in (6.1.5) contributing no nonzero terms as $e_{abcd}\vert_{{\cal D}^\perp} = 0$ since ${\cal D}^\perp$ is three dimensional. Hence
$$\Omega\bigl([X],[Y]\bigr) := e_{abcd}X^aY^bn^c\ell^d,\eqno(6.1.51)$$
is a well defined skew form on each ${\cal Q}_v$. Explicitly
$$\eqalignno{\Omega\bigl([X],[Y]\bigr) &= e_{abcd}(\zeta\tilde m^a + \tilde\zeta m^a)(\xi\tilde m^b + \tilde\xi m^b)n^c\ell^d\cr
&= (\zeta\tilde\xi - \tilde\zeta\xi)e_{abcd}\ell^a\tilde m^bn^cm^d\cr
&= \zeta\tilde\xi - \tilde\zeta\xi.&(6.1.52)}$$
Now, in each ${\cal M}_v$, with $T^a :=(\tilde m^a - m^a)/\sqrt2$ and $S^a := (\tilde m^a + m^a)/\sqrt2$, $\{T^a,S^a\}$ is a $\Psi$-ON basis, and if one chooses the orientation class $[T,Z] = [\tilde m,m]$ on ${\cal M}_v$, then one readily checks that the volume form on ${\cal M}_v$ is $m_a \wedge \tilde m_b = m_a \otimes \tilde m_b - \tilde m_a \otimes m_b$, whence the area of the parallelogram spanned by $\zeta \tilde m^a + \tilde\zeta m^a$ and $\xi \tilde m^a + \xi m^a$ is $\zeta\tilde\xi - \tilde\zeta\xi$. By (6.1.21), one may take this expression as the volume for $({\cal Q}_v,h_v)$, i.e., $\Omega$ is the volume form for $({\cal Q}_v,h_v)$.

Now consider the automorphism $F_v$ of $M_v$ defined by $F_v: \tilde m^a \mapsto -\tilde m^a$, $:m^a \mapsto m^a$, i.e., $F_v$ is a product structure on each $M_v$. Observe that
$$s_v\bigl(F_v(\zeta \tilde m^a + \tilde\zeta m^a),\xi\tilde m^a + \tilde\xi m^a\bigr) = (-\zeta\tilde m^a + \tilde\zeta m^a)(\xi\tilde m_a + \tilde\xi m_a) = \zeta\tilde\xi - \tilde\zeta\xi.$$
One readily checks that $F_v$ is anti-orthogonal, equivalently (as $F_v^2 = 1$) skew-adjoint. Clearly, $F_v$ induces an automorphism (also denoted $F_v$) on ${\cal Q}_v$ such that $F_v^2 = 1$, $F_v$ is an anti-orthogonal automorphism of $({\cal Q}_v,h_v)$ (whence skew-adjoint) and $\Omega(\ ,\ ) = h_v(F_v\ ,\ )$, i.e., $(h_v,F_v,\Omega)$ is a (linear) paraHermitian structure on ${\cal Q}_v$.

When $\cal D$ is auto-parallel, with $V^a$ and $W^a$ again orthogonal connecting vector fields along $C_0$, then from (6.1.48) one notes that at any point 
$$\Sigma_v(V,W) = {\rho-\tilde\rho \over 2}\Omega\bigl(\pi_v(V),\pi_v(W)\bigr).\eqno(6.1.53)$$
Note that neither $h$ nor $\Omega$ (or the equivalents in $M_v$) are $D$-parallel along $C_0$, though $F_v$ is by construction.

With $X^a$ and $Y^a$ in ${\cal D}^\perp$ as above, since $\ell^b\nabla_b\ell^a = 0$ and with $\ell^b\nabla_a\ell_b = 0$, the expression $X^aY^b\nabla_a\ell_b$ depends only on $[X]$ and $[Y]$ and defines a tensor on ${\cal Q}_v$:
$$T_v\bigl([X],[Y]\bigr) := X^aY^b\nabla_a\ell_b.\eqno(6.1.54)$$
Suppressing the subscript $v$, $T_{ab} = T_{(ab)} + T_{[ab]}$, where $T_{[ab]} \propto \Omega$, since ${\cal Q}_v$ is two dimensional. Indeed, evaluating (6.1.54) yields
$$\eqalign{T\bigl([X],[Y]\bigr) &= (\zeta\tilde m^a + \tilde\zeta m^a)(\xi\tilde m^b + \tilde\xi m^b)\nabla_a\ell_b\cr
&= \zeta\xi\tilde m^a\tilde m^b\nabla_a\ell_b + \zeta\tilde\xi\tilde m^a m^b\nabla_a\ell_b + \tilde\zeta\xi m^a \tilde m^b\nabla_a\ell_b + \tilde\zeta\tilde\xi m^a m^b\nabla_a\ell_b\cr
&= -(\zeta\xi\tilde\sigma + \zeta\tilde\xi\rho +\tilde\zeta\xi\tilde\rho + \tilde\zeta\tilde\xi\sigma).\cr}$$
It follows that one can represent $T_{ab}$ in $M_v$ by
$$-(\sigma\tilde m_a\tilde m_b + \tilde\sigma m_am_b) + {\rho+\tilde\rho \over 2}s_{ab} + {\rho-\tilde\rho \over 2}\Omega_{ab}.$$
Denoting the first tensor by $S_{ab}$ it is easy to see that it induces a well defined symmetric and trace free tensor (which I shall also denote by $S_{ab}$) on ${\cal Q}_u$. Thus, one can write
$$T_{ab} = S_{ab} + {\rho+\tilde\rho \over 2}h_{ab} + {\rho-\tilde\rho \over 2}\Omega_{ab}.\eqno(6.1.55)$$
This tensorial expression is an alternative formulation of the information contained in (6.1.25--27).

Let ${\cal D}_1$ and ${\cal D}_2$ be two auto-parallel null distributions of type I in $(M,g)$ with a common integral curve $C$. Assume Frobenius coordinates have been chosen for each distribution so that integral curves are affinely parametrized. Since $v_1$ and $v_2$ are both affine parametrizations for $C$, $v_1 = av_2+b$, for constants $a$ and $b$, and one can replace $v_2$ by $av_2+b$ everywhere on its domain in order to employ a single affine parameter $v$ for $C$ valid for both distributions; whence $\ell^a_1 = \ell^a_2$ on $C$. Furthermore, when imposing the A-P condition, one can choose the null tetrads for the two distributions to coincide along $C$. Hence, one obtains a common identification of $T_vM$ along $C$ with ${\bf R}^4$, of $\left({\cal D}_i^\perp\right)_v$ with ${\bf R}^3$, and of $({\cal Q}_i)_v$ with ${\bf R}^2$, $i=1,$ 2 and, indeed, identifications of ${\cal D}^\perp_1$ with ${\cal D}^\perp_2$ and ${\cal Q}_1$ with ${\cal Q}_2$, along $C$. Along $C$, each distribution defines $P_{{\cal D}_i}(v)$ as in (6.1.22), and these may be viewed as endomorphisms of a common space ${\cal M} \cong {\bf R}^{1,1}$ (obtained by identifying each $M_v$ using the $D$-parallel $\tilde m^a$ and $m^a$) which represents $\cal Q$ on $C$ (which results from identifying the ${\cal Q}_v$ with each other, as before).

With these arrangements, it follows that ${\cal D}_1$ and ${\cal D}_2$ have a common orthogonal connecting vector field iff $DX = P_{{\cal D}_1}X = P_{{\cal D}_2}X$, with $X$ now as in (6.1.45). Hence, $\det(P_{{\cal D}_1}-P_{{\cal D}_2}) = 0$ along $C$ is a necessary condition. Now
$$\det(P_{{\cal D}_1}-P_{{\cal D}_2}) = (\rho_1-\rho_2)(\tilde\rho_1-\tilde\rho_2) - (\sigma_1-\sigma_2)(\tilde\sigma_1-\tilde\sigma_2).$$
By direct computation, using (6.1.38),
$$D\bigl(\det(P_{{\cal D}_1}-P_{{\cal D}_2})\bigr) = -(\rho_1+\tilde\rho_1+\rho_2+\tilde\rho_2)\det(P_{{\cal D}_1}-P_{{\cal D}_2}).\eqno(6.1.56)$$
Hence, if $\det(P_{{\cal D}_1}-P_{{\cal D}_2})$ vanishes at any point of $C$, it vanishes all along $C$. Assuming $P_{{\cal D}_1} \not= P_{{\cal D}_2}$, then the nullity of $P_{{\cal D}_1}-P_{{\cal D}_2}$ is constant along $C$, either one or zero. So, when $\det(P_{{\cal D}_1}-P_{{\cal D}_2}) = 0$, $\ker(P_{{\cal D}_1}-P_{{\cal D}_2})$ is a well defined line bundle along $C$ whose sections are the common solutions of (6.1.22) for the two distributions. A common solution $q^a$ may be obtained as the projection of unique orthogonal connecting vector fields $V^a_i$ for ${\cal D}_i$ by choosing $V_i^a = q^a + \eta_i\ell^a$ where $\eta_i\ell^a$ is determined by the appropriate first order ODE in (6.1.19) subject to the same initial condition at $v=0$. Since $\pi_1(V^a_1) = \pi_2(V^a_2)$, and is represented by $q^a \in {\cal M}_v$, $q^a$ connects $C$ to another integral curve common to both distributions (at least locally), though the two distributions will not have a common parametrization for this other integral curve unless $V^a_1=V^a_2$.

Note that the condition $\det(P_{{\cal D}_1}-P_{{\cal D}_2}) = 0$ is precisely that $P_{{\cal D}_1}$ and $P_{{\cal D}_2}$ are null separated in the standard split-Quaternion structure of ${\bf R}(2)$. Thus, at any point $p$ on $C$, ${\cal P} := {\bf R}(2) \cong {\bf R}^{2,2}$ represents the possible one-dimensional null distributions containing $C_0$. For a given such $\cal D$, the different values of $P_{\cal D}$ at different points on $C$ may be viewed as different coordinatizations of $\cal P$. With this interpretation, (6.1.56) indicates that only the conformal geometry of $\cal P$ plays a role. As such, the conformal compactification of ${\cal P} \cong {\bf R}^{2,2}$ may be useful for describing null congruences for which $\rho$, $\tilde\rho$, $\sigma$, or $\tilde\sigma$ become infinite on $C$.
\vskip 24pt
\noindent {\bf 6.2 Null Distributions of Type II}
\vskip 12pt
Let $N$ be a null distribution of type II on $(M,g)$, i.e., each $N_p$, $p \in M$, is totally null. The quadric Grassmannian ${\cal Q}_2({\bf R}^{2,2})$, as is well known, consists of two components ${\cal Q}_2^+({\bf R}^{2,2})$ and ${\cal Q}_2^-({\bf R}^{2,2})$, whose elements are called $\alpha$-planes and $\beta$-planes respectively, see [17] \S 2 for further details.

I will suppose $(M,g)$ is orientable so that $N$ is either a distribution of $\alpha$-planes or a distribution of $\beta$-planes; in other words, the bivectors formed from local frames for the distribution are self-dual (SD) or anti-self-dual (ASD) respectively. An $\alpha$-plane has a spinor description as
$$Z_{[\pi]} := \{\,\eta^A\pi^{A'}:\hbox{ some fixed } [\pi^{A'}] \in {\bf P}S';\ \eta^A \in S\,\}$$
while a $\beta$-plane has a spinor description as
$$W_{[\eta]} := \{\,\eta^A\pi^{A'}:\hbox{ some fixed } [\eta^A] \in {\bf P}S;\ \pi^{A'} \in S'\,\}.$$
Clearly, the $\alpha$-plane $Z_{[\pi]}$ at $p \in M$ is determined by the projective spinor $[\pi^{A'}]$ and, at least locally, a distribution of $\alpha$-planes is equivalent to a section of the bundle ${\bf P}S'$ of projective primed spinors (globally so when $(M,g)$ is ${\bf SO^\bfplus}$-orientable). The (local) field $[\pi^{A'}]$ is called the projective spinor field associated to the distribution $N$ which, therefore, may be denoted $Z_{[\pi]}$. Similar remarks apply to the case of distributions of $\beta$-planes. By a choice of orientation, one may restrict attention to distributions of $\alpha$-planes, which for convenience I shall call $\alpha$-distributions.

For an $\alpha$-distribution $Z_{[\pi]}$, a local spinor field whose projectivization equals $[\pi^{A'}]$ is called a {\sl local scaled representative\/} (LSR) for the distribution. In terms of an LSR, the condition for integrability of the $\alpha$-distribution $Z_{[\pi]}$ is
$$\pi_{A'}\pi^{B'}\nabla_{BB'}\pi^{A'} = 0,\eqno(6.2.1)$$
which is valid for all LSRs if for any. In the context of Lorentzian geometry, (6.2.1) characterizes shear-free null geodesic congruences, see [23] \S 7.3--4. Pleba\'nski and Hacyan [26] recognized that in the context of complex general relativity (6.2.1) relates to totally null two-surfaces rather than null curves; see also [23] (7.3.18). Four-dimensional neutral geometry is more analogous to complex general relativity than to Lorentzian geometry. 

Note that (6.2.1) is also the condition for the distribution to be auto-parallel. I now restrict attention to integrable, therefore auto-parallel, $\alpha$-distributions. By (6.1.7), the integrable surfaces, called $\alpha$-{\sl surfaces}, are totally geodesic and $\nabla_a$ induces a linear connection within $\alpha$-surfaces.

Working locally, pick an LSR $\pi^{A'}$ as the first element of a primed spin frame: $o^{A'} := \pi^{A'}$. The freedom in construction of spin frames satisfying this choice is
$$\vcenter{\openup1\jot \halign{$\hfil#$&&${}#\hfil$&\qquad$\hfil#$\cr
o^{A'} &\mapsto \tilde\lambda o^{A'} & \iota^{A'} &\mapsto \tilde\lambda^{-1}\iota^{A'} + \tilde\mu o^{A'}\cr
o^A &\mapsto \alpha o^A + \beta \iota^A & \iota^A &\mapsto \gamma o^A + \delta \iota^A.\cr}}\eqno(6.2.2)$$
With $A := \left({\alpha \atop \beta}{\gamma \atop \delta}\right) \in {\bf SL(2;R)}$, the freedom is the subgroup
$$\left\{\,\left(A,\pmatrix{\tilde\lambda&\tilde\mu\cr 0&\tilde\lambda^{-1}\cr}\right): A \in {\bf SL(2;R)},\ \tilde\lambda \in {\bf R}^*,\ \tilde\nu \in {\bf R}\,\right\}$$
of ${\bf SL(2;R)} \times {\bf SL(2;R)}$. The resulting freedom in the corresponding null tetrad is
$$\vcenter{\openup1\jot \halign{$\hfil#$&&${}#\hfil$&\qquad$\hfil#$\cr
\ell^a &\mapsto \tilde\lambda(\alpha\ell^a + \beta\tilde m^a) & n^a &\mapsto \tilde\mu(\gamma\ell^a + \delta\tilde m^a) + \tilde\lambda^{-1}(\gamma m^a + \delta n^a)\cr
\tilde m^a &\mapsto \tilde\lambda(\gamma\ell^a + \delta\tilde m^a) & m^a &\mapsto \tilde\mu(\alpha\ell^a + \beta\tilde m^a) + \tilde\lambda^{-1}(\alpha m^a + \beta n^a).\cr}}\eqno(6.2.3)$$
From (2.9), (6.2.1) is equivalent to
$$\tilde\kappa = \tilde\sigma = 0\eqno(6.2.4)$$
for the above choice of spin frames. It then follows from (3.4)(\~b) that $\tilde\Psi_0 = 0$, i.e., $[\pi^{A'}]$ is a Weyl principal spinor (WPS). 

Recalling that WPSs for neutral metrics in four dimensions may be complex valued, see [16], one may ask whether complex-valued solutions of (6.2.1) yield complex WPSs? Recall first (from [16]) that a complex WPS $\pi^{A'} = \gamma^{A'} + i\delta^{A'}$, $\gamma^{A'}$, $\delta^{A'} \in S'$, is required to satisfy $\pi^{A'}\overline\pi_{A'} = -2i\gamma^{A'}\delta_{A'} \not= 0$ in order that it not be trivially complex, i.e., not just a complex scalar multiple of a real WPS.  In order to emphasize the neutral signature in my treatment, I will employ only spin frames from $S'$ and $S$ when taking components of {\bf C}-valued spinors; using complex frames from ${\bf C}S'$ and ${\bf C}S$ in effect means complexifying, whence ignoring the signature. In particular, the spin coefficients are always {\bf R}-valued quantities. Thus, I will refrain from choosing a {\bf C}-valued solution of (6.2.1) as an element of a spin frame.
\vskip 24pt
\noindent {\bf 6.2.5 Interpretation of Complex Solutions of (6.2.1)}\hfil\break
Just as an {\bf R}-valued spinor field defines a rank-two subbundle of the tangent bundle $TM$ (i.e., a two-dimensional distribution on $M$), a {\bf C}-valued spinor field $\pi^{A'}$ defines a rank-two subbundle $H := \{\,\eta^A\pi^{A'}: \eta^A \in {\bf C}S\,\}$ of the complexified tangent bundle ${\bf C}TM$ (a two-(complex-)dimensional smooth complex distribution on $M$). The condition that (local) sections of $H$ be closed under the taking of Lie brackets, which I denote by $[H,H] \leq H$, is still (6.2.1). Since ${\bf C}TM = H \oplus \overline H$, the Newlander-Nirenberg theorem ensures that the almost complex structure $J$, induced on $M$ by ${\bf C}TM = H \oplus \overline H$ in the standard manner, is integrable. In particular, there exist local coordinates $(x^1,\ldots,x^4)$ such that $z^j := x^j + ix^{2+j}$, $j=1$ , 2 are holomorphic coordinates with respect to this complex structure; i.e., $H = \langle \partial_{z^1},\partial_{z^2} \rangle_{\bf C}$ and is identified with the holomorphic tangent space in the standard fashion. In particular, there are complex spinor fields $\zeta_j^A$ such that $\zeta^A_j\pi^{A'} = \partial_{z^j}$, i.e., are holomorphic.

Of course, $\pi^{A'}$ is a solution of (6.2.1) iff $\overline\pi^{A'}$ is, the latter determining the distribution $\overline H$ and the conjugate complex structure; these obvious facts are consistent with the fact that $\pi^{A'}$ is a WPS iff $\overline\pi^{A'}$ is.
\vskip 24pt
With {\bf K} denoting either {\bf R} or {\bf C}, consider {\bf K}-valued solutions of (6.2.1); (6.2.1) is equivalent to each of the following:
$$S_a := \pi_{B'}\nabla_a\pi^{B'} = \omega_A\pi_{A'} \hskip 1.25in T_a := \pi^{B'}\nabla_{AB'}\pi_{A'} = \eta_A\pi_{A'},\eqno(6.2.6)$$
for some spinor fields $\omega_A$ and $\eta_A$. One may then write
$$\nabla_b\pi_{A'} = V_b\pi_{A'} + \omega_B\pi_{B'}\xi_{A'},\eqno(6.2.7)$$
where $\xi^{A'}$ is chosen to satisfy $\xi^{D'}\pi_{D'} = 1$ and $V_b$ is some vector satisfying $\pi^{B'}V_b = \eta_B$ (whence $V_b = \alpha_B\pi_{B'} - \eta_B\xi_{B'}$, for some $\alpha_B$). One then computes
$$(\nabla_b\pi_{A'})(\nabla^b\pi^{A'}) = 2(\eta^D\omega_D).\eqno(6.2.8)$$
\vskip 24pt
\noindent {\bf 6.2.9 Proposition}\hfil\break
Let $\pi^{A'}$ be a smooth {\bf K}-valued spinor field on $M$ satisfying (6.2.1). Then $[\pi^{A'}]$ is a WPS.

Proof. From (6.2.6),
$$\pi^{B'}\pi^{C'}\nabla^B_{C'}S_b = \pi^{B'}\pi^{C'}\nabla^B_{C'}(\omega_B\pi_{B'}) = 0$$
using (6.2.1). Hence,
$$0 = \pi^{B'}\pi^{C'}\nabla^B_{C'}(\pi_{A'}\nabla_b\pi^{A'}) = \pi_{A'}\pi^{B'}\pi^{C'}\nabla^B_{C'}\nabla_b\pi^{A'} + (\pi^{C'}\nabla^B_{C'}\pi_{A'})(\pi^{B'}\nabla_b\pi^{A'}).$$
The second term on the right is $(\eta^B\pi_{A'})(\eta_B\pi^{A'}) = 0$, so
$$\eqalign{0 &= \pi^{A'}\pi^{B'}\pi^{C'}\nabla_{BC'}\nabla^B_{B'}\pi_{A'}\cr
&= \pi^{A'}\pi^{B'}\pi^{C'}[\Square_{C'B'} + {1 \over 2}\epsilon_{C'B'}\Square]\pi_{A'}\qquad\hbox{using (A.5)}\cr
&= \pi^{A'}\pi^{B'}\pi^{C'}\Square_{C'B'}\pi_{A'}\cr
&= \pi^{A'}\pi^{B'}\pi^{C'}\bigl[\tilde\Psi_{C'B'A'E'}\pi^{E'} - \Lambda(\pi_{C'}\epsilon_{B'A'} + \epsilon_{C'A'}\pi_{B'})\bigr]\qquad\hbox{by (A.6)}\cr
&= \tilde\Psi_{A'B'C'D'}\pi^{A'}\pi^{B'}\pi^{C'}\pi^{D'},\cr}$$
i.e., $[\pi^{A'}]$ is a WPS.
\vskip 24pt
\noindent {\bf 6.2.10 Remark}\hfil\break
From (6.2.9), if at a point $p \in M$, for each $[\pi^{A'}] \in {\bf P}S'_p$ there is an integrable $\alpha$-distribution whose plane at $p$ is of the form $Z_{[\pi]}$, then $\tilde\Psi_{A'B'C'D'}$ vanishes at $p$.

In the context of complex general relativity, Penrose ([23] p. 168 and references cited there) observed that for each point $p$ of a complex space-time and for each $\alpha$-plane $Z$ at $p$, there is an $\alpha$-surface $\cal Z$ with $T_p{\cal Z} = Z$ iff $\tilde\Psi_{A'B'C'D'} = 0$. LeBrun \& Mason, [19] \S 3, provide a proof of this result in the context of four-dimensional neutral geometry. In fact, $\tilde\Psi_{A'B'C'D'}$ vanishes iff at each point $p$ and for each $[\pi^{A'}] \in {\bf P}S'_p$ there exists (locally) an integrable $\alpha$-distribution $Z$ such that $Z = Z_{[\pi]}$ at $p$.
\vskip 24pt
By (6.2.9), (6.2.1) is a sufficient condition for $[\pi^{A'}]$, real or complex, to be a WPS. Is it necessary?

In Lorentzian geometry (and complex GR), the relationship between solutions of (6.2.1) and WPS's is addressed by the (generalized) Goldberg-Sachs Theorem (see [23] (7.3.35) and [26]). With minor modifications, one can carry over [23] (7.3.35) to the current circumstances. Several preliminary lemmas will in fact be of subsequent utility, and since as a matter of principle it is important to distinguish the cases of real and complex spinors, I present the argument in full.
\vskip 24pt
\noindent {\bf 6.2.11 Lemma}\hfil\break
For any {\bf K}-valued solution $\pi^{A'}$ of (6.2.1), and noting (6.2.6),
$$\omega_A = \nabla_{AD'}\pi^{D'} - \eta_A.$$

Proof. One computes
$$\eqalign{\omega_A\pi_{A'} &= \pi_{B'}\nabla_a\pi^{B'}\cr
&= -\pi^{B'}\left(\nabla_{A(A'}\pi_{B')} + \nabla_{A[A'}\pi_{B']}\right)\cr
&= -{1 \over 2}\pi^{B'}\left(\nabla_{AA'}\pi_{B'} + \nabla_{AB'}\pi_{A'} + \epsilon_{A'B'}\nabla_{AD'}\pi^{D'}\right)\cr
&= {1 \over 2}\left(\omega_A\pi_{A'} - \eta_A\pi_{A'} + \pi_{A'}\nabla_{AD'}\pi^{D'}\right)\cr}$$
and the result follows immediately.
\vskip 24pt
\noindent {\bf 6.2.12 Lemma}\hfil\break
For any {\bf K}-valued solution $\pi^{A'}$ of (6.2.1),
$$\pi_{C'}\pi^{B'}\nabla_a\nabla_b\pi^{C'} = \omega_A(\eta_B - \omega_B)\pi_{A'}.$$

Proof. From (6.2.1),
$$\eqalign{0 &= \nabla_a(\pi_{C'}\pi^{B'}\nabla_b\pi^{C'})\cr
&= (\nabla_a\pi_{C'})(\pi^{B'}\nabla_b\pi^{C'}) + \pi_{C'}(\nabla_a\pi^{B'})(\nabla_b\pi^{C'}) + \pi_{C'}\pi^{B'}\nabla_a\nabla_b\pi^{C'}\cr
&= (\nabla_a\pi_{C'})(\eta_B\pi^{C'}) + (\nabla_a\pi^{B'})(\omega_B\pi_{B'}) + \pi_{C'}\pi^{B'}\nabla_a\nabla_b\pi^{C'}\cr
&= (-\omega_A\pi_{A'})\eta_B + (\omega_A\pi_{A'})\omega_B + \pi_{C'}\pi^{B'}\nabla_a\nabla_b\pi^{C'}\cr
&= \omega_A(\omega_B - \eta_B)\pi_{A'} + \pi_{C'}\pi^{B'}\nabla_a\nabla_b\pi^{C'}.\cr}$$
\vskip 24pt
\noindent {\bf 6.2.13 Lemma}\hfil\break
For any {\bf K}-valued spinor field $\pi^{A'}$
$$\displaylines{\hfill\pi^{B'}\nabla_b\nabla_a\pi^{C'} = \pi^{B'}\nabla_a\nabla_b\pi^{C'} + \pi_{A'}\Phi_{BA}{}^{C'}{}_{D'}\pi^{D'} + \epsilon_{BA}\tilde\Psi_{A'B'}{}^{C'}{}_{D'}\pi^{B'}\pi^{D'} - \Lambda\epsilon_{BA}\pi_{A'}\pi^{C'};\hfill\llap({\rm a})\cr
\hfill\pi^{B'}\nabla^B_{B'}\nabla_{BD'}\pi^{D'} = \pi^{B'}\nabla_{BD'}\nabla^B_{B'}\pi^{D'}.\hfill\llap({\rm b})\cr}$$
When $\pi^{A'}$ is a solution of (6.2.1):
$$\pi_{C'}\pi^{B'}\nabla_b\nabla_a\pi^{C'} = (\eta_B-\omega_B)\omega_A\pi_{A'}  - \pi_{A'}\Phi_{BAC'D'}\pi^{C'}\pi^{D'} - \epsilon_{BA}\tilde\Psi_{A'B'C'D'}\pi^{B'}\pi^{C'}\pi^{D'}\eqno({\rm c})$$
Proof. Formula (a) follows directly from (A.3--6). Formula (b) follows from (a) upon contracting $B$ on $A$ and $C'$ on $A'$. Formula (c) follows from (a) upon transvection by $\pi_{C'}$ and substitution of (6.2.12) for the first term on the right-hand side.
\vskip 24pt
\noindent {\bf 6.2.14 Lemma}\hfil\break
For any {\bf K}-valued solution $\pi^{A'}$ of (6.2.1):
$$\lambda_{A'} := \pi_{A'}\pi^{B'}\nabla_{BB'}\eta^B = \tilde\Psi_{A'B'C'D'}\pi^{B'}\pi^{C'}\pi^{D'},$$
whence $\lambda_{A'} = 0$ is equivalent to $[\pi^{A'}]$ being a multiple WPS.

Proof. From (6.2.6), $\eta^B(\pi^{B'}\nabla_{BB'}\pi_{A'}) = 0$, so
$$\eqalign{\lambda_{A'} &= \lambda_{A'} + \eta^B\pi^{B'}\nabla_{BB'}\pi_{A'}\cr
&= \pi^{B'}(\pi_{A'}\nabla_{BB'}\eta^B + \eta^B\nabla_{BB'}\pi_{A'})\cr
&= \pi^{B'}\nabla_{BB'}(\eta^B\pi_{A'})\cr
&= \pi^{B'}\nabla_{BB'}(\pi^{C'}\nabla^B_{C'}\pi_{A'})\cr
&= (\pi^{B'}\nabla_{BB'}\pi^{C'})(\nabla^B_{C'}\pi_{A'}) + \pi^{B'}\pi^{C'}\nabla_{BB'}\nabla^B_{C'}\pi_{A'}\cr
&= (\eta_B\pi^{C'})\nabla^B_{C'}\pi_{A'} + \pi^{B'}\pi^{C'}(\Square_{B'C'} + {1 \over 2}\epsilon_{B'C'}\Square)\pi_{A'}\qquad\hbox{using (A.5)}\cr
&= \eta_B\eta^B\pi_{A'} + \pi^{B'}\pi^{C'}\Square_{B'C'}\pi_{A'},\cr}$$
and the desired result follows now from (A.6).
\vskip 24pt
\noindent {\bf 6.2.15 Lemma}\hfil\break
For any {\bf K}-valued solution $\pi^{A'}$ of (6.2.1):
$$\mu_{A'} := \pi_{A'}\pi^{B'}\nabla^B_{B'}\nabla_{BD'}\pi^{D'} - \pi_{A'}\eta^B\nabla_{BD'}\pi^{D'} = \tilde\Psi_{A'B'C'D'}\pi^{B'}\pi^{C'}\pi^{D'} = \lambda_{A'}.$$

Proof. From the proof of (6.2.12), one obtains upon transvecting by $\epsilon^{AB}$
$$\eqalignno{0 &= \epsilon^{AB}\nabla_a(\pi_{C'}\pi^{B'}\nabla_b\pi^{C'})\cr
&= \epsilon^{AB}\bigl[\omega_A(\omega_B - \eta_B)\pi_{A'} + \pi_{C'}\pi^{B'}\nabla_a\nabla_b\pi^{C'}\bigr]\cr
&= (\omega^D\eta_D)\pi_{A'} + \pi_{C'}\pi^{B'}\epsilon^{AB}\nabla_a\nabla_b\pi^{C'}.&(*)\cr}$$
Now,
$$\eqalign{\pi^{B'}\pi^{C'}\nabla_{BA'}\nabla^B_{B'}\pi_{C'} &= \pi^{B'}\epsilon_{A'}{}^{D'}\pi^{C'}\nabla_{BD'}\nabla^B_{B'}\pi_{C'}\cr
&= \pi^{B'}\epsilon^{D'F'}\epsilon_{A'F'}\pi^{C'}\nabla_{BD'}\nabla^B_{B'}\pi_{C'}\cr
&= \pi^{B'}\epsilon^{D'F'}[\epsilon_{F'A'}\pi_{C'}\nabla_{BD'}\nabla^B_{B'}\pi^{C'}]\cr
&= \pi^{B'}\epsilon^{D'F'}[\pi_{F'}\nabla_{BD'}\nabla^B_{B'}\pi_{A'} - \pi_{A'}\nabla_{BD'}\nabla^B_{B'}\pi_{F'}]\cr
&= \pi^{B'}\pi^{D'}\nabla_{BD'}\nabla^B_{B'}\pi_{A'} - \pi^{B'}\pi_{A'}\nabla_{BD'}\nabla^B_{B'}\pi^{D'}\cr
&= \pi^{B'}\pi^{D'}(\Square_{D'B'} + {1 \over 2}\epsilon_{D'B'}\Square)\pi_{A'} - \pi^{B'}\pi_{A'}\nabla_{BD'}\nabla^B_{B'}\pi^{D'}\qquad\hbox{using (A.5)}\cr
&= \tilde\Psi_{A'B'C'D'}\pi^{B'}\pi^{C'}\pi^{D'} - \pi^{B'}\pi_{A'}\nabla_{BD'}\nabla^B_{B'}\pi^{D'}\qquad\hbox{using (A.6)}\cr
&= \tilde\Psi_{A'B'C'D'}\pi^{B'}\pi^{C'}\pi^{D'} - \pi^{B'}\pi_{A'}\nabla^B_{B'}\nabla_{BD'}\pi^{D'},\qquad\hbox{using (6.2.13b)}.\cr}$$
Substituting from $(*)$ into this last equation and rearranging terms yields
$$\pi^{B'}\pi_{A'}\nabla^B_{B'}\nabla_{BD'}\pi^{D'} + (\omega^D\eta_D)\pi_{A'} = \tilde\Psi_{A'B'C'D'}\pi^{B'}\pi^{C'}\pi^{D'}.\eqno(\dag)$$
Transvecting (6.2.11) by $\eta^A$ and substituting for the second term on the left-hand side of $(\dag)$ yields the desired result.
\vskip 24pt
\noindent {\bf 6.2.16 Lemma}\hfil\break
For a {\bf K}-valued solution of (6.2.1), the necessary and sufficient condition for a {\sl smooth} (local) solution $\phi$ of
$$\pi^{B'}\nabla_{AB'}\phi = \gamma_A,\eqno(6.2.16{\rm a})$$
for a given {\sl smooth} spinor field $\gamma_A$ (without loss of generality, take $\gamma^A$ to be {\bf R}-valued when $\pi^{A'}$ is {\bf R}-valued), is
$$\pi^{A'}\pi^{B'}\nabla^A_{B'}\gamma_A = \gamma_A\pi^{B'}\nabla^A_{B'}\pi^{A'}.\eqno(6.2.16{\rm b})$$

Proof. This lemma amounts to a classical integrability result, which I demonstrate by a variation of the proof in [23] (7.3.20). When $\pi^{A'}$ is {\bf R} valued, $Z_{[\pi]}$ is an integrable distribution of $TM$ and the Frobenius theorem provides local coordinates $(x^1,\ldots,x^4)$, say, such that $Z_{[\pi]} = \langle \partial_{x^1},\partial_{x^2} \rangle_{\bf R}$. When $\pi^{A'}$ is {\bf C} valued, $Z_{[\pi]}$ is, as described in (6.2.5), a two-(complex)-dimensional distribution in ${\bf C}TM$ and the Newlander-Nirenberg theorem provides local holomorphic coordinates $(z^1,z^2)$, say, such that $Z_{[\pi]} = \langle \partial_{z^1},\partial_{z^2} \rangle_{\bf C}$. To treat the two cases together, for $i=1$, 2, let $\partial_i$ denote $\partial_{x^i}$ and $\partial_{z^i}$ in the real and complex cases respectively, and let $X^i$ stand for $x^i$ in the real case and $z^i$ in the complex case. Then, $\partial_i = \zeta_i^A\pi^{A'}$, for some spinors $\zeta^A_i$, which must be linearly independent. (6.2.16a) is equivalent to the pair of PDEs: $\partial_i\phi = \zeta^A_i\gamma_A$, $i=1$, 2; i.e., to the single equation $\diamondsuit\phi = \psi$, where $\diamondsuit$ denotes the exterior derivative in the two coordinates $(X^1,X^2)$, i.e., the exterior derivative on integral surfaces of $Z_{[\pi]}$, in the real case and the $\partial$ operator in the complex case, and where $\psi = (\zeta^A_1\gamma_A)dX^1 + (\zeta^A_2\gamma_A)dX^2$. Thus solvability of the equation reduces to the relevant Poincar\'e lemma (in the complex case, the conjugate of the Dolbeault-Grothendieck lemma), i.e., to the $\diamondsuit$-closure of $\psi$, i.e., to $\partial_1(\zeta^A_2\gamma_A) = \partial_2(\zeta^A_1\gamma_A)$. Write out this condition as
$$(\zeta_1^B\pi^{B'}\nabla_b\zeta^A_2)\gamma_A + (\zeta_1^B\pi^{B'}\nabla_b\gamma_A)\zeta^A_2 = (\zeta_2^B\pi^{B'}\nabla_b\zeta^A_1)\gamma_A + (\zeta_2^B\pi^{B'}\nabla_b\gamma_A)\zeta^A_1$$
and then rearrange terms to obtain, with $\zeta_2 \cdot \zeta_1 := \zeta^A_2\zeta^B_1\epsilon_{BA}$,
$$\eqalignno{(\zeta_1^B\pi^{B'}\nabla_b\zeta^A_2 - \zeta_2^B\pi^{B'}\nabla_b\zeta^A_1)\gamma_A &= \zeta^A_1\zeta_2^B\pi^{B'}\nabla_b\gamma_A - \zeta^A_2\zeta_1^B\pi^{B'}\nabla_b\gamma_A\cr
&= (\zeta^A_1\zeta^B_2 - \zeta_2^A\zeta_1^B)\pi^{B'}\nabla_b\gamma_A\cr
&= (\zeta_2 \cdot \zeta_1)\epsilon^{AB}\pi^{B'}\nabla_b\gamma_A\cr
&= (\zeta_2 \cdot \zeta_1)\pi^{B'}\nabla^A_{B'}\gamma_A.&(*)\cr}$$
But,
$$\eqalign{0 &= [\zeta_1^B\pi^{B'},\zeta_2^C\pi^{C'}]^a\cr 
&= \zeta^B_1\pi^{B'}\nabla_b(\zeta^A_2\pi^{A'}) - \zeta^B_2\pi^{B'}\nabla_b(\zeta^A_1\pi^{A'})\cr
&= (\zeta^B_1\pi^{B'}\nabla_b\zeta^A_2 - \zeta^B_2\pi^{B'}\nabla_b\zeta^A_1)\pi^{A'} + (\zeta_2^A\zeta_1^B - \zeta_1^A\zeta_2^B)\pi^{B'}\nabla_b\pi^{A'}\cr
&= (\zeta^B_1\pi^{B'}\nabla_b\zeta^A_2 - \zeta^B_2\pi^{B'}\nabla_b\zeta^A_1)\pi^{A'} + (\zeta_1 \cdot \zeta_2)\epsilon^{AB}\pi^{B'}\nabla_b\pi^{A'}\cr
&= (\zeta^B_1\pi^{B'}\nabla_b\zeta^A_2 - \zeta^B_2\pi^{B'}\nabla_b\zeta^A_1)\pi^{A'} + (\zeta_1 \cdot \zeta_2)\pi^{B'}\nabla^A_{B'}\pi^{A'}.\cr}$$
Transvecting through the last equation by $\gamma_A$ yields
$$(\zeta^B_1\pi^{B'}\nabla_b\zeta^A_2 - \zeta^B_2\pi^{B'}\nabla_b\zeta^A_1)\gamma_A\pi^{A'} = (\zeta_2 \cdot \zeta_1)\gamma_A\pi^{B'}\nabla^A_{B'}\pi^{A'}.$$
Thus, the integrability condition $(*)$ is equivalent to
$$(\zeta_2 \cdot \zeta_1)\gamma_A\pi^{B'}\nabla^A_{B'}\pi^{A'} = (\zeta_2 \cdot \zeta_1)\pi^{A'}\pi^{B'}\nabla^A_{B'}\gamma_A,$$
i.e., to (6.2.16b).
(As an aside, (6.2.16b) can be directly obtained as a necessary condition as follows: 
$$\eqalign{\pi^{A'}\pi^{B'}\nabla^A_{B'}\gamma_A &= \pi^{A'}\pi^{B'}\nabla^A_{B'}(\pi^{C'}\nabla_{AC'}\phi)\cr
&= \pi^{A'}(\pi^{B'}\nabla^A_{B'}\pi^{C'})\nabla_{AC'}\phi + \pi^{A'}\pi^{B'}\pi^{C'}\nabla^A_{B'}\nabla_{AC'}\phi\cr
&= \pi^{A'}(\eta^A\pi^{C'})\nabla_{AC'}\phi + 0\qquad\hbox{since the second summand is skew in the subcripts $B'$ and $C'$}\cr
&= \eta^A\pi^{A'}\gamma_A\cr
&= \gamma_A\pi^{B'}\nabla^A_{B'}\pi^{A'}.)\cr}$$
\vskip 24pt
\noindent {\bf 6.2.17 Generalized Goldberg-Sachs Theorem (GGST) }\hfil\break
On $(M,g)$, with $\tilde\Psi_{A'B'C'D'}$ nonzero, consider the following three statements for a (possibly local) smooth {\bf K}-valued projective spinor field $[\pi^{A'}]$:\hfil\break
(i) $[\pi^{A'}]$ is a WPS of multiplicity $p$ ($2 \leq p \leq 4$);\hfil\break
(ii) each LSR $\pi^{A'}$ of $[\pi^{A'}]$ satisfies (6.2.1);\hfil\break
(iii) for each LSR $\pi^{A'}$ of $[\pi^{A'}]$, 
$$\underbrace{\pi^{A'}\ldots\pi^{C'}}_{5-q}\nabla^{DD'}\tilde\Psi_{A'B'C'D'} = 0.$$
Then:
$$\displaylines{({\rm i})+({\rm ii})\ \Rightarrow\ ({\rm iii})\hbox{ (with } $q=p$)\cr
({\rm i})+({\rm iii})\hbox{ (with } $q=p$)\ \Rightarrow\ ({\rm ii})\cr
({\rm ii})+({\rm iii})\hbox{ (with } $q=2$)\ \Rightarrow\ ({\rm i})\hbox{ (for some $p$, $2 \leq p \leq 4$)}.\cr}$$
Proof. Condition (i) entails $0 = \tilde\Psi_{A'B'C'D'}\pi^{A'}\pi^{B'}\pi^{C'}$, whence
$$0 = \nabla^{DD'}(\tilde\Psi_{A'B'C'D'}\pi^{A'}\pi^{B'}\pi^{C'}) = \pi^{A'}\pi^{B'}\pi^{C'}\nabla_{DD'}\tilde\Psi_{A'B'C'D'} + 3\tilde\Psi_{A'B'C'D'}\pi^{A'}\pi^{B'}\nabla_{DD'}\pi^{C'}.\eqno(1)$$
For convenience, take $p=2$ in (i); the following argument requires only simple modification for $p=3$, 4. With $p=2$, (i) is equivalent to
$$\tilde\Psi_{A'B'C'D'} = \cases{\Psi\pi_{(A'}\pi_{B'}\beta_{C'}\gamma_{D')},&$\pi^{A'}$ real;\cr \Psi\pi_{(A'}\pi_{B'}\overline\pi_{C'}\overline\pi_{D')},&$\pi^{A'}$ complex;\cr}\eqno(2)$$
for some nonzero real scalar $\Psi$, where in the real case real spinors $\beta_{A'}$ and $\gamma_{A'}$ satisfy $\pi^{D'}\beta_{D'} \not= 0$, $\pi^{D'}\gamma_{D'} \not= 0$, while in the complex case $\pi^{D'}\overline\pi_{D'} \not= 0$. Note that a complex $\pi^{A'}$ can only be a multiple WPS with $p=2$. From (2),
$$\tilde\Psi_{A'B'C'D'}\pi^{A'}\pi^{B'} = \cases{{\Psi \over 6}(\pi^{E'}\beta_{E'})(\pi^{F'}\gamma_{F'})\pi_{C'}\pi_{D'},&$\pi^{A'}$ real;\cr {\Psi \over 6}(\pi^{E'}\overline\pi_{E'})^2\pi_{C'}\pi_{D'},&$\pi^{A'}$ complex;\cr}$$
in either case, $\tilde\Psi_{A'B'C'D'}\pi^{A'}\pi^{B'} = \lambda\pi_{C'}\pi_{D'}$, for some {\sl nonzero} (real) scalar $\lambda$. Substituting into (1) yields
$$0 = \pi^{A'}\pi^{B'}\pi^{C'}\nabla_{DD'}\tilde\Psi_{A'B'C'D'} + 3\lambda\pi_{C'}\pi_{D'}\nabla^{DD'}\pi^{C'}.$$
The first two implications of the theorem follow immediately.

Now consider the third implication. By (6.2.9), (ii) entails $[\pi^{A'}]$ is a WPS. Hence, $\tilde\Psi_{A'B'C'D'}\pi^{A'}\pi^{B'}\pi^{C'} = \lambda\pi_{D'}$, for some $\lambda$; the aim is to prove $\lambda = 0$. Differentiation of the last equation yields
$$\pi^{A'}\pi^{B'}\pi^{C'}\nabla^{DD'}\tilde\Psi_{A'B'C'D'} + 3\tilde\Psi_{A'B'C'D'}\pi^{A'}\pi^{B'}\nabla^{DD'}\pi^{C'} = \pi_{D'}\nabla^{DD'}\lambda + \lambda\nabla^{DD'}\pi_{D'}.\eqno(3)$$
The first term on the left-hand side of (3) vanishes by (iii)($q=2$). Writing $\tilde\Psi_{A'B'C'D'} = \pi_{(A'}\alpha_{B'}\beta_{C'}\gamma_{D')}$, for some $\alpha_{A'}$, $\beta_{A'}$, and $\gamma_{A'}$, then
$$\lambda = {\pi^{A'}\pi^{B'}\pi^{C'}\alpha_{(A'}\beta_{B'}\gamma_{C')} \over 4},\eqno(4)$$
and the second term on the left-hand side of (3) is
$$\displaylines{3\left[{6 \over 24}\bigl(\alpha_{(A'}\beta_{B'}\gamma_{D')}\pi_{C'} + \alpha_{(A'}\beta_{B'}\gamma_{C')}\pi_{D'}\bigr)\right]\pi^{A'}\pi^{B'}\nabla^{DD'}\pi^{C'}\hfill\cr
\hfill = {3 \over 4}\left[\pi^{A'}\pi^{B'}\pi^{D'}\omega^D\alpha_{(A'}\beta_{B'}\gamma_{D')} - \pi^{A'}\pi^{B'}\pi^{C'}\eta^D\alpha_{(A'}\beta_{B'}\gamma_{C')}\right]= 3\lambda[\omega^D - \eta^D]\hfill\cr}$$
having used (6.2.6) for the first equality and (4) for the second. Substituting this expression into (3), and using (6.2.11) to rewrite the second term of the right-hand side of (3), gives
$$3\lambda(\omega^D - \eta^D) = \pi_{D'}\nabla^{DD'}\lambda - \lambda(\omega^D+ \eta^D).$$
Allowing for the possibility that $\lambda$ may be {\bf C}-valued, if $\lambda$ is nonzero at some point, one can choose a branch of $\ln\lambda$ on some neighbourhood of that point. On that neighbourhood, one can rewrite the last equation as
$$\pi^{D'}\nabla_{DD'}(\ln\lambda) = 2\eta_D - 4\omega_D = 6\eta_D - 4\nabla_{DD'}\pi^{D'}\eqno(5)$$
on using (6.2.11) again. But (5) is of the form (6.2.16a) and must therefore satisfy (6.2.16b) with $\gamma_A := 6\eta_A - 4\nabla_{AD'}\pi^{D'}$, i.e., 
$$\eqalign{\pi^{A'}\pi^{B'}\nabla^A_{B'}(6\eta_A - 4\nabla_{AD'}\pi^{D'}) &= (6\eta_A - 4\nabla_{AD'}\pi^{D'})\pi^{B'}\nabla^A_{B'}\pi^{A'}\cr
&= (6\eta_A - 4\nabla_{AD'}\pi^{D'})\eta^A\pi^{A'}\qquad\hbox{by (6.2.6)}\cr
&= -4\pi^{A'}\eta^A\nabla_{AD'}\pi^{D'}.\cr}$$
Using the definitions in (6.2.14--15), the last equation can be written
$$-6\lambda^{A'} = 6\pi^{A'}\pi^{B'}\nabla^A_{B'}\eta_A 
= 4\pi^{A'}\pi^{B'}\nabla^A_{B'}\nabla_{AD'}\pi^{D'} - 4\pi^{A'}\eta^A\nabla_{AD'}\pi^{D'} = 4\mu^{A'}.$$
By (6.2.15), it follows that $0 = \lambda_{A'}$, whence by (6.2.13) $0 = \tilde\Psi_{A'B'C'D'}\pi^{B'}\pi^{C'}\pi^{D'}$, i.e., $\lambda = 0$ after all and $[\pi^{A'}]$ is indeed a {\sl multiple} WPS.
\vskip 24pt
\noindent {\bf 6.2.18 Corollary}\hfil\break
For any {\bf K}-valued solution $[\pi^{A'}]$ of (6.2.1), the condition $\Phi_{ABA'B'}\pi^{A'}\pi^{B'} = 0$ implies $[\pi^{A'}]$ is a multiple WPS.

Proof. Observe that, for any solution of (6.2.1),
$$\eqalign{\pi^{C'}\nabla^A_{C'}(\Phi_{ABA'B'}\pi^{A'}\pi^{B'}) &= (\pi^{C'}\nabla^A_{C'}\Phi_{ABA'B'})\pi^{A'}\pi^{B'} + 2(\pi^{C'}\nabla^A_{C'}\pi^{A'})\Phi_{ABA'B'}\pi^{B'}\cr
&= (\nabla^{D'}_B\tilde\Psi_{A'B'C'D'})\pi^{A'}\pi^{B'}\pi^{C'} + 2\eta^A\Phi_{ABA'B'}\pi^{A'}\pi^{B'},\cr}$$
upon using a spinor Bianchi identity (A.7) and (6.2.6). Hence, $\Phi_{ABA'B'}\pi^{A'}\pi^{B'} = 0$ implies condition (iii) of GGST (with $q=2$), whence $[\pi^{A'}]$ is a multiple WPS by the GGST.
\vskip 24pt
\noindent {\bf 6.2.19 Corollary}\hfil\break
If $(M,g)$ is Einstein, (i) and (ii) of the GGST are equivalent, i.e., the multiple WPSs of $\tilde\Psi_{A'B'C'D'}$ are precisely the solutions of (6.2.1). Explicitly: an Einstein $(M,g)$ admits an integrable $\alpha$-distribution $Z_{[\pi]}$ iff $[\pi^{A'}]$ is a (real) multiple WPS; and admits a complex solution $[\pi^{A'}]$ of (6.2.1) iff $[\pi^{A'}]$ is a multiple WPS, i.e., iff $\tilde\Psi_{A'B'C'D'} \propto \pi_{(A'}\pi_{B'}\overline\pi_{C'}\overline\pi_{D')}$ ([16]).

Proof. The Einstein condition is $\Phi_{ABA'B'} = 0$, in which case the Bianchi equation (A.7) entails (iii) of the GGST is trivial.
\vskip 24pt
\noindent {\bf 6.2.20 Remark}\hfil\break
If a compact Einstein $(M,g)$ admits a solution of (6.2.1), real or complex, then $\tilde\Psi_{A'B'C'D'}$ is algebraically special as just noted; in particular, $\tilde\Psi_{A'B'C'D'}$ is not of type $\{1\overline 111\}$Ib, whence $\chi(M) \leq 3\tau(M)/2$ ([15], [20]). If $\Psi_{ABCD}$ is not of type $\{1\overline 111\}$Ib either, then in fact $\chi(M) \leq -3\vert \tau(M)\vert/2 \leq 0$.
\vskip 24pt
In Lorentzian geometry there is a significant connexion between solutions of (6.2.1) and massless fields. Though the physical interpretation of massless fields is lacking, the analogous constructions in the neutral context may still be of geometrical interest. These constructions carry over from [23] routinely but I include them for completeness.
\vskip 24pt
\noindent {\bf 6.2.21 Definition}\hfil\break
A {\sl primed massless field} of rank $n$ on $(M,g)$ is a smooth section $\phi_{A'B'\ldots L'}$ of 
$$\underbrace{{\bf C}S' \odot \cdots \odot {\bf C}S'}_n\qquad\hbox{satisfying}\qquad \nabla^{AA'}\phi_{A'B'\ldots L'} = 0.$$
Such fields may, of course, be real valued. One may define {\sl unprimed massless fields} similarly.
\vskip 24pt
\noindent {\bf 6.2.22 The Buchdahl Constraint}\hfil\break
H. A. Buchdahl [3] discovered an algebraic consistency condition for the existence of massless fields, derived as follows:
$$\eqalign{0 &= \nabla^{B'}_A\nabla^{AA'}\phi_{A'B'\ldots L'}\cr
&= (\nabla^{(B'}_A\nabla^{A')A} + \nabla^{[B'}_A\nabla^{A']A})\phi_{A'B'\ldots L'}\cr
&= (\Square^{B'A'} + {1 \over 2}\epsilon^{B'A'}\Square)\phi_{A'B'\ldots L'}\qquad\hbox{by (A.5)}\cr
&= \Square^{A'B'}\phi_{A'B'\ldots L'}\cr
&= -\tilde X^{A'B'}{}_{A'}{}^{M'}\phi_{M'B'\ldots L'} - \tilde X^{A'B'}{}_{B'}{}^{M'}\phi_{A'M'\ldots L'} - \tilde X^{A'B'}{}_{C'}{}^{M'}\phi_{A'B'M'\ldots L'} \cdots - \tilde X^{A'B'}{}_{L'}{}^{M'}\phi_{A'B'\ldots M'},\cr}$$
by [22] (4.9.12). The first two terms of the last expression vanish by [22] (4.6.6). By [22] (4.6.35), one then obtains, for $n \geq 2$,
$$0 = (n-2)\phi_{A'B'C'(D'\ldots K'}\tilde\Psi_{L')}{}^{A'B'C'}.$$
\vskip 24pt
\noindent {\bf 6.2.23 Lemma}\hfil\break
If $\pi^{A'}$ is a principal spinor (PS, see [22], p. 162; [16]) of multiplicity $k > 1$ of a massless field $\phi_{A'B'\ldots L'}$, then $\pi^{A'}$ solves (6.2.1). (The result is local.)

Proof. By assumption:
$$\phi_{A'\ldots G'H'\ldots L'}\underbrace{\pi^{H'}\ldots\pi^{L'}}_{n-k} = \nu\pi_{A'}\ldots\pi_{G'},\qquad\nu \not= 0, \hskip 1in \phi_{A'\ldots F'G'\ldots L'}\underbrace{\pi^{G'}\ldots\pi^{L'}}_{n-k+1} = 0.$$
Hence,
$$\eqalign{0 &= \nabla^{FF'}(\phi_{A'\ldots F'G'\ldots L'}\pi^{G'}\ldots\pi^{L'})\qquad\hbox{(note that $n-(n-k+1) = k-1 \geq 1$)}\cr
&= \phi_{A'\ldots F'G'\ldots L'}\nabla^{FF'}(\pi^{G'}\ldots\pi^{L'})\cr
&= (n-k+1)\phi_{A'\ldots F'G'\ldots L'}\pi^{H'}\ldots\pi^{L'}\nabla^{FF'}\pi^{G'}\cr
&= \nu(n-k+1)\pi_{A'}\ldots\pi_{G'}\nabla^{FF'}\pi^{G'}.\cr}$$
Noting that $n - (n-k) = k \geq 2$, (6.2.1) follows.
\vskip 24pt
\noindent {\bf 6.2.24 Remark}\hfil\break
Hence, if one can construct massless fields with a multiple PS $\pi^{A'}$ on $(M,g)$, then $\pi^{A'}$ solves (6.2.1), whence is a WPS by (6.2.9), and multiple according to the dictates of the GGST. A particular case is the following result.
\vskip 24pt
\noindent {\bf 6.2.25 Lemma}\hfil\break
If $(M,g)$ admits a null massless field ([22], [16]) $\phi_{A'B'\ldots L'} = \chi\pi_{A'}\pi_{B'}\ldots\pi_{L'}$, of rank $n > 2$, then $[\pi^{A'}]$ is a multiple WPS.

Proof. In this case, the Buchdahl constraint is $0 = \pi_{A'}\pi_{B'}\pi_{C'}\pi_{(D'}\ldots\pi_{K'}\tilde\Psi_{L')}{}^{A'B'C'}$ which, by [22] (3.5.15), is equivalent to $\tilde\Psi_{A'B'C'D'}\pi^{A'}\pi^{B'}\pi^{C'} = 0$.
\vskip 24pt
A converse of (6.2.25) is provided by the analogue of results of Robinson and Sommers, see [23] (7.3.14).
\vskip 24pt
\noindent {\bf 6.2.26 Proposition}\hfil\break
Let $\pi^{A'}$ be a {\bf K}-valued solution of (6.2.1). If $(n-2)\tilde\Psi_{A'B'C'D'}\pi^{A'}\pi^{B'}\pi^{C'} = 0$, then there exists a null massless field $\phi_{A'B'\ldots L'} = \chi\pi_{A'}\pi_{B'}\ldots\pi_{L'}$, $\chi$ a smooth function (positive when $\pi^{A'}$ is real), of rank $n$. The freedom in $\chi$ is multiplication by $\exp(f)$, where $f$ is any smooth function satisfying $\pi^{B'}\nabla_{BB'}f = 0$. When $\pi^{A'}$ is real, $f$ is constant on $\alpha$-surfaces. When $\pi^{A'}$ is complex, all (complex) gradients from $H$ of $f$ vanish; in particular, $\partial_{z^i} f = 0$, $i=1$, 2, i.e., $f$ is anti-holomorphic with respect to the induced complex structure described in (6.2.5).

Proof. Take $\phi_{A'B'\ldots L'}$ as defined with $\chi$ to be determined. The massless field equations dictate
$$\eqalign{0 = \nabla^{AA'}(\chi\pi_{A'}\ldots\pi_{L'}) &= \pi_{B'}\ldots\pi_{L'}\pi_{A'}\nabla^{AA'}\chi + \chi\pi_{B'}\ldots\pi_{L'}\nabla^{AA'}\pi_{A'}\cr 
&+ \chi \pi_{C'}\ldots\pi_{L'}\pi_{A'}\nabla^{AA'}\pi_{B'} + \ldots + \chi\pi_{B'}\ldots\pi_{K'}\pi_{A'}\nabla^{AA'}\pi_{L'}\cr
&= \pi_{B'}\ldots\pi_{L'}\pi_{A'}\nabla^{AA'}\chi + \chi\bigl(\nabla^{AA'}\pi_{A'} - (n-1)\eta^A\bigr)\pi_{B'}\ldots\pi_{L'}\cr}$$
upon using (6.2.6). Hence,
$$\pi^{A'}\nabla_{AA'}(\ln\chi) = -\nabla_{AA'}\pi^{A'} - (n-1)\eta_A.\eqno(*)$$
To determine whether this equation is solvable, one checks the integrability condition of (6.2.16), which here is:
$$\pi^{A'}\pi^{B'}\nabla^A_{B'}\bigl(-\nabla_{AD'}\pi^{D'} - (n-1)\eta_A\bigr) = -\bigl(\nabla_{AD'}\pi^{D'} + (n-1)\eta_A\bigr)\pi^{B'}\nabla^A_{B'}\pi^{A'}.$$
By (6.2.6), the right-hand side simplifies to $-\pi^{A'}\eta^A\nabla_{AD'}\pi^{D'}$, whence one can write the integrability condition as
$$-\pi^{A'}\eta^A\nabla_{AD'}\pi^{D'} = -\pi^{A'}\pi^{B'}\nabla^A_{B'}\nabla_{AD'}\pi^{D'} - (n-1)\pi^{A'}\pi^{B'}\nabla^A_{B'}\eta_A,$$
i.e., as
$$\pi^{A'}\pi^{B'}\nabla^A_{B'}\nabla_{AD'}\pi^{D'} - \pi^{A'}\eta^A\nabla_{AD'}\pi^{D'} = (n-1)\pi^{A'}\pi^{B'}\nabla_{AB'}\eta^A.$$
By (6.2.15), the left-hand side is $\mu^{A'}$ while, by (6.2.14), the right-hand side is $(n-1)\lambda^{A'}$. Thus, the integrability condition is $\mu^{A'} = (n-1)\lambda^{A'}$, which by (6.2.14) amounts to $(n-2)\tilde\Psi_{A'B'C'D'}\pi^{B'}\pi^{C'}\pi^{D'}$. Granted this condition, $(*)$ admits a solution $F$, say, whence $\chi = \exp(F)$. The freedom in $F$ is $F \mapsto F+f$, where $\pi^{B'}\nabla_{BB'}f = 0$, i.e., the freedom in $\chi$ is multiplication by $\exp(f)$.
\vskip 24pt
\noindent {\bf 6.2.27 Remark}\hfil\break
Guillemin and Sternberg [10] characterized null rank two massless fields on the conformal compactification of ${\bf R}^{2,2}$ while Hughston and Mason [12] generalized (6.2.25--26) to the conformal compactification of ${\bf C}^{2n}_{\rm E}$.
\vskip 24pt
I will now largely concentrate on {\bf R}-valued solutions of (6.2.1); the treatment of {\bf C}-valued solutions on the basis of (6.2.5) will be presented elsewhere. Therefore, let $[\pi^{A'}]$ be a (local) section of ${\bf P}S'$whose corresponding $\alpha$-distribution $Z_{[\pi]}$ is integrable. Then, as noted previously: any LSR satisfies (6.2.1); $Z_{[\pi]}$ is also auto-parallel; with respect to the spin frames (6.2.2), $\tilde\kappa = \tilde\sigma = 0$ (6.2.4); $S^a$, $T^a \in Z_{[\pi]}$. In general
$$o_{B'}\nabla_a o^{B'} = \tilde\tau\ell_a - \tilde\rho\tilde m_a + \tilde\kappa n_a - \tilde\sigma m_a \hskip 1in o^{B'}\nabla_{AB'}o_{A'} = \tilde\beta\ell_a - \tilde\epsilon\tilde m_a - \tilde\kappa n_a + \tilde\sigma m_a\eqno(6.2.28)$$
so, with respect to the spin frames (6.2.2)
$$\omega_A = \tilde\tau o_A - \tilde\rho\iota_A \hskip .75in \eta_A = \tilde\beta o_A - \tilde\epsilon\iota_A,\qquad\hbox{whence}\qquad \omega^D\eta_D = \tilde\tau\tilde\epsilon - \tilde\rho\tilde\beta.\eqno(6.2.29)$$
Under change of LSR, $\pi^{A'} \mapsto \tilde\lambda\pi^{A'}$,
$$S_a \mapsto \tilde\lambda^2 S_a,\hskip .5in \omega_A \mapsto \tilde\lambda\omega_A \hskip .5in T_a \mapsto \tilde\lambda^2 T_a + \tilde\lambda\pi_{A'}\pi^{B'}\nabla_{AB'}\tilde\lambda \hskip .5in \eta_A \mapsto \tilde\eta_A + \pi^{B'}\nabla_{AB'}\tilde\lambda.\eqno(6.2.30)$$
It will also be convenient to record the behaviour of relevant spin coefficients under the change (6.2.2) in spin frames:
$$\displaylines{\pmatrix{\tilde\kappa\cr \tilde\sigma\cr} \mapsto \tilde\lambda^3\left({^\tau\! A}\right)\pmatrix{\tilde\kappa\cr \tilde\sigma\cr} \hskip 1in \pmatrix{\tilde\rho\cr \tilde\tau\cr} \mapsto \tilde\lambda\left({^\tau\! A}\right)\pmatrix{\tilde\rho\cr \tilde\tau\cr} + \tilde\lambda^2\tilde\mu\left({^\tau\! A}\right)\pmatrix{\tilde\kappa\cr \tilde\sigma\cr}\cr
\noalign{\vskip -6pt}
\hfill\llap(6.2.31)\cr
\noalign{\vskip -6pt}
\pmatrix{\tilde\epsilon\cr \tilde\beta\cr} \mapsto \tilde\lambda\left({^\tau\! A}\right)\pmatrix{\tilde\epsilon\cr \tilde\beta\cr} - \tilde\lambda^2\tilde\mu\left({^\tau\! A}\right)\pmatrix{\tilde\kappa\cr \tilde\sigma\cr} + \left({^\tau\! A}\right)\pmatrix{D\tilde\lambda\cr \triangle\tilde\lambda\cr},\cr}$$
the first and second pair being of the same form as in (6.1.6). Notice that the vanishing of $\tilde\kappa$ and $\tilde\sigma$ is invariant, as expected; the vanishing of $\tilde\rho$ and $\tilde\tau$ is invariant when $\tilde\kappa = \tilde\sigma = 0$, whence in the present circumstances; but the vanishing of $\tilde\beta$ and $\tilde\epsilon$, even when $\tilde\kappa = \tilde\sigma = 0$, is only invariant under rescalings $\tilde\lambda$ which are constant on $\alpha$-surfaces. These facts reflect the scaling behaviour recorded in (6.2.30).

With this information at hand, now consider the significance of the quantities $T_a$ and $S_a$. The behaviour of $T_a$ in (6.2.30) indicates its vanishing is related to a choice of LSR rather than the underlying geometry of the $\alpha$-distribution. Indeed:
\vskip 24pt
\noindent {\bf 6.2.32 Proposition}\hfil\break
For a {\bf K}-valued solution $[\pi^{A'}]$ of (6.2.1), there is an LSR satisfying $T_a := \pi^{B'}\nabla_{AB'}\pi_{A'} = 0$ iff $\tilde\Psi_{A'B'C'D'}\pi^{B'}\pi^{C'}\pi^{D'} = 0$, i.e., iff $[\pi^{A'}]$ is a multiple WPS.

Proof. For a solution of (6.2.1), $T_a = 0$ iff $\eta_A = 0$. Thus, $T_a = 0$ implies $\tilde\Psi_{A'B'C'D'}\pi^{B'}\pi^{C'}\pi^{D'} = 0$ by (6.2.14).

Conversely, given any LSR $\pi^{A'}$, there is an LSR for which $T_a = 0$ iff there is a function $f$ satisfying $0 = (f\pi^{B'})\nabla_{AB'}(f\pi_{A'})$, i.e., iff $0 = f\eta_A\pi_{A'} + \pi^{B'}(\nabla_{AB'}f)\pi_{A'}$, i.e., iff $\pi^{B'}\nabla_{AB'}(\ln f) = -\eta_A$. By (6.2.16), the integrability condition for this equation is $\pi^{A'}\pi^{B'}\nabla^A_{B'}(-\eta_A) = -\eta_A\pi^{B'}\nabla^A_{B'}\pi^{A'}$. The right-hand side is $-\eta_A\eta^A\pi^{A'} = 0$, while the left-hand side, by (6.2.14), is just $\tilde\Psi^{A'}{}_{B'C'D'}\pi^{B'}\pi^{C'}\pi^{D'}$.
\vskip 24pt
Restricting attention again to real solutions of (6.2.1), recall that (6.2.1) implies $\alpha$-surfaces are totally geodesic. A more geometrical condition one might consider than the vanishing of $T_a$ (i.e., the existence of an LSR parallel on $\alpha$-surfaces) is the flatness of the induced connection within $\alpha$-surfaces. The induced connection is determined by $D$ and $\triangle$. Since the unprimed spin frame in (6.2.2) is arbitrary, any local frame field of $Z_{[\pi]}$ can be represented as $\{\ell^a,\tilde m^a\}$ for an appropriate choice of unprimed spin frame and LSR $\pi^{A'}$ (see [17], 2.4). Hence, $D\ell^a = D\tilde m^a = \triangle\ell^a = \triangle\tilde m^a = 0$ for some such choice is equivalent to flatness of the induced connection. Given that $\tilde\kappa = \tilde\sigma = 0$ already and noting (2.11), from (2.10) $\ell^a$ and $\tilde m^a$ are parallel on $\alpha$-surfaces iff
$$\kappa = \rho = \alpha = \epsilon = \tau' = \sigma' = \tilde\epsilon = \tilde\beta = 0.\eqno(6.2.33)$$
From (3.4) one deduces: $\Phi_{00} = 0$ from (a); $\Phi_{10} = 0$ from (h); $\Phi_{20} = 0$ from $({\rm e}')$; $\tilde\Psi_1 = 0$ from (\~i); and we already know $\tilde\Psi_0 = 0$ from (\~ b). As $\tilde\Psi_{A'B'C'D'}\pi^{A'}\pi^{B'}\pi^{C'}\pi^{D'} = 0$, one expects the invariant expressions for the condition for flatness of the induced connection on the $\alpha$-surfaces to be:
$$\Phi_{ABA'B'}\pi^{A'}\pi^{B'} = 0; \hskip 1.25in \tilde\Psi_{A'B'C'D'}\pi^{B'}\pi^{C'}\pi^{D'} = 0.\eqno(6.2.34)$$
These conditions hold for any Walker metric, as does (6.2.33) with respect to the Walker spin frames.

One can confirm that (6.2.34) does characterize flatness of the induced connection by computing the components of the curvature of $(M,g)$ which determine the curvature of the induced connection. Using the frame $\{\ell^a,\tilde m^a\}$ for the $\alpha$-distribution, $\{n_a,m_a\}$ is a frame of one-forms on an $\alpha$-surface $\cal Z$. The components of $\nabla_b\ell^a$ and of $\nabla_b\tilde m^a$ on $\cal Z$, i.e., $D\ell^a$, $\triangle\ell^a$, $D\tilde m^a$ and $\triangle\tilde m^a$ may be written
$$\displaylines{\left[(\epsilon+\tilde\epsilon)\ell^a + \kappa\tilde m^a\right]n_b - \left[(\alpha+\tilde\beta)\ell^a+\rho\tilde m^a\right]m_b\cr
\left[-\tau'\ell^a + (\tilde\epsilon - \epsilon)\tilde m^a\right]n_b - \left[\sigma'\ell^a + (\tilde\beta-\alpha)\tilde m^a\right]m_b.\cr}$$
Using the convention in [22], \S 4.13, for the staggering of indices on the connection one-forms, one has
$$(\omega_{\bf j}{}^{\bf i}) = \pmatrix{(\epsilon+\tilde\epsilon)n_b-(\alpha+\tilde\beta)m_b&-\tau'n_b-\sigma'm_b\cr \kappa n_b - \rho m_b&(\tilde\epsilon-\epsilon)n_b - (\tilde\beta-\alpha)m_b\cr}$$
where superscripts label rows and subscripts columns. Hence,
$$(\omega_{\bf j}{}^{\bf i}) \wedge (\omega_{\bf k}{}^{\bf j})= \pmatrix{(\tau'\rho+\sigma'\kappa)n \wedge m&-2(\sigma'\epsilon + \tau'\alpha)n \wedge m\cr -2(\kappa\alpha - \epsilon\rho)n \wedge m&-(\kappa\sigma' + \rho\tau')n \wedge m\cr}.$$
Using (4.1), and noting (6.1.24), on $\cal Z$ one computes:
$$\eqalign{d\omega_1{}^1 &= \left[D(\epsilon+\tilde\epsilon)n_a - \triangle(\epsilon+\tilde\epsilon)m_a\right]\wedge n_b + (\epsilon+\tilde\epsilon)dn_b - \left[D(\alpha+\tilde\beta)n_a - \triangle(\alpha+\tilde\beta)m_a\right] \wedge m_b - (\alpha+\tilde\beta)dm_b\cr
&= \left[(D\alpha - \triangle\epsilon) + (D\tilde\beta - \triangle\tilde\epsilon) + (\epsilon + \tilde\epsilon)(\alpha + \tau' + \tilde\beta) - (\alpha + \tilde\beta)(\gamma' - \rho + \tilde\epsilon)\right]m \wedge n;\cr
d\omega_1{}^2 &= \left[(D\kappa)n_a - (\triangle\kappa) m_a\right] \wedge n_b + \kappa dn_b - \left[(D\rho) n_a - (\triangle\rho) m_b\right] \wedge m_b - \rho dm_b\cr
&= \left[D\rho - \triangle\kappa + \kappa(\alpha+ \tau'+\tilde\beta) - \rho(\tilde\epsilon - \epsilon - \rho)\right]m \wedge n\cr
d\omega_2{}^1 &= -\left[(D\tau')n_a - (\triangle\tau')m_a\right] \wedge n_b - \tau'dn_b - \left[(D\sigma')n_a - (\triangle\sigma')m_a\right] \wedge m_b - \sigma' dm_b\cr
&= \left[\triangle\tau' + D\sigma' - \tau'(\alpha+\tau'+\tilde\beta) - \sigma'(\tilde\epsilon - \epsilon - \rho)\right]m \wedge n\cr
d\omega_2{}^2 &= \left[D(\tilde\epsilon-\epsilon)n_a - \triangle(\tilde\epsilon-\epsilon)m_a\right] \wedge n_b + (\tilde\epsilon-\epsilon)dn_b - \left[D(\tilde\beta-\alpha)n_a - \triangle(\tilde\beta-\alpha)m_a\right]\wedge m_b - (\tilde\beta-\alpha)dm_b\cr
&= \left[\triangle\epsilon - D\alpha + D\tilde\beta - \triangle\tilde\epsilon + (\tilde\epsilon - \epsilon)(\alpha + \tau' + \tilde\beta) - (\tilde\beta - \alpha)(\tilde\epsilon - \epsilon - \rho)\right]m \wedge n\cr}$$
Using (3.4), the second structure equation for the induced connection then gives, for the curvature of the induced connection:
$$\eqalign{\Omega_1{}^1 &= d\omega_1{}^1 + \omega_k{}^1\wedge\omega_1{}^k = -\tilde\Psi_1 + \Phi_{10}\cr
\Omega_1{}^2 &= d\omega_1{}^2 + \omega_k{}^2\wedge\omega_1{}^k =\Phi_{00}\cr
\Omega_2{}^1 &= d\omega_2{}^1 + \omega_k{}^1\wedge\omega_2{}^k = -\Phi_{20}\cr
\Omega_2{}^2 &= d\omega_2{}^2 + \omega_k{}^2\wedge\omega_2{}^k = \tilde\Psi_1 - \Phi_{10}\cr}$$
\vskip 24pt

Now consider the quantity $S_a$. Its simple scaling behaviour (6.2.30) under change of LSR suggests $S_a$ is a more interesting quantity than $T_a$; indeed, $S_a$ is precisely the obstruction to $(M,g)$ being a Walker metric when $[\pi^{A'}]$ is real. What can one say when $S_a \not= 0$?

If $S^a \not= 0$, $[\pi^{A'}]$ defines a null distribution ${\cal D} := \langle S^a \rangle_{\bf R}$ of type I, contained within the $\alpha$-distribution $Z_{[\pi]}$. Since $S^a = \omega^A\pi^{A'}$, in accord with \S 6.1, select $o^A = \omega^A$ and $o^{A'} = \pi^{A'}$ in the choice of spin frames. If one restricts $A$ in (6.2.2) to the form
$$A = \pmatrix{\tilde\lambda&\mu\cr 0&\tilde\lambda^{-1}\cr},\eqno(6.2.35)$$
then, with $\lambda = \tilde\lambda$, (6.2.2) coincides with (6.1.3) and (6.2.31) with (6.1.6), and these restrictions guarantee the scaling behaviour of $\omega^A$ and $\pi^{A'}$ in (6.2.30). With these restrictions on the choice of spin frame, one has, noting (6.2.29),
$$\tilde\kappa = \tilde\sigma = \tilde\rho = 0 \hskip 1in \tilde\tau = 1.\eqno(6.2.36)$$
From (6.2.31), these conditions are invariant under the remaining freedom in the choice of spin frames.

By (6.1.9), $\cal D$ is auto-parallel iff $\tilde\kappa = \kappa = 0$. As $\tilde\kappa = 0$ already, it remains to investigate the vanishing of $\kappa$.
\vskip 24pt
\noindent {\bf 6.2.37 Lemma}\hfil\break
With the choice of spin frames specified in the preceding paragraphs, $\kappa  = \Phi_{00} = \Phi_{ABA'B'}\omega^A\omega^B\pi^{A'}\pi^{B'}$.

Proof. With the specified choice of spin frames $\kappa = \omega_BD\omega^B$. Observe that 
$$\eqalign{(\omega^AD\omega_A)\pi_{A'} &= \omega^AD(\omega_A\pi_{A'}) - \omega^A(D\pi_{A'})\omega_A\cr
&= \omega^A\omega^B\pi^{B'}\nabla_bS_a\cr
&= \omega^A\omega^B\pi^{B'}\nabla_b(\pi_{C'}\nabla_a\pi^{C'})\cr
&= \omega^A\omega^B\pi^{B'}(\nabla_b\pi_{C'})(\nabla_a\pi^{C'}) + \omega^A\omega^B\pi^{B'}\pi_{C'}\nabla_b\nabla_a\pi^{C'}\cr
&= \omega^A\omega^B(\eta_B\pi_{C'})(\nabla_a\pi^{C'}) + \omega^A\omega^B\pi^{B'}\pi_{C'}\nabla_b\nabla_a\pi^{C'}\cr
&= \omega^A(\omega^D\eta_D)(\omega_A\pi_{A'}) + \omega^A\omega^B\pi^{B'}\pi_{C'}\nabla_b\nabla_a\pi^{C'}\cr
&=  \omega^A\omega^B\pi^{B'}\pi_{C'}\nabla_b\nabla_a\pi^{C'}\cr
&= -\Phi_{ABC'D'}\omega^A\omega^B\pi^{C'}\pi^{D'}\pi_{A'}\qquad\hbox{using (6.2.13)(c)}\cr}$$
\vskip 24pt
One can choose Frobenius coordinates $(u,w,z,f)$ for the nested distributions ${\cal D} \leq Z_{[\pi]}$ so that $\alpha$-surfaces are given by setting $z$ and $f$ constant, and integral curves of $\cal D$ by further setting $w$ constant. Given such coordinates, choose an LSR $\pi^{A'}$ of $[\pi^{A'}]$ so that $\partial_u = S^a$ (changing the sign of $u$ if necessary). If we suppose that $\cal D$ is auto-parallel, i.e., $\Phi_{ABA'B'}\omega^A\omega^B\pi^{A'}\pi^{B'} = 0$ (which condition is obviously independent of the LSR of $[\pi^{A'}]$ chosen), then, as in \S 6.1, one can replace $u$ by an affine parameter for the integral curves of $\cal D$ so that $(v,w,z,f)$  are still Frobenius coordinates for the nested distributions. Now $\partial_v = (du/dv)\partial_u = (du/dv)S^a$, and $du/dv = \left[\exp(\int\Gamma^1_{11})\right]^{-1} > 0$, so by rescaling $\pi^{A'}$ and $\omega^A$ according to (6.2.30) (which is allowed by the freedom remaining in the choice of spin frames) by $\tilde\lambda = \sqrt{(du/dv)}$, $\partial_v = S^a$. But (6.2.30) forces $\lambda = \tilde\lambda$ in (6.1.6); to maintain $\partial_v = S^a$, where $v$ is any affine parameter, then $\tilde\lambda$ must be constant. Granted these choices, one has
$$\tilde\kappa = \tilde\sigma = \tilde\rho = \kappa = \epsilon + \tilde\epsilon = 0;\qquad \tilde\tau = 1,\eqno(6.2.38)$$
all of which are invariant under the remaining freedom in the choice of spin frames.

The above conditions cannot be maintained if one wishes to impose the A-P condition however. For example, $\tilde\epsilon \mapsto \tilde\lambda^2\tilde\epsilon$ under the current restrictions, so cannot be made zero if $\tilde\lambda$ is restricted to being constant. Imposing the A-P condition as done in \S 6.1 requires that one removes the constraint $\lambda = \tilde\lambda$ in (6.1.6) but maintains $\lambda = \tilde\lambda^{-1}$ to preserve $\partial_v = \ell^a$. Doing so violates (6.2.30), however. Assuming the situation reached in the previous paragraph, use a subscript 0 to indicate values of quantities under those conditions. Write $\ell^a = \partial_v = S^a_0 = \omega^A_0\pi^{A'}_0$. Let $\tilde\lambda$ be the necessary scaling of $\pi^{A'}_0$ and $\lambda = \tilde\lambda^{-1}$ of $\omega^A_0$ to impose the A-P condition. After this rescaling,
$$S^a_0 = \partial_v = \ell^a = o^A\pi^{A'} = (\tilde\lambda^{-1}\omega^A_0)(\tilde\lambda\pi^{A'}_0) = \tilde\lambda^{-2}\omega^A\pi^{A'} = \tilde\lambda^{-2}S^a$$
(where quantities without the subscript zero are the values after the rescaling). Using (6.1.6) with $\lambda = \tilde\lambda^{-1}$, one finds after the rescaling
$$\tilde\kappa = \tilde\sigma = \tilde\rho = \kappa = \epsilon = \tilde\epsilon = 0\qquad \tilde\tau = \tilde\lambda^2.\eqno(6.2.39)$$
For orthogonal connecting vector fields $V^a$ along integral curves of $\cal D$, (6.1.22) implies $\tilde\zeta$ is constant, i.e., as measured by $m^a$, orthogonal connecting vector fields of integral curves of $\cal D$ have a constant component pointing out of $\alpha$-surfaces; with respect to Frobenius coordinates, nearby integral curves of $\cal D$ in different $\alpha$-surfaces have constant separation in the direction of $m^a$. I shall return to these nested distributions in \S 6.3.

Now consider the situation when $S_a = 0$, i.e., $(M,g)$ is a four-dimensional Walker geometry, see [17]. We therefore return to the conditions (6.2.2) for the freedom in choice of the spin frames. By (6.2.4) and (6.2.29), with respect to such spin frames, Walker geometry is characterized by
$$\tilde\kappa = \tilde\sigma = \tilde\rho = \tilde\tau = 0.\eqno(6.2.40)$$
The Walker spin frames for given Walker coordinates are a specific choice; see \S 5. In (5.11), I considered the spin coefficient field equations (3.4) for the Walker spin frames. More generally, for the spin frames (6.2.2) when $S_a = 0$, one deduces from (3.4):
$$\vcenter{\openup2\jot \halign{$\hfil#$&&${}#\hfil$&\qquad$\hfil#$\cr
\hbox{(\~ a)}\ &\Rightarrow\ \Phi_{00} = 0 & \hbox{(\~ b)}\ &\Rightarrow\ \tilde\Psi_0 = 0;\cr
\hbox{(\~ e)}\ &\Rightarrow\ \Phi_{20} = 0 & \hbox{(\~ f)}\ &\Rightarrow\ \tilde\Psi_2 + 2\Lambda = 0,\cr}}\eqno(6.2.41)$$
\vskip -12pt
$$\hbox{(\~ c)} + \hbox{(\~ d)}\ \Rightarrow\ \tilde\Psi_1 = \Phi_{10} = 0.$$
which is equivalent to [17], 2.5. Thus, for a Walker geometry, $[\pi^{A'}]$ is a multiple WPS and a PS of $\Phi_{ABA'B'}$. In [17], \S 3, we saw that if the Walker geometry $(M,g,[\pi^{A'}])$ admits a parallel LSR, then $[\pi^{A'}]$ is WPS of multiplicity four and a multiple PS of $\Phi_{ABA'B'}$ (when the relevant curvature spinors are nonzero). Moreover, in this case we found a single function $\vartheta$ in terms of which the metric, whence the curvature, can be expressed; this formulation slightly generalizes Dunajski's [7] result and reduces to the real neutral-signature version of Pleba\'nski's [25] second heavenly form of the metric when there is a complementary parallel $\alpha$-distribution.

Since $\Phi_{ABA'B'}\pi^{A'}\pi^{B'} = 0$ already for a Walker geometry $(M,g,[\pi^{A'}])$, a natural condition to consider is
$$\Phi_{ABA'B'}\pi^{B'} = 0,\eqno(6.2.42)$$
which is obviously independent of the choice of LSR and weaker than the existence of a parallel LSR. I shall refer to such as a {\sl Ricci-null Walker geometry}. As in [17], write the metric components with respect to a choice of Walker coordinates $(u,v,x,y)$ as
$$g_{\rm ab} = \pmatrix{{\bf 0}_1&{\bf 1}_2\cr {\bf 1}_2&W\cr} \hskip 1in W = \pmatrix{a&c\cr c&b\cr}.\eqno(6.2.43)$$
In [17], (2.32), (6.2.42) is equivalent to $B_{AB} = 0$, i.e., to $\vartheta = \mu = \nu = 0$ [17] (2.33), which, by [17], A1.8, is equivalent to
$$a_{11} = b_{22} \hskip .5in b_{12} = -c_{11} \hskip .5in a_{12} = -c_{22},\eqno(6.2.44)$$
where the numerical subscripts 1, 2, 3, and 4 correspond to the coordinates $u$, $v$, $x$, and $y$ respectively. These conditions are explicitly implied by the stronger conditions [17], (3.2), see [17], (3.3--4), so it is natural to enquire whether one can still derive, as in [17], (3.7), a potential function $\vartheta$ for the metric.
\vskip 24pt
\noindent{\bf 6.2.45 Construction}\hfill\break
Let $(M,g,[\pi^{A'}])$ be a Ricci-null Walker geometry and $(u,v,x,y)$ a set of Walker coordinates. From (6.2.44), one can write
$$a_1 = -c_2 + f(u,x,y) \hskip 1.25in b_2 = - c_1 + g(v,x,y),\eqno(6.2.45{\rm a})$$
for certain functions $f$ and $g$ with the indicated dependence. But $a_{11} = b_{22}$ implies $f_1 = g_2$, whence $f_1$ and $g_2$ must be independent of $u$ and $v$:
$$f_1 = g_2 = h(x,y);\eqno(6.2.45{\rm b})$$
in particular, all second order partial derivatives of $f$ and $g$ with respect to $u$ and $v$ vanish.

Choosing specific antiderivatives $r := \int cdu$ and $s := \int cdv$, the equations $k_1 = s$ and $k_2 = r$ satisfy the classical integrability condition $k_{12} = k_{21}$; hence there is a solution $k$, with the freedom $k \mapsto k + z(x,y)$, $z$ an arbitrary function of $x$ and $y$. Define
$$\vartheta := {1 \over 2}k(u,v,x,y) + m(u,x,y) + n(v,x,y),\eqno(6.2.45{\rm c})$$
where, at this stage, $m$ and $n$ are arbitrary functions of their arguments. Now $2\vartheta_{12} = k_{12} = c$. Also $-2\vartheta_{22} = -k_{22} - 2n_{22} = -r_2 - 2n_{22}$. But $r_2$ is a specific antiderivative, with respect to $u$, of $c_2$. Using (6.2.45a), $r_2 = \int(f(u,x,y) - a_1)du = -a + F(u,x,y) + p(v,x,y)$, where $F(u,x,y)$ is an antiderivative, with respect to $u$, of $f(u,x,y)$ chosen to be independent of $v$, and $p$ is some specific function of $(v,x,y)$. Hence, $a = -r_2 + F + p = -2\vartheta_{22} + 2n_{22} + F + p$. Similarly, $b = -2\vartheta_{11} + 2m_{11} + G(v,x,y) + q(u,x,y)$, where $G$ is an antiderivative, with respect to $v$, of $g$ chosen to be independent of $u$, and $q$ is some specific function of $(u,x,y)$. Now choose $m$ and $n$ to eliminate $q$ and $p$, respectively, from this pair of equations. Hence, one now has
$$\displaylines{a = -2\vartheta_{22} + F(u,x,y),\qquad F_1 = f; \hskip 1in b = -2\vartheta_{11} + G(v,x,y);\qquad G_2 = g;\cr
\noalign{\vskip -6pt}
\hfill\llap(6.2.45{\rm d})\cr
\noalign{\vskip -6pt}
f_1 = g_2 = h(x,y) \hskip 1.25in c = 2\vartheta_{12}.\cr}$$
Conversely, given an arbitrary function $\vartheta(u,v,x,y)$, functions $f(u,x,y)$ and $g(v,x,y)$ satisfying $f_1 = g_2 = h(x,y)$, for some $h$, defining $F(u,x,y)$ and $G(v,x,y)$ so that $F_1 = f$ and $G_2 = g$, and then defining $a$, $b$, and $c$ as in (6.2.45d) yields a Walker metric (6.2.43) with 
$$W = -2\pmatrix{\vartheta_{22}&-\vartheta_{12}\cr -\vartheta_{12}&\vartheta_{11}\cr} + \pmatrix{F(u,x,y)&0\cr 0&G(v,x,y)\cr}.\eqno(6.2.45{\rm e})$$
The conditions (6.2.44) are automatically fulfilled so the metric is Ricci-null Walker. This description therefore provides a (local) characterization of Ricci-null Walker metrics which generalizes [17], 3.7.
\vskip 24pt
I now determine the curvature of Ricci-null Walker metrics in terms of the functions $\vartheta$, $f$, and $g$; in particular, I seek a generalization of the expression in [17], (3.12), which proves so useful for specifying $\tilde\Psi_{A'B'C'D'}$ and $\Phi_{ABA'B'}$. First, note that
$$S = a_{11} + b_{22} + 2c_{12} = (a_1 + c_2)_1 + (b_2 + c_1)_2 = f_1 + g_2 = 2f_1 = 2g_2 = 2h.\eqno(6.2.46)$$
The  nonzero components of $\tilde\Psi_{A'B'C'D'}$ with respect to the Walker spin frames are, from [17], (2.20), $\tilde\Psi_2 = S/12$, $\tilde\Psi_3 = -(B+Sc)/8$ and $\tilde\Psi_4 = (6Bc-A+S(3c^2-1)/24$, where $B$ and $A$ are as in [17], (A1.23--24). From these expressions, one computes directly
$$B+Sc = 2(a_{14} + c_{24} - b_{23} - c_{13}) = 2(f_4 - g_3),\eqno(6.2.47)$$
and
$$A - 6Bc - S(3c^2-1) = A - 12c(f_4 - g_3) + 2(3c^2+1)f_1.\eqno(6.2.48)$$
By [17], (3.9), in a Walker geometry one has, for an arbitrary function $H$,
$$\Square H = -aH_{11} - 2cH_{12} - bH_{22} + 2H_{13} + 2H_{24} - (a_1+c_2)H_1 - (b_2+c_1)H_2,$$
which reduces to
$$\Square H = -aH_{11} - 2cH_{12} - bH_{22} + 2H_{13} + 2H_{24} - fH_1 - gH_2,\eqno(6.2.49)$$
in the Ricci-null Walker case. As in [17], (3.12), define
$$P := \vartheta_{13} + \vartheta_{24} + \vartheta_{11}\vartheta_{22} - (\vartheta_{12})^2.\eqno(6.2.50)$$
Additionally, define
$$\displaylines{Q := {g\vartheta_2 - G\vartheta_{22} + f\vartheta_1 - F\vartheta_{11} - h\vartheta \over 2}\cr
\noalign{\vskip -6pt}
\hfill\llap(6.2.51)\cr
\noalign{\vskip -6pt}
T := -{(vF_4 + uG_3) \over 4} \hskip 1in R := P + Q + T.\cr}$$
A somewhat tedious calculation yields
$$\eqalignno{A - 6Bc - S(3c^2 - 1) &= 12\Square R + 24fR_1 + 24gR_2\cr
&- 3(fG_3 + gF_4) - 6(F_{44} + G_{33}) + 3(vff_4 + ugg_3) + 6(vf_{34} + ug_{34})\cr
&- 3vh_4(F - 2\vartheta_{22}) - 3uh_3(G- 2\vartheta_{11}).&(6.2.52)\cr}$$
Note the manifest invariance under the Walker symmetry [17], (A1.7), when in addition to $u\ \leftrightarrow\ v$, $x \leftrightarrow\ y$, and $a \leftrightarrow b$ one requires, for consistency, $f\ \leftrightarrow\ g$, whence $F\ \leftrightarrow\ G$, and $h \leftrightarrow h$.

Also note that the condition for $(M,g,[\pi^{A'}])$ to admit a parallel LSR is, by [17], (3.2), $a_1 + c_2 = 0 = b_2 + c_1$, i.e., $f = 0 = g$, whence $h=0$. Then $F$ and $G$ are functions of $x$ and $y$ alone. To recover the construction in [17], 3.7, one can absorb $F$ into $p$ and $G$ into $q$; in effect, one takes $F$ and $G$ zero. Then it is clear that not only is $B=0$, but $Q = T = 0$, $R=P$, and (6.2.52) reduces to [17], (3.11), and (6.2.45e) to [17], (3.8).

Considering now the ASD Weyl curvature, the components of $\Psi_{ABCD}$ with respect to the Walker spin frames are, from [17], (2.25), and using (6.2.44) and (6.2.45d):
$$\displaylines{\Psi_0 = {b_{11} \over 2} = -\vartheta_{1111}\cr
\Psi_1 = {b_{12} \over 2}= -\vartheta_{1112}\cr
\hfill\Psi_2 = {b_{22} - 2c_{12} \over 6} = {-2\vartheta_{1122} + G_{22} - 4\vartheta_{1212} \over 6} = -\vartheta_{1122} + {h \over 6}\hfill\llap(6.2.53)\cr
\Psi_3 = {a_{12} \over 2} = -\vartheta_{1222}\cr
\Psi_4 = {a_{22} \over 2} = -\vartheta_{2222},\cr}$$
equations which differ from [17], (3.13), only by the additional term in $h$ in $\Psi_2$. 

It is convenient to introduce, for oriented Walker coordinates and their associated Walker spin frames, the spinor operator
$$\eqalignno{\delta_A &:= \pi^{A'}\nabla_{AA'} = \alpha_A\tilde m^b\nabla_b - \beta_A\ell^b\nabla_b = \alpha_A\triangle - \beta_AD\cr
&= \alpha_A\partial_2 - \beta_A\partial_1,\qquad\hbox{when acting on functions}.&(6.2.54{\rm a})}$$
This operator represents, in effect, the induced linear connection within $\alpha$-surfaces. Since, by (5.8), the Walker spin frames are parallel on $\alpha$-surfaces, writing the Walker spin frames as $\epsilon_{\bf A}{}^A$ and $\epsilon_{\bf A'}{}^A$, then for any spinor $\chi^{S\ldots X'\ldots}_{L\ldots P'\ldots}$,
$$\delta_A\chi^{S\ldots X'\ldots}_{L\ldots P'\ldots} = \left(\delta_A\chi^{{\bf S}\ldots{\bf X}'\ldots}_{{\bf L}\ldots{\bf P}'\ldots}\right)\epsilon_{\bf S}{}^S\ldots\epsilon_{{\bf X}'}{}^{X'}\ldots\epsilon_L{}^{\bf L}\ldots\epsilon_{P'}{}^{{\bf P}'}\ldots\eqno(6.2.54.{\rm b})$$
which is consistent with the earlier observation that, in Walker geometry, the induced linear connection within $\alpha$-surfaces is flat. Furthermore, the components of $\delta_A\ldots\delta_K\chi^{S\ldots X'\ldots}_{L\ldots P'\ldots}$ with respect to the Walker spin frames are just
$$\delta_{\bf A}\ldots\delta_{\bf K}\chi^{\bf S\ldots X'\ldots}_{\bf L\ldots P'\ldots} = (-1)^q{\partial^k\chi^{\bf S\ldots X'\ldots}_{\bf L\ldots P'\ldots} \over \partial u^q\partial v^{k-q}},\eqno(6.2.54{\rm c})$$
where there are $k$ applications of $\delta_A$ and $q$ 1's amongst the ${\bf A},\ldots,{\bf K}$. It follows that $\delta_A$ and $\delta_B$ commute, whence $\delta_A\delta_B\ldots\delta_K\chi^{S\ldots X'\ldots}_{L\ldots P'\ldots}$ is totally symmetric in $A,\ldots,K$ (recall from \S 5 that $D$ and $\triangle$ commute on arbitrary spinors). In particular, for $k$ applications of $\delta_A$ to a function $F$:
$$\eqalignno{\delta_A\ldots\delta_KF &= \delta_A\ldots\delta_J(\alpha_KF_v-\beta_KF_u)\cr
&= \delta_A\ldots\delta_I(\alpha_K\delta_JF_v-\beta_K\delta_JF_u)\cr
&= \sum_{q=0}^k\,(-1)^q{k \choose q}\underbrace{\alpha_{(A}\ldots\alpha_F}_{k-q}\underbrace{\beta_G\ldots\beta_{K)}}_{q}{\partial^kF \over \partial u^q\partial v^{k-q}},&(6.2.54{\rm d})\cr}$$
whence $\delta_A\ldots\delta_KF = 0$ means all partial derivatives, with respect to $u$ and $v$, of $F$ of order $k$ vanish.

Similarly to [17], 3.7, one finds
$$\Psi_{ABCD} = -\delta_A\delta_B\delta_C\delta_D \left(\vartheta - {u^2v^2h \over 24}\right),\eqno(6.2.55)$$
which obviously generalizes [17], (3.15). Note that $\Psi_{ABCD}$ vanishes iff $\vartheta - (u^2v^2h/24)$ is polynomial in $u$ and $v$ with each term of total degree in $u$ and $v$ at most three.

Finally, the Ricci-null condition amounts to $\Phi_{ABA'B'} = A_{AB}\pi_{A'}\pi_{B'}$, with $A_{AB}$ as in [17], (2.32--33), but with the conditions (6.2.44) in force. Actually, it is simpler to first invoke the equivalent conditions $\mu = \nu = \theta = 0$ (i.e., $B_{AB} = 0$ in [17], (2.32--33)). Hence,
$$A_{\rm AB} = {1 \over 2}\pmatrix{\Upsilon&-\eta\cr -\eta&\zeta\cr}.\eqno(6.2.56)$$
Calculating the three quantities $\Upsilon$, $\eta$, and $\zeta$ from [17], A1.7, under (6.2.44), one finds:
$$\displaylines{A_{00} = {\Upsilon \over 2} = (P + Q)_{11}\cr
\hfill A_{01} = -{\eta \over 2} = (P + Q + T)_{12}\hfill\llap(6.2.57)\cr
A_{11} = {\zeta \over 2} = (P + Q)_{22}.\cr}$$
While (6.2.57) does generalize [17], (3.17), it is not quite in the form of that equation.

The algebraic classification of $\tilde\Psi_{A'B'C'D'}$ for any Walker metric was given in [17], 2.6. Supposing first that $S \not= 0$, then $\tilde\Psi_{A'B'C'D'}$ is of type $\{2,2\}$Ia if $S^2 + AS + 3B^2 = 0$ and of type $\{211\}$II or $\{1\overline 12\}$II if $S^2 + AS + 3B^2 \not= 0$. In each case, $[\pi^{A'}]$ is, of course, the multiple WPS.

Now suppose $S=0$. In addition to (6.2.44), one then also has $a_{11} = b_{22} = -c_{12}$. From (6.2.46), $h = f_1 = g_2 = 0$, so $f$ and $g$ are functions of $x$ and $y$ alone and one can take 
$$F(u,x,y) = uf(x,y) \hskip .75in G(v,x,y) = vg(x,y) \hskip .75in h(x,y) = 0.\eqno(6.2.58)$$
The metric takes the Walker canonical form with
$$W = -2\pmatrix{\vartheta_{22}&-\vartheta_{12}\cr -\vartheta_{12}&\vartheta_{11}\cr} + \pmatrix{uf(x,y)&0\cr 0&vg(x,y)\cr}.\eqno(6.2.59)$$
The nonzero curvature components  of $\tilde\Psi_{A'B'C'D'}$ are:
$$\displaylines{\tilde\Psi_3 = -{B \over 8} = {g_3 - f_4 \over 4};\cr
\noalign{\vskip -6pt}
\hfill\llap(6.2.60)\cr
\noalign{\vskip -6pt}
\tilde\Psi_4 = {6Bc - A \over 24} = -{\Square R \over 2} - fR_1 - gR_2 + {(ug-vf)(f_4-g_3) \over 8} + {u(f_4-g_3)_4 - v(f_4-g_3)_3 \over 4}.\cr}$$
By [17] 2.6, if $B = 2(f_4-g_3) = 0$, but $A \not= 0$, $\tilde\Psi_{A'B'C'D'}$ is of type $\{4\}$II; if $B \not = 0$, $\tilde\Psi_{A'B'C'D'}$ is of type $\{31\}$III; in each case $[\pi^{A'}]$ is the multiple WPS.
The remaining curvature is
$$\displaylines{\Psi_{ABCD} = -\delta_A\delta_B\delta_C\delta_D\vartheta\cr
\noalign{\vskip -6pt}
\hfill\llap(6.2.61)\cr
\noalign{\vskip -6pt}
\Phi_{ABA'B'} = A_{AB}\pi_{A'}\pi_{B'} \hskip .75in A_{AB} = \delta_A\delta_BR,\cr}$$
noting that under (6.2.58)
$$T = -{uv(f_4+g_3) \over 4} \hskip 1in Q = {g(\vartheta_2 - v\vartheta_{22}) + f(\vartheta_1 - u\vartheta_{11}) \over 2},\eqno(6.2.62)$$
the form of $T$ explaining why (6.2.57) now reduces to the simpler form given in (6.2.61), which permits a simple characterization of the Einstein condition $\Phi_{ABA'B'} = 0$.
\vskip 24pt
\noindent {\bf 6.2.63 Proposition}\hfil\break
A Ricci-null, Ricci-scalar-flat Walker geometry $(M,g,[\pi^{A'}])$ is Einstein iff the function $R$ is affine in $u$ and $v$, i.e.,
$$\vartheta_{13} + \vartheta_{24} + \vartheta_{11}\vartheta_{22} - (\vartheta_{12})^2 + {g(\vartheta_2 - v\vartheta_{22}) + f(\vartheta_1 - u\vartheta_{11}) \over 2} - {uv(f_4+g_3) \over 4} = uA(x,y) +  vB(x,y) + C(x,y).$$
\vskip 24pt
As it happens, the last equation is a real, neutral-signature version of the {\sl hyperheavenly equation for nonexpanding $\cal HH$-spaces} in complex general relativity, see [27], [28], [8], and [2]. The results presented here arise naturally in the context of null geometry of (four-dimensional) neutral geometry, specifically of $\alpha$-distributions (the focus on $\alpha$-distributions is also paramount in the work cited on complex GR), and takes advantage of Walker's [32] results on parallel distributions and the spinor formulation presented in [17]. The equation presented in (6.2.63) differs from the hyperheavenly equation presented in [2], (5.10), due to the different derivation presented here. $\cal HH$-spaces, i.e., hyperheavens, are four-dimensional holomorphic Einstein spaces with algebraically special SD Weyl curvature. Pleba\'nski and co-workers, in the citations listed above, also considered what they called {\sl expanding} $\cal HH$ spaces, which in the neutral context amounts to an Einstein $(M,g)$ with an integrable $\alpha$-distribution (recall that in the Einstein case, the GGST says an integrable $\alpha$-distribution $Z_{[\pi]}$ is equivalent to $[\pi^{A'}]$ being a multiple WPS) for which $S_a \not=0$, i.e., not Walker. The case when $(M,g)$ admits an integrable $\alpha$-distribution $Z_{[\pi]}$ for which $[\pi^{A'}]$ is a multiple WPS but $(M,g,[\pi^{A'}])$ is not Walker is treated in [18], where it is shown that such geometries are locally conformal to a Walker geometry; in particular, an independent derivation of the (real neutral-signature version of the) hyperheavenly equation for expanding ${\cal HH}$-spaces is obtained there.

Although there was a smattering of work in the second half of the twentieth century exploiting Walker's [32] canonical form for the metric of a pseudo-Riemannian manifold with a parallel totally null distribution, see for example the bibliography in [6], the utility of Walker's canonical form seems only to have been more widely appreciated in this century, beginning perhaps with [9]. The presentation given here underscores the fact that Walker's [32] results lay largely unexploited for fifty years, apparently playing no role in the work done in complex GR in the 70s, though the related work on recurrent tensors, see [31], was utilized in general relativity itself, see for example [11].
\vskip 24pt
\noindent {\bf 6.3 Null Distributions of Type III}
\vskip 12pt
Let $\cal H$ be a null distribution of type III and ${\cal D} := {\cal H}^\perp$. With $\ell^a$ a local section of $\cal D$ then, as in \S 6.1 one can, at least locally, write $\ell^a = o^Ao^{A'}$ and then form spin frames $\{o^A,\iota^A\}$ and $\{o^{A'},\iota^{A'}\}$ and the associated null tetrad $\{\ell^a,n^a,m^a,\tilde m^a\}$. Then, $\{\ell^a, \tilde m^a,m^a\}$ is a local frame field for $\cal H$. The freedom in these choices is again (6.1.3--5).

The first question to address is the integrability of $\cal H$. Let $\cal I(\cal H)$ denote the ideal of forms, within the exterior algebra $\Lambda_\bullet(M)$, which vanish on $\cal H$. Locally, $\cal I(\cal H)$ is generated by $\ell_a$, so one has the following standard result.
\vskip 24pt
\noindent {\bf 6.3.1 Lemma}\hfil\break
$$\eqalign{{\cal H}\hbox{ is integrable}\ &\Leftrightarrow\ d\ell_a \in \cal I(\cal H)\cr
&\Leftrightarrow\ \ell_a \propto dF,\qquad\hbox{some function } F\cr
&\Leftrightarrow\ d\ell_a \wedge \ell_b = 0.\cr}$$
Consequently, the integral manifolds of an integrable $\cal H$ are, locally, the level surfaces of a function $F$ satisfying $\ell_a \propto dF$.

Proof. See, for example, [30], pp. 293--295
\vskip 24pt
\noindent {\bf 6.3.2 Lemma}\hfil\break
A null distribution $\cal H$ of type III is integrable iff ${\cal D} = {\cal H}^\perp$ is auto-parallel and, with respect to the spin frames (6.1.3), $\rho = \tilde\rho$.

Proof. The condition $d\ell \wedge \ell = 0$ is equivalent to $(\nabla_{[b}\ell_c)\ell_{a]} = 0$. The argument of [23], p. 180, carries over to the case of neutral signature: with $*$ the Hodge (${\bf R}^{2,2}$)-star operator on {\sl two\/}-forms acting on the indices $b$ and $c$, $*(\nabla_{[b}\ell_c)\ell_{a]} = 0$, which in turn is equivalent, by [22] (3.4.26), to $(*\nabla_{[a}\ell_{b]})\ell^a = 0$, i.e., to $\ell^{AA'}(\nabla_{AB'}\ell_{BA'} - \nabla_{BA'}\ell_{AB'}) = 0$, i.e., to $o^{A'}o_Bo^A\nabla_{AB'}o_{A'} = o^Ao_{B'}o^{A'}\nabla_{BA'}o_A$. Taking components with respect to the relevant spin frames yields, from (2.9): $\tilde\kappa = \kappa = 0$ and $\rho = \tilde\rho$. The assertion then follows from (6.1.9).
\vskip 24pt
Recall from (6.1.40) that, when the A-P condition is imposed on an auto-parallel null distribution of type I, the condition $\rho = \tilde\rho$ is not a restriction on curvature. 

I will henceforward assume $\cal H$ is integrable, in which case one has nested integrable distributions: ${\cal D} \leq {\cal H}$. One can therefore choose Frobenius coordinates $(u,x,y,f)$ such that: level hypersurfaces of $f$ are integral manifolds of $\cal H$; integral curves of $\cal D$ are given by setting $x$, $y$, $f$ constant and $\partial_u$ is a local section of $\cal D$.

Since $\cal D$ is auto-parallel, one can repeat the procedure in \S 6.1, first replacing $u$ by an affine parameter $v$ for the integral curves of $\cal D$ but so that $(v,x,y,f)$ are still Frobenius coordinates for the nested distributions ${\cal D} \leq {\cal H}$; then, with $\partial_v = :\ell^a = o^Ao^{A'}$, use the freedom in (6.1.3), without changing $\ell^a$, to rescale each of $o^A$ and $o^{A'}$ to be $D$ parallel; and then finally choose $\iota^A$ and $\iota^{A'}$ to be $D$ parallel; i.e., impose the A-P condition on $\cal D$.

Integral curves of $\cal D$ lying within an integral hypersurface of $\cal H$ are called {\sl null geodesic generators\/} of that hypersurface. Connecting vector fields between null geodesic generators within a given integral hypersurface of $\cal H$ must lie in $\cal H$, i.e., are orthogonal connecting vector fields of $\cal D$ and thus described by (6.1.22) and (6.1.46). Note from (6.1.48) that $\Sigma$ vanishes identically on these orthogonal connecting vector fields. To understand non-orthogonal connecting vector fields of $\cal D$, consider how Frobenius coordinates $(u,x,y,f)$ relate to the geometry.
\vskip 24pt
\noindent {\bf 6.3.3 Lemma}\hfil\break
$N^a := g^{ab}(df)_b \in {\cal D}$, whence null; moreover, integral curves of $\cal D$ parametrized so that $N^a$ is their tangent vector are null geodesics.

Proof. Any level hypersurface of the function $f$ has tangent space equal to $\ker(df)$, i.e., to $\langle N \rangle_{\bf R}^\perp$, i.e., ${\cal H} = \langle N \rangle_{\bf R}^\perp$; equivalently, ${\cal D} = \langle N \rangle_{\bf R}$. Moreover,
$$N^a\nabla_aN_b = N^a\nabla_a\nabla_bf = N^a\nabla_b\nabla_af = N^a\nabla_bN_a = {1 \over 2}\nabla_b(N^aN_a) = 0,$$
as $N^a$ is null.
\vskip 24pt
\noindent {\bf 6.3.4 Observation}\hfil\break
Begin with any Frobenius coordinates $(u,x,y,f)$ for the nested distributions ${\cal D} \leq {\cal H}$. By (6.3.3), $\partial_u = hN^a$, for some function $h$, whence $g(\partial_u,\partial_f) = g(hN,\partial_f) = hdf(\partial_f) = h$.  As ${\cal H} = \langle \partial_u,\partial_x,\partial_y \rangle_{\bf R}$ and ${\cal D} = {\cal H}^\perp$, then, with respect to the Frobenius coordinates
$$(g_{\bf ab}) = \pmatrix{0&\matrix{0&0\cr}&h\cr \matrix{0\cr 0\cr}&{\cal A}&\cr h&&\cr},\eqno(6.3.4{\rm a})$$
for some $({\cal A}) \in {\bf R}_+(3)$. When $u$ is replaced by an affine parameter $v$, $h \equiv \hbox{constant}$, and one can choose $v$ so that $h = 1$, i.e., $\ell^a = \partial_v = N^a$, whence $g(\ell^a,\partial_f) = g(N^a,\partial_f) = 1$. Noting that ${\cal H} = {\cal D}^\perp = \langle \ell^a, m^a, \tilde m^a \rangle_{\bf R}$, one now has
$$\vcenter{\openup1\jot \halign{$\hfil#$&&${}#\hfil$&\qquad$\hfil#$\cr
\partial_v &= \ell^a = N^a & \partial_f &=: L\ell^a + Mm^a + \tilde M \tilde m^a + n^a\cr
\partial_x &=: A\ell^a + Bm^a + \tilde B \tilde m^a & \partial_y &=: D\ell^a + Cm^a + \tilde C \tilde m^a,\cr}}\eqno(6.3.4{\rm b})$$
for some functions $L$, $M$, $\tilde M$, $A$, $B$, $\tilde B$, $D$, $C$, and $\tilde C$, whence
$${\cal A} = \pmatrix{-2B\tilde B&-(B\tilde C + \tilde BC)&A - (B\tilde M + \tilde BM)\cr -(B\tilde C + \tilde BC)&-2C\tilde C&D - (C\tilde M + \tilde C M)\cr A - (B\tilde M + \tilde BM)&D - (C\tilde M + \tilde C M)&2(L - M\tilde M)\cr}.\eqno(6.3.4{\rm c})$$
\vskip 24pt
With respect to the spin frames constructed in association with $\cal D$ (6.1.3), in addition to $\tilde\kappa = \kappa = 0$ and $\rho=\tilde\rho$ (integrability of $\cal H$) and $\epsilon + \tilde\epsilon = 0$ (affine parametrization of the integral curves of $\cal D$), one obtains the following result from $\ell^a = N^a$.
\vskip 24pt
\noindent {\bf 6.3.5 Lemma}\hfill\break
With the choices in (6.3.4), specifically an affine parameter $v$ so that $\ell^a = \partial_v = N^a$,
$$\tau + \tilde\alpha + \beta = 0 = \tilde\tau + \alpha + \tilde\beta.$$

Proof. Since $\ell_a = N_a = df$, $d\ell_a = 0$ and the result follows immediately from the first equation in (4.1).
\vskip 24pt
Supposing now that the A-P condition has been imposed on the spin frames (6.1.3) associated with $\cal D$, then altogether one has
$$\tilde\kappa = \kappa = \tilde\epsilon = \epsilon = \tau' = \tilde\tau' = 0 = \tau + \tilde\alpha + \beta = \tilde\tau + \alpha + \tilde\beta\qquad \rho = \tilde\rho.\eqno(6.3.6)$$
The equations (6.1.38) simplify in that the first {\sl pair} are identical to the last  {\sl pair} and the equations of the second {\sl pair} are identical, so the independent equations are:
$$\vcenter{\openup1\jot \halign{$\hfil#$&&${}#\hfil$&\qquad$\hfil#$\cr
D\sigma &= -2\rho\sigma - \Psi_0 & D\tilde\sigma &= -2\tilde\rho\tilde\sigma - \tilde\Psi_0\cr
D\tau &= -(\tau\rho + \tilde\tau\sigma) - \Psi_1 - \Phi_{01} & D\tilde\tau &= -(\tilde\tau\tilde\rho + \tau\tilde\sigma) - \tilde\Psi_1 - \Phi_{10}\cr
D\tilde\rho &= -(\tilde\rho^2 + \sigma\tilde\sigma) - \Phi_{00} &
D(\gamma + \tilde\gamma) &= -2\tau\tilde\tau + \Psi_2 + \tilde\Psi_2 + 2\Phi_{11} - 2\Lambda.\cr}}\eqno(6.3.7)$$
Since $\cal D$ is auto-parallel, all connecting vector fields (6.1.16) have $\nu$ constant, see (6.1.19). Those with $\nu = 0$ lie in $\cal H$ and `connect' null geodesic generators of a given integral hypersurface. By (6.3.4b), a nonorthogonal connecting vector field $V^a$ with $\nu = \nu_0$ has the form
$$V^a = \eta\ell^a + \zeta\tilde m^a + \tilde\zeta m^a + \nu_0 n^a = H^a + \nu_0\partial_f$$
for some $H^a \in {\cal H}$, i.e., $V^a$ `connects' a null geodesic generator in a level hypersurface $f = d$ with null geodesic generators in $f = d+ \nu_0$. Note that (6.1.19) reduces to
$$D\pmatrix{\eta\cr \zeta\cr \tilde\zeta\cr \nu_0\cr} = \pmatrix{0&-\tilde\tau&-\tau&\gamma+\tilde\gamma\cr 0&\rho&\sigma&\tau\cr 0&\tilde\sigma&\rho&\tilde\tau\cr 0&0&0&0\cr}.\pmatrix{\eta\cr \zeta\cr \tilde\zeta\cr \nu_0\cr}.\eqno(6.3.8)$$
One also sees from (6.1.42) that $\Sigma$ in fact vanishes on all connecting vector fields of $\cal D$.
\vskip 24pt
\noindent {\bf 6.3.9 Lemma}\hfil\break
A null distribution $\cal H$ of type III is auto-parallel iff
$$\kappa = \tilde\kappa = \sigma = \tilde\sigma = \rho = \tilde\rho = 0,\eqno(6.3.9{\rm a})$$
and parallel iff, in addition to (6.3.9a), $\tau = \tilde\tau = 0$; of course, $\cal H$ is parallel iff ${\cal D} = {\cal H}^\perp$ is parallel, whence one recovers (6.1.12a).

If $(M,g)$ admits an auto-parallel null distribution of type III, then
$$\Psi_0 = \tilde\Psi_0 = \Phi_{00} = 0 \qquad \Psi_1 = \Phi_{01} \qquad \tilde\Psi_1 = \Phi_{10};\eqno(6.3.9{\rm b})$$
in particular, $o^{A'}$ and $o^A$ are each WPSs. If $(M,g)$ admits a parallel null distribution of type III, equivalently a parallel null distribution of type I, then
$$\Psi_0 = \Psi_1 = 0 = \tilde\Psi_0 = \tilde\Psi_1 \qquad \Phi_{00} = \Phi_{10} = \Phi_{01} = \Phi_{20} = \Phi_{02} = 0 \qquad \Psi_2 = -2\Lambda = \tilde\Psi_2;\eqno(6.3.9{\rm c})$$
in particular, $o^A$ and $o^{A'}$ are each WPSs of multiplicity at least two and both are PSs of $\Phi_{ABA'B'}$.

Proof. By (6.1.7), $\cal H$ is auto-parallel iff $\nabla_YX \in {\cal H}$ for all local sections $Y$ and $X$ of $\cal H$, i.e., iff $\nabla_YX \perp {\cal D}$, i.e., $\ell_a\nabla_YX^a = 0$. Expressing $X$ and $Y$ in terms of $\ell^a$, $m^a$ and $\tilde m^a$ yields (6.3.9a). Allowing $X$ to be an arbitrary vector field yields the additional two conditions. (6.3.9b--c) follow from (3.4); also see (6.1.38).
\vskip 24pt
\noindent {\bf 6.3.10 Remarks}\hfil\break
Walker [33] gave a canonical coordinate form for the metric $h$ of a pseudo-Riemannian manifold $(N,h)$ admitting a parallel distribution of degenerate subspaces. For null distributions of type III in neutral four-manifolds $(M,g)$, Walker's canonical form is (6.3.4) with $h=1$ and the condition of being parallel additionally allowing one to choose the Frobenius coordinates $(v,x,y,f)$ so that ${\cal A}_{11}$, ${\cal A}_{12}$, ${\cal A}_{13}$, ${\cal A}_{22}$ and ${\cal A}_{23}$ are independent of $v$. This form coincides with the Walker [32] canonical form (6.1.12b) for the parallel distribution ${\cal D} = {\cal H}^\perp$ (and, as noted there, this form does entail that $v$ is indeed affine).

When $\cal H$ is auto-parallel, orthogonal connecting vector fields along null geodesic generators of integral hypersurfaces are represented, according to (6.1.20--22), by $q^a$   with constant components with respect to $m^a$ and $\tilde m^a$. When $\cal H$ is parallel, all connecting vector fields along null geodesic generators have a particularly simple form:
$$V^a = \eta(v)\ell^a + \zeta_0\tilde m^a + \tilde\zeta_0 m^a + \nu_0n^a,$$
with $\zeta_0$, $\tilde\zeta_0$, and $\nu_0$ constant and $\eta(v)$ satisfying $D\eta = (\gamma + \tilde\gamma)\nu_0$.
\vskip 24pt
Clearly even the auto-parallel condition on $\cal H$ is quite strong. As $\kappa = \tilde\kappa = 0$ for integrability of $\cal H$, the vanishing of each of $\sigma$, $\tilde\sigma$, and $\rho=\tilde\rho$ is well defined, see (6.1.6). Of particular interest is the vanishing of $\sigma$ and/or $\tilde\sigma$. The distribution $\cal H$ contains two null distributions of type II, a distribution of $\beta$-planes $W_{[o]}$ and a distribution of $\alpha$-planes $Z_{[\tilde o]}$ with local frame fields $\{\ell^a,m^a\}$ and $\{\ell^a,\tilde m^a\}$, respectively. By (6.2.4), integrability of $W_{[o]}$ is equivalent to $\kappa = \sigma = 0$, while integrability of $Z_{[\tilde o]}$ is equivalent to $\tilde\kappa = \tilde\sigma = 0$. By (3.4)(b) and (\~ b), or (6.1.39), these conditions imply $\Psi_0 = 0$ and $\tilde\Psi_0 = 0$ respectively, i.e., $o^A$ is a WPS and $o^{A'}$ a WPS, respectively.

Of special interest is when $\sigma = \tilde\sigma = 0$, i.e., both $W_{[o]}$ and $Z_{[\tilde o]}$ are integrable and the integral hypersurfaces of $\cal H$ are foliated by $\beta$-surfaces and by $\alpha$-surfaces. In this case, (6.1.22) reduces to $D\zeta = \rho\zeta$ and $D\tilde\zeta = \rho\tilde\zeta$, with solution given by the first equation of (6.1.30), i.e., orthogonal connecting vector fields along a null geodesic generator suffer dilation only. Moreover, (6.1.55) reduces to $T_{ab} = \rho h_{ab}$, i.e., for $X$, $Y$ vector fields in $\cal H$,
$$X^aY^b\nabla_a\ell_b = \rho h([X],[Y]) = \rho X_aY^b,\eqno(6.3.11)$$
with notation as in \S 6.1.

Alternatively, in (6.1.25--27), one sees that integrability of $\cal H$ and affine parametrization of the null geodesic generators yields $\nabla_{[a}\ell_{b]} = 0$, whence $\nabla_a\ell_b = \nabla_{(a}\ell_{b)}$. Since $g_{ab} = 2(\ell_{(a}n_{b)} - m_{(a}\tilde m_{b)})$, the induced metric on ${\cal Q} = {\cal D}^\perp/{\cal D} = {\cal H}/{\cal D}$ can be represented by $-2m_{(a}\tilde m_{b)}$. When $\sigma = \tilde\sigma = 0$, the expression for $\nabla_{(a}\ell_{b)}$ given in \S 6.1 has as its only nonzero term $-2\rho m_{(a}\tilde m_{b)}$ consistent with (6.3.11).

In fact, the induced metric on an integral hypersurface $\cal S$ is effectively $-2m_{(a}\tilde m_{b)}$. Moreover, as $\ell^a \in {\cal H}$,
$$\eqalign{\Lie_\ell\left(g\vert_{\cal S}\right) &= \left(\Lie_\ell g\right)\vert_{\cal S}\cr
&= (\nabla_a\ell_b)\vert_{\cal S}\cr
&= \rho\left(g\vert{\cal S}\right)\cr}$$
as $\sigma = \tilde\sigma = 0$, i.e., $\ell^a$ acts as a `conformal Killing vector' of (the degenerate) $g\vert_{\cal S}$.

Let $B$ denote the set of null geodesic generators of $\cal S$. If $(v,x,y,f)$ are Frobenius coordinates for $\cal H$ such that $f=k$ is the intersection of the coordinate domain with $\cal S$, then $(x,y)$ provide local coordinates on a neighbourhood of the point $b \in B$ labelling the null geodesic generator $(v,0,0,k)$. Assume $B$ can be covered by coordinates obtained from such Frobenius coordinates for $\cal H$ to provide $B$ with a smooth structure. The projection $\pi:{\cal S} \to B$ is a submersion (locally, $\pi:(v,x,y) \mapsto (x,y)$). Hence, for $p \in \pi^\dashv(\{b\})$, $\cal X$, ${\cal Y} \in T_bB$, one can choose $X$, $Y \in T_p{\cal S}$ such that $D\pi_p(X) = {\cal X}$ and $D\pi_p(Y) = {\cal Y}$. Of course, $D\pi_p:{\cal Q}_p \to T_pB$ is an isomorphism. More generally, then, any tangent vector $p\partial_x + q\partial_y \in T_bB$ lifts to vector fields of the form $V := p\partial_x + q\partial_y + \eta\partial_v$ along the null geodesic generator $C_b$ labelled by $b$. As $\ell^a = \partial_v$, $[\ell,V]^a = \partial_v(\eta)\ell^a$. Taking $\eta$ constant along $C_b$ yields an orthogonal connecting vector field $V^a$ along $C_b$. Thus, if one lifts $\cal X$ and $\cal Y$ to orthogonal connecting vector fields $V^A_1$ and $V^a_2$ along $C_b$, respectively, then, in terms of (6.1.16), at $p$,
$$h\bigl([V_1],[V_2]\bigr) = g(V_1,V_2) = -2\zeta\tilde\zeta.\eqno(6.3.12)$$
From (6.1.32), $D(\zeta\tilde\zeta) = 2\rho(\zeta\tilde\zeta)$. So,
$$h\bigl([V_1(v)],[V_2(v)]\bigr) = \left[\exp\int2\rho\right]h\bigl([V_1(p)],[V_2(p)]\bigr).$$
Hence, when $\rho$ also vanishes, $h$ induces a metric on $T_bB$ (with $(T_bB,h) \cong {\bf R}^{1,1}$), but when $\rho$ does not vanish, consistent with the conformal Killing behaviour of $\ell^a$ one only obtains a conformal structure on $T_bB$, i.e., only $h\bigl([V_1],[V_2]\bigr) = 0$ is well defined. Evaluating $h\bigl([V_1],[V_2]\bigr)$ on a smooth section $s:B \to {\cal S}$ yields a metric representative of the conformal structure.
\vskip 24pt
\noindent {\bf 6.3.13 Example}\hfil\break
Let $(M,g)$ consist of a single Walker coordinate system $(u,v,x,y)$ as in (6.1.33). Take ${\cal D} = \langle \partial_u \rangle_{\bf R}$ as before and ${\cal H} = {\cal D}^\perp = \langle \ell^a,m^a, \tilde m^a\rangle_{\bf R}$. The Walker spin frames are a suitable choice to describe the geometry of $\cal H$. As $\kappa = \tilde\kappa = \rho = \tilde\rho = 0$, $\cal H$ is integrable. In fact, $\ell_a = dx$, so the integral hypersurfaces are the level hypersurfaces of the function $x$, i.e., $(u,v,y,x)$ are Frobenius coordinates for the nested integrable distributions ${\cal D} \leq {\cal H}$, with $u$ here an affine parameter for the null geodesic generators. As we know from (6.1.33), $\cal D$ is indeed auto-parallel. The metric form in (6.3.4) is just the Walker canonical form with the coordinates rearranged as $(u,v,y,x)$. This example confirms (6.3.5). The vanishing of $\Sigma$ for $\cal D$ is consistent with (6.1.44). As $\sigma \not= 0$, $\cal H$ is not auto-parallel; $\tilde\sigma = 0$ of course, corresponding to the integrability of the distribution $Z_{[\pi]}$, so one does have the nested integrable distributions ${\cal D} \leq Z_{[\pi]} \leq {\cal H}$ and the hypersurfaces $x = {\rm constant}$ are foliated by the $\alpha$-surfaces of the Walker geometry and have the null geodesics of $\partial_u$ as null geodesic generators. As $\sigma \not= 0$, the distribution $W_{[\alpha]}$ is not integrable. The null distribution $\langle \partial_v \rangle^\perp_{\bf R}$ of type III has similar properties. Its integral hypersurfaces are given by $y = {\rm constant}$. When $a=b=c=0$ one obtains a description of ${\bf R}^{2,2}$. The hypersurfaces $x = {\rm constant}$ are null hyperplanes, foliated both by $\alpha$-planes (spanned by $\ell^a$ and $\tilde m^a$)and $\beta$-planes (spanned by $\ell^a$ and $m^a$).
\vskip 24pt
Returning to a topic in \S 6.2, let $Z_{[\pi]}$ be an integrable $\alpha$-distribution in $(M,g)$, with $S_b \not= 0$ and consider the nested distributions ${\cal D} := \langle S^b \rangle_{\bf R} \leq Z_{[\pi]} \leq {\cal H} := {\cal D}^\perp$. By (6.2.37), $\cal D$ is auto-parallel iff $\Phi_{ABA'B'}\omega^A\omega^B\pi^{A'}\pi^{B'} = 0$, in particular, if $\Phi_{ABA'B'}\pi^{A'}\pi^{B'} = 0$. In fact, this assumption also suffices for the integrability of $\cal H$, yielding a very natural neutral-geometry formulation of Kerr's lemma [23] (7.3.44).
\vskip 24pt
\noindent {\bf 6.3.14 Lemma}\hfil\break
For an integrable $\alpha$-distribution $Z_{[\pi]}$, with $S^a \not= 0$, ${\cal H} = \langle S^a \rangle^\perp_{\bf R}$ is integrable if $\Phi_{ABA'B'}\pi^{A'}\pi^{B'} = 0$.

Proof. $\tilde\kappa = 0$ by integrability of $Z_{[\pi]}$, and $\kappa = 0$ by (6.2.37); by (6.3.2) it remains to show $\rho=\tilde\rho$. Choosing $o^A = \omega^A$ as in \S 6.2 forces $\tilde\rho=0$, see (6.2.36), and $\rho = \omega_B\triangle\omega^B$. Now $\omega^A\triangle S_a = (\omega^A\triangle\omega_A)\pi_{A'}$, so
$$\eqalign{(\omega^A\triangle\omega_A)\pi_{A'} &= \omega^A\triangle(\pi_{C'}\nabla_a\pi^{C'})\cr
&= \omega^A(\triangle\pi_{C'})(\nabla_a\pi^{C'}) + \omega^A\pi_{C'}\triangle\nabla_a\pi^{C'}\cr
&= \omega^A\iota^B(\pi^{B'}\nabla_b\pi_{C'})(\nabla_a\pi^{C'}) + \omega^A\pi_{C'}\triangle\nabla_a\pi^{C'}\cr
&= \omega^A\iota^B(\eta_B\pi_{C'})(\nabla_a\pi^{C'}) + \omega^A\pi_{C'}\triangle\nabla_a\pi^{C'}\cr
&= \iota^B\eta_B\omega^A\omega_A\pi_{A'} + \omega^A\pi_{C'}\triangle\nabla_a\pi^{C'}\cr
&= \omega^A\pi_{C'}\triangle\nabla_a\pi^{C'}.\cr}$$
Hence, and using (6.2.13)(c),
$$\eqalign{(\omega_A\triangle\omega^A)\pi_{A'} &= -\omega^A\iota^B\pi_{C'}\pi^{B'}\nabla_b\nabla_a\pi^{C'}\cr
&= -\omega^A\iota^B[(\eta_B-\omega_B)\omega_A\pi_{A'} - \Phi_{BAC'D'}\pi^{C'}\pi^{D'}\pi_{A'} - \epsilon_{BA}\tilde\Psi_{A'B'C'D'}\pi^{B'}\pi^{C'}\pi^{D'}]\cr
&= \Phi_{ABC'D'}\iota^A\omega^B\pi^{C'}\pi^{D'} + (\omega^A\iota_A)\tilde\Psi_{A'B'C'D'}\pi^{B'}\pi^{C'}\pi^{D'}\cr
&= 0,\cr}$$
by the assumption and noting (6.2.18).
\vskip 24pt
\noindent {\bf 6.3.15 Remarks}\hfil\break
By virtue of (6.2.18), one does not need to explicitly assume $[\pi^{A'}]$ is a multiple WPS as is usually done in the Kerr lemma. The result is automatic when $(M,g)$ is Einstein. The lemma has been formulated in terms of conditions on $[\pi^{A'}]$; just as the weaker condition $\Phi_{ABA'B'}\omega^A\omega^B\pi^{A'}\pi^{B'} = 0$ suffices for $\cal D$ to be auto-parallel the desired result follows if only $\Phi_{ABC'D'}\omega^B\pi^{C'}\pi^{D'}$, provided one then also assumes that $[\pi^{A'}]$ is a multiple WPS.
\vskip 24pt
Consider then an $(M,g,[\pi^{A'}])$, with $[\pi^{A'}]$ a solution of (6.2.1), satisfying $\Phi_{ABA'B'}\pi^{A'}\pi^{B'} = 0$; in particular $[\pi^{A'}]$ is a multiple WPS. I shall call this geometry a {\sl Ricci-aligned $\alpha$-geometry}. It contains nested integrable distributions
$${\cal D} = \langle S^a \rangle_{\bf R} \leq Z_{[\pi]} \leq {\cal H} = {\cal D}^\perp,\eqno(6.3.16)$$
with $\cal D$ auto-parallel by (6.2.37). One can choose Frobenius coordinates $(u,w,z,f)$ for these nested distributions so that ${\cal H} = \langle \partial_u,\partial_w,\partial_z \rangle_{\bf R}$, $Z_{[\pi]} = \langle \partial_u,\partial_w \rangle_{\bf R}$, and ${\cal D} = \langle \partial_u \rangle_{\bf R}$. Choose spin frames with $o^{A'}$ an LSR of $[\pi^{A'}]$ and $o^A = \omega^A$ as in \S 6.2, see (6.2.35--36). By (6.2.37) and (6.3.14), $\kappa = \rho = 0$, which by (6.1.6) are invariant conditions under the remaining freedom in the choice of spin frames. As in \S 6.2, see the discussion leading to (6.2.38), replace $u$ by an affine parameter $v$ for the integral curves of $\cal D$. The coordinates $(v,w,z,f)$ then yield the metric component form of (6.3.4) with $h=1$, with (6.3.4b) taking the form
$$\vcenter{\openup1\jot \halign{$\hfil#$&&${}#\hfil$&\qquad$\hfil#$\cr
\partial_v &= \ell^a = S^a & \partial_f &=: L\ell^a + Mm^a + \tilde M \tilde m^a + n^a\cr
\partial_w &=: A\ell^a + \tilde B \tilde m^a & \partial_z &=: D\ell^a + Cm^a + \tilde C \tilde m^a,\cr}}\eqno(6.3.17)$$
i.e., $B = 0$ in (6.3.4c). As $\tilde\tau = 1$ (which, recall, is equivalent to $\omega^A = o^A$), by (6.3.9) $\cal H$, equivalently $\cal D$, cannot be parallel, so (6.3.4a) does not satisfy Walker's further specialization of these coordinates described in (6.3.10). Note that $S_a = \ell_a = N_a = \nabla_af$.

Now (6.2.38) is in effect, whence (3.4)(\~ c) entails $\tilde\tau(\epsilon - \tilde\epsilon) = \Phi_{10} + \tilde\Psi_1 = 0$, as each summand vanishes. Since $\tilde\tau = 1$, then $\epsilon = \tilde\epsilon$. But from (6.2.38), it follows that $\epsilon = \tilde\epsilon= 0$, i.e., $\omega^A$ and $\pi^{A'}$ are in fact $D$ parallel. There is, therefore, no obstacle to imposing the A-P condition while retaining all assumptions so far; specifically $\ell^a = S^a = \partial_v$ and $\omega^A = o^A$. Furthermore:
\vskip 24pt
\noindent {\bf 6.3.18 Lemma}\hfil\break
For a Ricci-aligned $\alpha$-geometry, with $o^{A'} = \pi^{A'}$, $o^A = \omega^A$, $S^a =\ell^a = \partial_v$, $v$ an affine parameter, then
$$\tilde\beta = 0 \qquad \alpha = -1 \qquad \tilde\Psi_2 + 2\Lambda = 2\tilde\alpha.$$
All three geometric quantities in (6.1.25--27) vanish, as does $\Sigma$ in (6.1.42).

Proof. From (3.4)(\~ e) one obtains $0 = \tilde\tau^2 - \tilde\tau(\tilde\beta - \alpha) + \Phi_{20}$. By the Ricci-alignment condition and (6.3.5), one therefore deduces $\tilde\beta = 0$. So, by (6.3.5), $1 = \tilde\tau = -\alpha$. (3.4)(\~ f) and (6.3.5) give $\tilde\Psi_2 + 2\Lambda = \tilde\tau(\tilde\alpha - \beta - \tau) = 2\tilde\tau\tilde\alpha = 2\tilde\alpha$.
\vskip 24pt
The Frobenius coordinates $(v,w,z,f)$, together with the metric components of (6.3.4) with (6.3.17), provide one approach to studying Ricci-aligned $\alpha$-geometries. Another approach will be considered elsewhere in the context of characterizing neutral geometries $(M,g)$ with algebraically special SD Weyl curvature.
\vskip 24pt
\noindent {\section Appendix}
\vskip 12pt
For convenience, I quote here some formulae cited on several occasions in the text which are derived exactly as in [22], \S\S 4.6 and 4.9, see also [17], Appendix Two. As already noted, my sign conventions for curvature coincide with those of [22] (see [17], Appendix One, for a full account) except my Ricci curvature is the negative of theirs; I accommodate this one difference by retaining their definitions of $\Phi_{ABA'B'}$ and $\Lambda$ and inserting appropriate negative signs in [22], (4.6.20--23), to obtain:
$$R_{ab} = 2\Phi_{ABA'B'} - 6\Lambda\epsilon_{AB}\epsilon_{A'B'} = 2\Phi_{ab} - 6\Lambda g_{ab},\eqno({\rm A}.1)$$
whence
$$R = -24\Lambda \hskip 1.25in 2\Phi_{ab} = R_{ab} - {R \over 4}g_{ab}.\eqno({\rm A}.2)$$
Thus, $2\Phi_{ab} = E_{ab}$ is the tensor often referred to as the Einstein tensor (at least in a non-GR context).

The only modifications needed to carry over [22], \S\S 4.6 \& 4.9, to the neutral-signature context are the obvious ones of replacing complex conjugation by tilde and taking into account the different Hodge star operator. In particular, one has exactly the same definitions of the curvature spinors and
$$\triangle_{ab} := 2\nabla_{[a}\nabla_{b]} = \epsilon_{A'B'}\,\Square_{AB} + \epsilon_{AB}\,\Square_{A'B'},\eqno({\rm A}.3)$$
where 
$$\Square_{AB} := \nabla_{X'(A}\nabla_{B)}{}^{X'} \hskip 1in \Square_{A'B'} := \nabla_{X(A'}\nabla_{B')}{}^X.\eqno({\rm A}.4)$$
Hence
$$\nabla_{BA'}\nabla^B_{B'} = \Square_{A'B'} + {1 \over 2}\epsilon_{A'B'}\Square \hskip 1in \nabla_{AB'}\nabla^{B'}_B = \Square_{AB} + {1 \over 2}\epsilon_{AB}\Square,\eqno({\rm A}.5)$$
where [22] (2.5.24) has been used to express the skew part.

The {\sl spinor Ricci identities\/} for arbitrary spinors $\kappa_A$ and $\tau_{A'}$ are:
$$\displaylines{\Square_{AB}\kappa_C = \Psi_{ABCE}\kappa^E - \Lambda(\kappa_A\epsilon_{BC} + \epsilon_{AC}\kappa_B) \hskip 1in \Square_{AB}\tau_{C'} = \Phi_{ABC'E'}\tau^{E'}\cr
\hfill\llap({\rm A}.6)\cr
\Square_{A'B'}\tau_{C'} = \tilde\Psi_{A'B'C'E'}\tau^{E'} - \Lambda(\tau_{A'}\epsilon_{B'C'} + \epsilon_{A'C'}\tau_{B'}) \hskip .7in \Square_{A'B'}\kappa_C = \Phi_{A'B'CE}\kappa^E\cr}$$
The algebraic Bianchi equation is equivalent to
$$\displaylines{\nabla^A_{B'}\Psi_{ABCD} = \nabla^{A'}_{(B}\Phi_{CD)A'B'} \hskip 1in \nabla^{A'}_B\tilde\Psi_{A'B'C'D'} = \nabla^A_{(B'}\Phi_{C'D')AB}\cr
\hfill\llap({\rm A}.7)\cr
\nabla^{CA'}\Phi_{CDA'B'} = -3\nabla_{DB'}\Lambda\cr}$$
\vskip 24pt
\noindent {\bf Acknowledgments}
I thank Yasuo Matsushita for comments on this paper and the referee for carefully reading the manuscript.
\vskip 24pt
\noindent {\section References}
\vskip 12pt
\newcount\q \q=0
\def\nref {\global\advance\q by1 \item{[\the\q]}}
\frenchspacing
\baselineskip=12pt
\vskip 1pt
\nref V.I. Arnold, Ordinary Differential Equations, MIT Press, Cambridge, MA, 1973.
\vskip 1pt
\nref C.P. Boyer, J.D. Finley III, J.F. Pleba\'nski, Complex General Relativity, $\cal H$ and ${\cal HH}$ Spaces-A Survey of One Approach, in General Relativity and Gravitation: One Hundred Years After the Birth of Albert Einstein, Vol. 2, A. Held (ed.), Plenum Press, New York, NY, 1980, 241--281.
\vskip 1pt
\nref H.A. Buchdahl, On the compatibility of relativistic wave equations for particles of higher spin in the presence of a gravitational field, Nuovo Cim. 10 (1958) 96--103.
\vskip 1pt
\nref C. Chevalley, The Algebraic Theory of Spinors, Columbia University Press, New York, 1954.
\vskip 1pt
\nref A. Derdzinski, Einstein Metrics in dimension four, in: Handbook of Differential Geometry, Vol. I, F.J.E. Dillen, L.C.A. Verstraelen (eds), Elsevier Science B. V., Amsterdam, 2000, pp. 419--707.
\vskip 1pt
\nref A. Derdzinski, W. Roter, Walker's theorem without coordinates, J. Math. Phys. 47 (2006) 062504 (8 pp.).
\vskip 1pt
\nref M. Dunajski, Anti-self-dual four-manifolds with a parallel real spinor, Proc. R. Soc. Lond. A 458 (2002) 1205--1222.
\vskip 1pt
\nref J.D. Finley III, J.F. Pleba\'nski, The intrinsic spinorial structure of hyperheavens, Journal of Mathematical Physics 17 (1976) 2207--2214.
\vskip 1pt
\nref R. Ghanam, G. Thompson, The holonomy Lie algebras of neutral metrics in dimension four, J. Math. Phys. 42 (2001) 2266--2284.
\vskip 1pt
\nref V. Guillemin, S. Sternberg, An Ultra-Hyperbolic Analogue of the Robinson-Kerr Theorem, Lett. Math. Phys. 12 (1986) 1--6.
\vskip 1pt
\nref G.S. Hall, W. Kay, Curvature structure in general relativity, J. Math. Phys. 29 (1988) 420--427.
\vskip 1pt
\nref L.P. Hughston, L.J. Mason, A generalised Kerr-Robinson theorem, Classical and Quantum Gravity 5 (1988) 275--285.
\vskip 1pt
\nref S. Kobayashi, K. Nomizu, Foundations of Differential Geometry. Vol. II. Interscience, New York NY, 1969.
\vskip 1pt
\nref M. Kossowski, Fold Singularities in Pseudo-Riemannian Geodesic Tubes, Proc. Amer. Math. Soc. 95 (1985) 463--469.
\vskip 1pt
\nref P.R. Law, Neutral Einstein metrics in four dimensions, J. Math. Phys. 32 (1991) 3039--3042.
\vskip 1pt
\nref P.R. Law, Classification of the Weyl curvature spinors of neutral metrics in four dimensions, J. Geo. Phys. 56 (2006) 2093--2108.
\vskip 1pt
\nref P.R. Law, Y. Matsushita, A Spinor Approach to Walker Geometry, Comm. Math. Phys. 282 (2008), 577--623; arXiv:math/0612804v4[math.DG] 7 Apr 2008.
\vskip 1pt
\nref P.R. Law, Y. Matsushita, Algebraically Special, Real Alpha-Geometries, arXiv:0808.2082v1[math.DG] 15 Aug 2008.
\vskip 1pt
\nref C. Lebrun, L.J. Mason, Nonlinear Gravitions, Null Geodesics, and Holomorphic Disks. Duke Math. J. 136 (2007) 205--273.
\vskip 1pt
\nref Y. Matsushita, P.R. Law, Hitchin-Thorpe-Type Inequalities for Pseudo-Riemannian 4-Manifolds of Metric Signature $(++--)$, Geom. Ded. 87 (2001) 65--89.
\vskip 1pt
\nref B. O'Neill, Semi-Riemannian Geometry; With Applications to Relativity, Academic Press, Orlando, FL, 1983.
\vskip 1pt
\nref R. Penrose, W. Rindler, Spinors and Space-Time, Vol. 1: Two-spinor calculus and relativistic fields, Cambridge University Press, Cambridge, 1984.
\vskip 1pt
\nref R. Penrose, W. Rindler, Spinors and Space-Time, Vol. 2: Spinor and twistor methods in space-time geometry, Cambridge University Press, Cambridge, 1986.
\vskip 1pt
\nref A.Z. Petrov, Einstein Spaces, Pergamon Press, Oxford, 1969.
\vskip 1pt
\nref J.F. Pleba\~nski, Some Solutions of Complex Einstein Equations, J. Math. Phys. 16 (1975) 2395--2402.
\vskip 1pt
\nref J.F. Pleba\~nski, S. Hacyan, Null geodesic surfaces and Goldberg-Sachs theorem in complex Riemannian spaces, J. Math. Phys. 16 (1975), 2403--2407.
\vskip 1pt
\nref J.F. Plebanski, I. Robinson, Left-Degenerate Vacuum Metrics, Physical Review Letters 37 (1976) 493--495.
\vskip 1pt
\nref J.F. Pleba\'nski, I. Robinson, The Complex Vacuum Metric with Minimally Degenerated Conformal Curvature: in Asymptotic Structure of Space-Time, F. P. Esposito \& L. Witten (eds), Plenum Press, New York \& London, 1977, pp. 361--406.
\vskip 1pt
\nref I.R. Porteous, Topological Geometry, Second Edition, Cambridge University Press, Cambridge, 1969.
\vskip 1pt
\nref M. Spivak, A Comprehensive Introduction to Differential Geometry, Vol. I, Second Edition, Publish or Perish, Berkeley, CA, 1979.
\vskip 1pt
\nref A.G. Walker, On parallel fields of partially null vector spaces, Quart. Journ. of Math. (Oxford) 20 (1949) 135--145.
\vskip 1pt
\nref A.G. Walker, Canonical form for a Riemannian space with a parallel field of null planes, Quart. J. Math. Oxford(2) 1 (1950) 69--79.
\vskip 1pt
\nref A.G. Walker, Canonical forms (II): Parallel partially null planes, Quart. J. Math. Oxford(2) 1 (1950) 147--152.
\vskip 1pt
\bye